\documentclass[10pt]{amsart}
\usepackage[T1]{fontenc}
\usepackage{txfonts}
\usepackage{eucal}

\usepackage{geometry}
\usepackage{fancyhdr}
\usepackage[parfill]{parskip}
\usepackage{lscape}
\usepackage{rotating}

\usepackage{amsmath, amsthm}
\usepackage{amssymb, amsfonts, amscd}
\usepackage{mathrsfs}
\usepackage{xfrac}

\numberwithin{equation}{section}

\usepackage{array}
\usepackage{tabularx}
\usepackage{multirow}
\usepackage{multicol}
\usepackage{booktabs}

\usepackage{xcolor}
\usepackage{float}
\usepackage{caption}
\usepackage{subfig}
\usepackage{mdframed}
\usepackage{mwe}

\usepackage{graphicx}
\usepackage{tikz}
\usepackage{pgf}
\usepackage{tikz-cd}
\usepackage{tikzscale}
\usepackage[all]{xy}

\usetikzlibrary{arrows.meta, positioning}
\usetikzlibrary{
    backgrounds,
    shapes,
    shapes.geometric,
    arrows,
    decorations.markings,
    intersections,
    patterns,
    calc
}

\usepackage{enumitem}
\usepackage{xspace}
\usepackage{comment}
\usepackage{blindtext}
\usepackage[ruled, vlined]{algorithm2e}

\definecolor{smoked}{RGB}{216, 212, 204}
\definecolor{mauve}{RGB}{200, 55, 171}
\definecolor{apricot}{RGB}{250, 144, 4}
\definecolor{sky}{RGB}{66, 169, 244}
\definecolor{plum}{RGB}{76, 0, 102}
\definecolor{forest}{RGB}{90, 145, 120}
\definecolor{sand}{RGB}{180, 160, 120}
\definecolor{sabbia}{RGB}{207, 185, 151}
\definecolor{aliceblue}{rgb}{0.9, 0.95, 1.0}

\usepackage{hyperref}
\hypersetup{
    colorlinks=true, 
    linktoc=all,     
    linkcolor=blue,  
    citecolor=blue,
    filecolor=blue,
    urlcolor=blue
}

\usepackage{cleveref}

\newcommand\C{{\mathbb C}}
\newcommand\DD{{\mathbb D}}

\newcommand\PP{{\mathbb P}}

\newcommand\R{{\mathbb R}}

\newcommand\Z{{\mathbb Z}}
\newcommand\RP{{\mathbb{RP}}}

\newcommand\rp[1]{{\R\PP^{#1}}}

\newcommand{\pslr}[1]{{\mathrm{PSL}(#1,\R)}}
\newcommand{\slr}[1]{{\mathrm{SL}(#1,\R)}}
\newcommand{\slc}[1]{{\mathrm{SL}(#1,\C)}}
\newcommand{\pglr}[1]{{\mathrm{PGL}(#1,\R)}}

\DeclareMathOperator{\hol}{hol}

\newcommand*{\defeq}{\mathrel{\vcenter{\baselineskip0.5ex \lineskiplimit0pt
                     \hbox{\scriptsize.}\hbox{\scriptsize.}}}%
                     =}

\theoremstyle{plain}                   
\newtheorem{thm}{Theorem}[section]
\newtheorem{thma}{Theorem}

\newtheorem{cora}[thma]{Corollary}

\newtheorem{lem}[thm]{Lemma}
\newtheorem{prop}[thm]{Proposition}
\newtheorem{cor}[thm]{Corollary}

\theoremstyle{definition}
\newtheorem{defn}[thm]{Definition}

\newtheorem{rmk}[thm]{Remark}

\theoremstyle{remark}
\newtheorem{ex}[]{Example}

\newtheorem{ques}[thm]{Question}

\newtheorem{claim}[thm]{Claim}  

\newenvironment{proofclaim}
  {\proof}
  {\endproof}

\tikzstyle{rb} = [rectangle, rounded corners, minimum width=3cm, minimum height=1cm, text width=3cm, text centered, draw=black, fill=blue!30]
\tikzstyle{sb} = [rectangle, minimum width=2cm, minimum height=1cm, text width=3cm, text centered, draw=black, fill=violet!30]

\tikzstyle{co} = [circle, radius=2cm, text width=3cm, text centered, draw=white, fill=orange!50]
\tikzstyle{cr} = [circle, radius=2cm, text width=3cm, text centered, draw=white, fill=red!40]
\tikzstyle{cg} = [circle, radius=2cm, text width=3cm, text centered, draw=white, fill=LimeGreen!40]
\tikzstyle{cb} = [circle, radius=2cm, text width=3cm, text centered, draw=white, fill=blue!30]
\tikzstyle{cv} = [circle, radius=2cm, text width=3cm, text centered, draw=white, fill=violet!30]
\tikzstyle{cy} = [circle, radius=2cm, text width=3cm, text centered, draw=white, fill=yellow!50]
\tikzstyle{cgg} = [circle, radius=2cm, text width=3cm, text centered, draw=white, fill=green!40]
\tikzstyle{cc} = [circle, radius=2cm, text width=3cm, text centered, draw=white, fill=cyan!30]

\title[Branched real projective structures on surfaces,\\ cyclic Higgs Bundles and Anosov representations]{Branched real projective structures on surfaces,\\ cyclic Higgs Bundles and Anosov representations}

\author{Gianluca Faraco}
\address[Gianluca Faraco]{School of Mathematics, Monash University, Clayton, VIC 3800, Australia}
\email{gianluca.faraco@monash.edu}

\author{Nicholas Rungi}
\address[Nicholas Rungi]{Instituto de Ciencias Matemáticas (CSIC-UAM-UC3M-UCM), Nicolás Cabrera 13--15, Cantoblanco, 28049 Madrid, Spain.}
\email{nicholas.rungi@icmat.es}

\date{\today}

\begin{document}

\subjclass[2020]{}%
\dedicatory{}

\begin{abstract}
We use the theory of \(3\)-cyclic \(\mathrm{SL}(3,\mathbb R)\)-Higgs bundles to study real branched projective structures on compact surfaces recently introduced by the authors. After fixing a Riemann surface structure on $S$, the Higgs bundle approach allows us, on the one hand, to construct branched structures with prescribed branching data and, on the other hand, to introduce new classes of such structures. The latter are directly related to singular hyperbolic affine spheres and carry a Blaschke metric with conical singularities. In the final part of the paper, using different methods, we prove that every Anosov surface group representation into $\mathrm{SL}(3,\mathbb{R})$ is either Hitchin or belongs to the Barbot component. As a by-product, we obtain that some of the representations associated with the branched structures constructed above are not Anosov.

\end{abstract}

\maketitle
\tableofcontents

\section{Introduction}\label{sec:intro}

In the present paper, we continue the study of branched real projective structures on closed oriented surfaces of genus $g\geq 2$, recently introduced by the authors as geometric structures locally modelled on the real projective geometry \((\slr3,\rp2)\) with branch points. The approach developed here is substantially different from the topological one adopted in \cite{faracorungi_geom}, where the main goal was to establish the geometrisation of representations into $\mathrm{SL}(3,\mathbb R)$. By fixing a Riemann surface structure $X$ on the underlying surface, we consider the geometrisation problem by using the theory of $\mathrm{SL}(3,\mathbb R)$-Higgs bundles, and in particular of a distinguished family of $3$-cyclic Higgs bundles. Although this requires the additional choice of a complex structure, the holomorphic nature of the construction provides more control over the branching data. At the same time, it naturally leads us to introduce a richer class of geometric structures, which we call \emph{branched projective structures of hyperbolic affine sphere type}. One of the first consequences of this approach is that the total degree of branching can be prescribed in a fairly large range. More precisely, we obtain the following. 
\begin{thma}\label{thm: A introduction}
    Let $S$ be a closed and oriented surface of genus $g\geq 2$ and let $X$ be a Riemann surface structure on $S$. For every integer $n$ satisfying $1\le n\le 3g-3$ and $n\neq 2g-2$, there is a family of branched projective structures of hyperbolic affine sphere type on $S$ having total branching order equal to $n$. Moreover the branching divisor can be prescribed arbitrarily provided $0<n<2g-2$.
\end{thma}

Notice that the last statement provides a more refined version of the existence result established in \cite{faracorungi_geom}. In the present note, given any effective divisor $D$ of degree strictly between $0$ and $2g-2$, our construction produces a family of branched projective structures of hyperbolic affine sphere type whose branching divisor is precisely $D$. In fact, these families admit an explicit parametrisation in terms of holomorphic data associated with the corresponding cyclic Higgs bundles. As we shall see later, the additional metric appearing in the construction is closely related to the construction of hyperbolic affine spheres with conical singularities, by providing a singular counterpart of the classical picture associated with Hitchin representations into $\mathrm{SL}(3,\R)$. Notably, the holonomies of the structures constructed above lie in the other connected components of the $\mathrm{SL}(3,\mathbb R)$-character variety. This naturally raises the question of whether these representations can nevertheless be Anosov, a dynamical property satisfied by all Hitchin representations. Using different techniques, we prove the following general result.
\begin{thma}\label{thm: Anosov Barbot Introduction}
    Every Anosov representation $\rho\colon\pi_1(S)\to\mathrm{SL}(3,\mathbb R)$ is either Hitchin or belongs to the Barbot component.
\end{thma}
Here, the Barbot component is the connected component containing the representations obtained by composing a Fuchsian representation $\rho_o\colon\pi_1(S)\to\mathrm{SL}(2,\mathbb R)$ with the reducible embedding $\mathrm{SL}(2,\mathbb R)\to\mathrm{SL}(3,\mathbb R)$. As a consequence, many, though not all, of the representations arising from our construction are necessarily non-Anosov. The question therefore remains open in the remaining cases. 

\medskip

Only after posting this paper on the arXiv did we become aware, through Kassel and Nolte, that the same result was obtained independently by two other teams. In particular, we were informed of independent work by Kassel--Nolte in \cite{KasselNolteTopologicalAspects} devoted to the study of topological invariants of Anosov representations in a substantially more general setting than ours. As a consequence of their main results, they also obtain an independent proof of Theorem \ref{thm: Anosov Barbot Introduction}, which is also based on Barbot's characterisation, but adopts a different approach. At the same time, Brum--Carvajales--Potrie \cite{BrumCarvajalesPotrie} independently proved Theorem \ref{thm: Anosov Barbot Introduction} by using an approach similar to ours but yet different.



\medskip

\subsection{Holomorphic data and branched geometric structures}\label{ssec:higgs bundles}
In our earlier work, see \cite{faracorungi_geom}, we adopted a topological approach to the geometrisation problem where branched projective structures are realised by gluing a finite collection of simplexes modelled on the real projective plane. In the present work we adopt a complex-differential approach based on the theory of Higgs bundles. Upon fixing a Riemann surface \(X\) of genus \(g\geq 2\), the now classical non-abelian Hodge correspondence identifies the moduli space of poly-stable \(G\)-Higgs bundles with the character variety of reductive representations of the surface group into a complex or real semi-simple Lie group \(G\). In a different and more recent direction, the holomorphic theory of Higgs bundles has been used to construct geometric \((G,\mathbb X)\)-structures directly on the underlying surface. In his doctoral dissertation \cite{baraglia2009g2}, Baraglia provided a pioneering contribution by showing, among other results, how a properly convex \(\mathbb{RP}^2\)-structure can be recovered from a stable \(\slr{3}\)-Higgs bundle lying in the image of the Hitchin section. Baraglia also developed an analogous construction for \(\mathrm{PSp}(4,\mathbb R)\), leading to the convex foliated projective structures studied by Guichard--Wienhard; see \cite{guichard08}. Related questions have subsequently been investigated for several other real Lie groups, including \(\mathrm{SO}(2,n)\) \cite{CTT}, \(\mathrm{SL}(3,\mathbb R)\) \cite{rungi2025complex}, in connection with other geometric structures as in \cite{nolte2025concave}, \(\mathrm{G}_2'\) \cite{evans2025geometric,davalo2025geometric}, and \(\mathrm{SO}(p,p+1)\) \cite{davalo2025einstein}.

\smallskip

\subsubsection{Geometrisation via 3-cyclic Higgs bundles}
 
In the framework of branched \((G,\mathbb X)\)-structures, much less is currently known. An early use of Higgs bundles in this setting appears in \cite{biswas2021uniformization}, where they are applied to hyperbolic metrics with conical singularities, originally studied in \cite{troyanov1991prescribing,mcowen1993prescribed}. More recently, \(\mathrm{PSL}(2,\mathbb R)\)-Higgs bundles have been used to study branched harmonic maps and domination problems, leading to constructions of branched anti-de Sitter manifolds; see \cite{diaf2025domination}. As we shall explain, our approach is closer in spirit to the former work, although some analogies with the latter remain. In fact, some of the branched structures we construct from \(\mathrm{SL}(3,\mathbb R)\)-Higgs bundles can be viewed as generalisations of those introduced in \cite{biswas2021uniformization}.

Let \(\mathcal{M}_3^s(X)\) denote the moduli space of stable \(\mathrm{SL}(3,\mathbb R)\)-Higgs bundles. We shall identify subspaces \(\mathcal{CS}_d^+(X)\) and \(\mathcal{CS}_d^-(X)\) of \(\mathcal{M}_3^s(X)\), indexed by integers \(d\) satisfying \(0<d\leq 2g-2\) and \(1-g\leq d<0\), respectively; see Equation \eqref{eq: subspace of cyclic higgs bundles}. As we shall see, these subspaces parametrise isomorphism classes of stable \(3\)-cyclic Higgs bundles; see \S\ref{ssec: cyclic higgs bundles}. In Propositions \ref{prop: parametrization of subspace higgs bundles} and \ref{prop: parametrization of subspace higgs bundles2}, we show that these spaces carry natural complex manifold structures and determine their topology. By using the non-abelian Hodge correspondence, we then study their respective images \(\mathcal R_d^+(X)\) and \(\mathcal R_d^-(X)\) in the character variety \(\mathfrak{X}_3(S)\). We prove that \(\mathcal R_d^+(X)\) is a properly embedded submanifold of \(\mathfrak{X}_3(S)\), while \(\mathcal R_d^-(X)\) is a properly embedded submanifold of its irreducible locus; see Theorems \ref{thm: proper_embedding} and \ref{thm: proper embedding d negative}. These submanifolds have real dimension \(6g-6+2d\). To the best of our knowledge, these proper embedding results have appeared explicitly in the literature only for \(d=2g-2\) and \(d=1-g\); see \cite{Hitchin1992} and \cite{bronstein2025anosov}, respectively. For \(d=2g-2\), the space \(\mathcal R_d^+(X)\) lies in the Hitchin component and is diffeomorphic to \(\mathcal{CS}_d^+(X)\), which is in turn naturally identified with the space of holomorphic cubic differentials on \(X\). For all other values of \(d\), the connected component containing the corresponding representation space depends on the parity of \(d\); see Lemma \ref{lem: rho_connected_component_Higgs_bundles}. When \(d=1-g\), this dependence can equivalently be expressed in terms of the parity of \(g\). 

\smallskip

The next result gives a geometric interpretation of \(\mathcal R_d^\pm(X)\) in terms of branched projective structures equipped with an additional metric datum, as in Definition \ref{def: branched conical projective}.  A branched real projective structure \(\sigma\) on \(S\) with branched points at \(\mathcal B=\{\,p_1,\dots,p_n\,\}\) yields a developing--holonomy pair \((\mathrm{dev},\rho)\), where \(\mathrm{dev}\colon\widetilde{S}\longrightarrow\rp2\) is a smooth map equivariant with respect to a representation \(\rho\colon\pi_1(S)\longrightarrow\pslr3\). 
We say that \(\sigma\) is \emph{of hyperbolic affine sphere type} if the developing map \(\mathrm{dev}\) admits a \(\rho\)-equivariant lift \(f\colon\widetilde{S}\longrightarrow\mathbb{R}^3\setminus\{\,0\,\}\) that is a hyperbolic affine sphere with conical singularities of order \(m_i\) at each lift of \(p_i\), see Definition \ref{def:singular affine sphere}.

\begin{thma}[Theorem \ref{thm: every rho is holonomy of conical branched}]\label{thm: INTRODUCTION - branched conical structures}
    Every representation in \(\mathcal{R}_d^+(X)\) and \(\mathcal{R}_d^-(X)\) arises as the holonomy of a branched projective structure of hyperbolic affine sphere type.
\end{thma}

We shall prove Theorem \ref{thm: INTRODUCTION - branched conical structures} in \S\ref{sec: from cyclic to geometric structres}. As we shall see, our construction produces the desired geometric structure and also allows us to prescribe the orders of its singularities. Before giving further details, we recall for the reader's convenience that an element of \(\mathcal{CS}_d^{\pm}(X)\) is an equivalence class of triples \((L,\gamma,\nu)\), where \(L\) is a holomorphic line bundle of degree \(d\), \(\gamma\in H^0(X,K\otimes L^{-1})\) is non-zero, and \(\nu\in H^0(X,K\otimes L^2)\) is a section which is required to be non-zero in the case of \(\mathcal{CS}_d^-(X)\). The associated real projective structure is a branched real projective structure with branched points on the support of the divisor \(D\defeq\operatorname{div}(\gamma)\) and it comes with a singular Riemannian metric given by \(h^{-1}|\gamma|^2\), where \(h\) is the Hermitian metric on \(L\) that solves the Hitchin equations. This metric is conformal to \(X\) and has singularities precisely at the points of \(\operatorname{supp}(D)\). 

\smallskip

As a direct consequence of Theorem \ref{thm: INTRODUCTION - branched conical structures}, we can solve the problem of prescribing branching data along a divisor \(D\), provided that \(\deg D\) satisfies a suitable inequality; see Theorems \ref{thm: prescribed problem data} and \ref{thm: prescribed problem data2}. These Theorem together yield Theorem \ref{thm: A introduction} as a direct corollary. The degree constraint plays a role analogous to that of the Gauss--Bonnet condition for branched hyperbolic structures in \cite[Theorem A]{troyanov1991prescribing}.

\subsubsection{Realisation into strata of branched real projective structures}\label{sssec:real into strata} Similarly to other class of geometric structures on a surface, branched real projective structures can be classified up to projective diffeomorphisms isotopic to the identity, thereby yielding a moduli space \(\mathcal{RP}^2(\,S\,)\) of marked branched projective structures. By virtue of the fact that a branched real projective structure is uniquely specified by its developing–holonomy pair, see Proposition \ref{prop:BPSasmaps}, assigning to each structure its corresponding holonomy representation leads to a well-defined map that boils down to the so-called \textit{holonomy map} defined as 
\begin{equation}\label{eq:holonomy-map intro}
    \hol \colon \mathcal{RP}^2(\,S\,) \longrightarrow \mathfrak X_3(\,S\,) = \mathrm{Hom}^+\Big(\,\pi_1(\,S\,), \pslr3\,\Big)\,\Big/\,\pslr3.
\end{equation}

The moduli space of branched real projective structures on a surface \(S\) admits a natural stratification into strata \(\mathcal{RP}^2(\,\mu\,)\), where \(\mu = (m_1, \dots, m_k)\) is a tuple of positive integers termed the \textit{signature}. A stratum comprises equivalence classes of branched projective structures possessing precisely \(k\) branch points of orders \(m_1, \dots, m_k\), where a branch point has order \(m\) if the local degree of any local chart around it is \(m+1\) (so that regular points correspond to order zero). The holonomy map defined in \eqref{eq:holonomy-map intro} restricts to a holonomy map on each stratum. In earlier work \cite{faracorungi_geom}, we established that the holonomy map is surjective, whereas a fairly general necessary condition arises in the subtle problem of realising representations in strata of branched projective structures \cite[Theorem C]{faracorungi_geom}. Although this condition is also sufficient for Hitchin representations, no complete characterisation has been established for non-Hitchin representations, thereby leaving open the problem of determining necessary and sufficient conditions for a representation to be realised in a given stratum \cite[Question 1.3]{faracorungi_geom}. As a consequence of our main result, we provide a partial answer in this direction.

\begin{cora}\label{cor: D introduction}
    Let $S$ be a closed oriented surface of genus $g\geq 2$ and let $X$ be a Riemann surface structure on $S$. Let $\mu=(m_1,\ldots,m_k)$ be a signature satisfying \(0<\deg(\,\mu\,)<2g-2\) and set $d=2g-2-\deg(\,\mu\,)$. Then there exists a family of representations in $\mathcal R_d^+(\,X\,)$ which arise as holonomies of branched real projective structures of hyperbolic affine-sphere type belonging to the stratum $\mathcal{RP}^2(\,\mu\,)$. More precisely, for every effective divisor $D$ on $X$ having signature $\mu$, such a family is parametrised by $H^0(X,K^3(-2D))$.
\end{cora}

The presence of branch points gives geometric structures on surfaces considerable flexibility through local deformations that preserve holonomy while moving branch points or splitting them into points of lower order. Such deformations are known as movements of branch points in the theory of branched complex projective structures, while the splitting operation is commonly called breaking up zeros in the theory of translation surfaces. These operations share a common local description. Recall that a branched chart can be viewed as a branched covering between topological discs, with the target disc contained in the model space. The surgery replaces this chart with another branched covering of the same degree, possibly with more branch points, which agrees with the original near the boundary after a suitable identification of the source discs. This replacement preserves the topology of the surface and, since the geometric structure remains unchanged outside the disc and near its boundary, also preserves the holonomy. Moreover, the Riemann–Hurwitz formula ensures that the total branching order remains unchanged, although the individual branching orders may vary. When an underlying complex structure is available, as in the setting of the present paper, these deformations need not preserve it. In the complex-analytic setting, movements of branch points are instances of Schiffer variations. The same local construction applies to branched real projective structures. In particular, a branch point of order \(m\) can be split into branch points of any prescribed positive orders \(m_1,\ldots,m_r\) satisfying \(m_1+\cdots+m_r=m\). Consequently, if a representation \(\rho\) occurs as the holonomy of a branched real projective structure in the stratum \(\mathcal{RP}^2(\,\mu\,)\), then it also occurs as the holonomy of a structure in \(\mathcal{RP}^2(\,\mu'\,)\) for every refinement \(\mu'\) of \(\mu\), obtained by replacing one or more entries of \(\mu\) with partitions of those entries. In particular, \(\deg(\mu)=\deg(\mu')\).

\smallskip

\subsection{Hyperbolic affine spheres with conical singularities}
Notably, any geometric structure in \(\mathcal{CS}_d^\pm(X)\) comes naturally equipped with a holomorphic cubic differential 
\begin{equation}
    q_3 \defeq \nu \gamma^2 \in H^0\big(X, K^3(-2D)\big), \qquad D= \operatorname{div}(\,\gamma\,)
\end{equation}
whose vanishing order at each zero of \(\gamma\) is at least twice the multiplicity of that point in the divisor \(D\). When the degree is maximal, \textit{i.e.}, \(d = 2g-2\), we recover the unbranched setting. In this case, \(\gamma\) is nowhere vanishing (so we may take \(\gamma = 1\in\mathbb{C}\)) and the cubic differential admits a precise geometric interpretation: up to a constant factor, its real part coincides with the Pick tensor of the \(\rho\)-equivariant hyperbolic affine sphere \(f\colon\widetilde S\to\mathbb{R}^3\) obtained by lifting the developing map of the projective structure \cite{nomizu1994affine}, which is properly convex in this setting; see \cite{loftin2001, labourie2007flat, benoist2013cubic}. Moreover, the Hermitian metric \(h^{-1}\) coincides, again up to a constant factor, with the Blaschke metric of \(f\). In particular, \(h^{-1}\) has strictly negative Gaussian curvature, and the norm of the cubic differential is bounded from above. These are by now classical facts that follow from Wang's equation, see,\textit{ e.g.}, \cite[\S 6.2]{li2019introduction}.

\smallskip

Extending this description to the branched setting is far from immediate, and it is natural to ask which properties of the unbranched case still hold. As already mentioned, in \S\ref{sec: affine sphere singular}, we introduce the notion of a hyperbolic affine sphere with conical singularities along a discrete subset \(\mathcal P\); see Definition \ref{def:singular affine sphere}. Roughly speaking, away from \(\mathcal P\) the map is a hyperbolic affine sphere in the classical sense, and we require the Blaschke metric, which is defined on the punctured surface by assumption, to extend to a Riemannian metric on the whole surface with integer conical singularities along \(\mathcal P\). An important feature of this construction is that the holomorphic cubic differential, naturally defined on the punctured Riemann surface, extends holomorphically across \(\mathcal P\), with vanishing order at each point of \(\mathcal P\) at least twice the corresponding order of the singular Blaschke metric, see Lemma \ref{lem: singular hyperbolic affine spheres}. The same type of condition is satisfied by the \(3\)-cyclic Higgs bundles belonging to the families \(\mathcal{CS}_d^+(X)\) and \(\mathcal{CS}_d^-(X)\) considered in the previous section. In fact, we obtain the following general result.
\begin{thma}[Theorem \ref{thm: developing map lifts to a singular affine spheres}]\label{thm: INTRODUCTION - singular affine sphere}
    The developing map of the branched projective structure of hyperbolic affine sphere type associated with a representation \(\rho\) in \(\mathcal{R}_d^+(X)\) or \(\mathcal{R}_d^-(X)\) lifts to a \(\rho\)-equivariant singular hyperbolic affine sphere \(f\colon\widetilde X\longrightarrow\mathbb{R}^3\), whose Blaschke metric \(G_B\) and Pick differential \(q_3\) satisfy:
    \begin{equation}
        G_B=2h^{-1}|\,\gamma\,|^2,\qquad q_3=-2\nu\gamma^2.
    \end{equation}
\end{thma}

By analysing the Hitchin equation \eqref{eq: global Hitchin equation}, we shall study the curvature of the singular Blaschke metric and the norm of the cubic differential. For \(0 < d < 2g-2\), and under suitable assumptions, this allows us to recover results analogous to those in the maximal-degree case. By contrast, when \(1-g \le d < 0\), a new and interesting phenomenon emerges. In what follows, \(K_{G_B}\) denotes the Gaussian curvature of the Blaschke metric of the conical hyperbolic affine sphere given by Theorem \ref{thm: INTRODUCTION - singular affine sphere}, and \(\|q_3\|^2_{G_B}\) denotes the squared norm of its cubic differential with respect to \(G_B\).

\begin{thma}[Corollary \ref{cor: Curvature Blaschke metric singular case}]
    If \(d\in(0,2g-2)\) and \(\operatorname{div}(\gamma)<\operatorname{div}(\nu)\), then \(K_{G_B}\) extends to a strictly negative continuous function on \(X\). Moreover, \(\|q_3\|^2_{G_B}<\frac{1}{2}\) on \(X\).
\end{thma}

Like in the classical case, the idea behind these results is to apply a maximum principle to a function depending on the norms of \(\gamma\) and \(\nu\). For such a principle to be applied, it is sufficient that the vanishing order of \(\nu\) controls that of \(\gamma\) locally, in the sense that \(\operatorname{div}(\nu) > \operatorname{div}(\gamma)\), a condition which can be imposed in the range \(0 < d < 2g-2\). By contrast, the case \(1-g \le d < 0\) is particularly interesting because a conclusion analogous to the one stated above no longer holds. In this range, the natural condition to impose on the holomorphic sections is \(\operatorname{div}(\nu) < \operatorname{div}(\gamma)\). In particular, in the limiting case \(d = 1-g\), this condition is automatically satisfied because \(\nu = 1 \in \mathbb{C}\). Let \(X_0 \defeq X \setminus Z(\gamma)\), where \(Z(\gamma)\) denotes the zero set of \(\gamma\). A similar approach shows that the curvature \(K_{G_B}\) is strictly positive on \(X_0\) and tends to \(+\infty\) along the (possibly empty) set \(Z(\gamma) \setminus Z(\nu)\), as well as at every point \(p \in X\) satisfying \(\operatorname{ord}_p(\nu) < \operatorname{ord}_p(\gamma)\). By the assumption on the divisors, at least one such point always exists. Likewise, one proves that \(\|q_3\|^2_{G_B} > \frac{1}{2}\) on \(X_0\), while its norm diverges to \(+\infty\) under the same conditions as the curvature. Notably, even without imposing any condition on the divisors, a direct application of the Gauss--Bonnet formula for conical Riemannian metrics, see \cite[Proposition 1]{troyanov1991prescribing} yields
\begin{equation}
    \int_{X \setminus Z(\gamma)} K_{G_B} \mathrm{d}A_{G_B} = 2\pi\big(\chi(S) + \operatorname{deg}(D)\big) = -2\pi d > 0,
\end{equation}
where \(D = \operatorname{div}(\gamma)\) and \(\operatorname{deg}(D) = 2g-2-d\).

In particular, the branched \(\rp2\)-structures we construct for \(1-g \le d < 0\) cannot be deformed into projectivisations of hyperbolic conical structures through deformations preserving the total branching order.

\subsection{Anosov representations and topological invariants}
Several works in the literature are devoted to the study of Anosov surface group representations into $\mathrm{SL}(3,\mathbb R)$; see, for instance, \cite{barbot2001flag, Labourie2006Anosov, barbot2010three, kim2016boundary, canary2020topological, bronstein2025anosov, yao2026non, riestenberg2026certifying}.
As already mentioned, the holonomies of the branched projective structures constructed above lie outside the Hitchin component. This naturally raises the question of whether such representations can nevertheless be Anosov. Rather than working directly with the dynamical definition of Anosov representations, our approach relies on a topological characterisation due to Barbot. More precisely, in his work \cite{barbot2010three}, Barbot proved that an Anosov representation into \(\mathrm{SL}(3,\mathbb{R})\) is Hitchin if and only if its associated limit curve in \(\rp2\) is null-homotopic. As we shall explain, this characterisation yields a topological obstruction to the existence of Anosov representations in certain connected components of the character variety. The main result of \S\ref{sec: main result Anosov} can be stated as follows.

\begin{thma}\label{thm: Anosov Stiefel Whitney Introduciton}
    Let $S$ be a closed oriented surface of genus $g\geq 2$, and let $\rho:\pi_1(S)\to\mathrm{SL}(3,\R)$ be an Anosov representation which is not Hitchin. Then $$w_2(\rho)\equiv g-1 \ (\text{mod} \ 2),$$ where $w_2(\rho)$ denotes the second Stiefel--Whitney number of the associated flat bundle.
\end{thma}

It readily follows that Theorem \ref{thm: Anosov Barbot Introduction} is a corollary of the formula above. As we shall explain, the proof of Theorem \ref{thm: Anosov Stiefel Whitney Introduciton} combines Barbot's characterisation with a computation in the cohomology of the unit tangent bundle \(T^1S\) due to Guichard and Wienhard \cite{guichard2010}. Let us denote by \(\pi \colon T^1S \to S\) the fibre bundle projection. Starting from the boundary maps associated with an Anosov representation, we construct three real line bundles \(\mathcal{L}_+\), \(\mathcal{L}_0\), and \(\mathcal{L}_-\) over \(T^1S\), whose direct sum is isomorphic to the pull-back, via \(\pi\), of the flat bundle associated with \(\rho\). Notably, from the first Stiefel--Whitney class \(\alpha \defeq w_1(\mathcal{L}_+)\) we can define an invariant \(\varepsilon(\alpha) \in \mathbb{Z}_2\), which detects whether the restriction of \(\mathcal{L}_+\) to a fibre of \(T^1S\) is orientable. The key observation is that Barbot's characterisation allows us to determine this invariant directly from the topology of the limit curve (Lemma \ref{lem: non orientability on fibres}). More precisely, \(\varepsilon(\alpha)\) vanishes for Hitchin representations and is non-zero for non-Hitchin Anosov representations, see Remark \ref{rem:topological-invariant-anosov}. The remaining part of the proof consists of relating this invariant to the second Stiefel--Whitney class of the flat bundle associated with \(\rho\). To this end, we use the results of Guichard and Wienhard \cite{guichard2010} together with standard algebraic topology.

\smallskip

As an application, for values of \(d\) such that \(d \not\equiv g-1 \pmod{2}\) and \(d \neq 2g-2\), we conclude that any representation \(\rho\) in either space \(\mathcal{R}_d^+(X)\) or \(\mathcal{R}_d^-(X)\) is necessarily non-Anosov, see Corollary \ref{cor:non-anosov-cyclic-slices}. In the complementary case, our result does not determine whether these representations are Anosov, leaving the problem open in general. Specifically, for \(1-g < d < 2g-2\) with \(d \equiv g-1 \pmod{2}\), the question remains open.

\begin{ques}
If \(1-g < d < 2g-2\) and \(d \equiv g-1 \pmod{2}\), is every representation \(\rho\) in \(\mathcal{R}_d^+(X)\) and \(\mathcal{R}_d^-(X)\) non-Anosov?
\end{ques}

\subsection{Structure of paper}
We outline the organisation of the present paper. In \S\ref{sec:rpstructures}, we introduce the notion of a branched projective structure and recall the decomposition of the \(\mathrm{SL}(3,\mathbb{R})\)-character variety into connected components. In \S\ref{sec:holomorphic point of view}, we introduce Higgs bundles and study the subspaces \(\mathcal{CS}_d^\pm(X)\) consisting of \(3\)-cyclic Higgs bundles. In particular, \S\ref{ssec:geo&top of slices} is devoted to their topology and geometry. \S\ref{sec:cyclic Higgs bundles and BGS} contains the main geometric constructions of the paper. We introduce branched projective structures satisfying an additional condition on the developing map, and in \S\ref{sec: from cyclic to geometric structres} we construct such structures using Higgs bundles. Furthermore, \S\ref{sec: prescribed data problem} is devoted to the problem of prescribing the branching data. In \S\ref{sec:analysis_and_geometry}, we analyse the Hitchin equation and derive consequences for the curvature of the singular Blaschke metric and the norm of the Pick differential. Finally, \S\ref{sec: Anosov} is entirely devoted to the proof of our results concerning Anosov representations.

\subsection*{Acknowledgements} The idea of using cyclic Higgs bundles to study branched projective structures originated from discussions between the second author (N.R.) and Andrea Tamburelli about two years ago. This project was initiated following a visit by G.F. to the University of Turin in February 2026 to deliver a seminar, during which he met N.R. The first named author is grateful to Andrea Seppi for the invitation and to the Department of Mathematics of the University of Turin for the warm hospitality during his stay. 
The second author warmly thanks Colin Davalo and Andrea Seppi for many helpful discussions, and Andrea Tamburelli for pointing out an inaccuracy in the statement of Corollary \ref{cor: D introduction} in an earlier version of the paper. N.R. was funded by the European Union (ERC, GENERATE, 101124349). Views and opinions expressed are however those of the author(s) only and do not necessarily reflect those of the European Union or the European Research Council Executive Agency. Neither the European Union nor the granting authority can be held responsible for them.  

\subsection*{GenAI declaration}
The second author (N.R.) used a large language model (LLM) to search for references in the literature concerning Anosov representations into $\mathrm{SL}(3,\mathbb{R})$, beyond the already known reference \cite{barbot2010three}. In particular, the model pointed to \cite{guichard2010}, which was previously unknown to N.R. The entire manuscript, as well as all the ideas underlying the proofs, are the work of the authors. No LLM was used for any of these purposes.

\smallskip

\section{Branched real projective structures on surfaces}\label{sec:rpstructures} 

Let \(\rp2\) denote the real projective plane and \(\pslr3\) denote the group of projective transformations. For the readers' convenience, we recall that in our setting \(\pslr3=\slr3\) and that the latter is naturally isomorphic to \(\pglr3\). Real projective structures, or simply \(\rp2\)-structures, on a closed surface is a geometric structure locally modelled on the real projective plane \(\rp2\), with transition maps in the projective group \(\pslr3\). In the present section, we aim to recall the notion of real branched projective structures on surfaces introduce in the preceding paper \cite{faracorungi_geom}.

\subsection{Basic definitions}\label{ssec:basicdefinitions} Let \(S\) be an oriented and connected surface. A \textit{real projective structure} on \(S\) is the datum of a maximal atlas of \(\rp2\)-valued charts that differ by some projective transformation on their overlaps, see \cite[\S2]{goldman1990}. In order to extend the definition to branched structures, it is first necessary to extend the notion of a chart. We begin with the following

\begin{defn}\label{def:branchechart}
    Let \(S\) be a surface. A \textit{branched real projective chart} is a pair, say \((U,\, \varphi)\), where \(U \subseteq S\) is an open subset and \(\varphi \colon U \to \varphi(\,U\,) \subseteq \rp2\) is a finite degree branched covering map. More specifically, for every \(p \in S\), there exist an open neighbourhood \(U \subset S\) of \(p\), an open neighbourhood \(V \subset \rp2\) of \(\varphi(\,p\,)\), and homeomorphisms \(f\colon U \to \DD\) and \(h\colon V \to \DD\), where \(\DD\) is the unit disc in \(\R^2\), such that \(f(p)=h\circ\varphi(p)=0\) and \(h\circ\varphi\circ f^{-1}|_{\mathbb D\setminus \{0\}}\) is a finite degree covering of punctured discs. By using the polar coordinates \((r,\theta)\) on the disc \(\DD\), the map \(h\circ\varphi\circ f^{-1}\) is isotopic to the map
    \begin{equation}
        (r,\theta)\longmapsto \Big(\,r^{m+1},\,(m+1)\theta\,\Big)
    \end{equation}
    for some integer \(m \geq 1\). Moreover, \(p\) is said to be a \textit{branch point of order} \(m\). We shall say that two branched charts, say \(\varphi_1\) and \(\varphi_2\), at \(p\) are equivalent of there exists a projective transformation \(g\in\pslr{3}\) such that the identity \(\varphi_2=g\circ\varphi_1\) holds.
\end{defn}

\begin{defn}\label{def:branchedprojestructure}
    A \textit{branched \(\rp2\)-structure} on \(S\), denoted by \(\sigma\), is the datum of a maximal atlas, say \(\mathcal U\), of (possibly branched) real projective charts. More precisely, a \(\rp2\)-structure consists of the following data: 
    \begin{itemize}
        \item[1.] a discrete collection of points, say \(\mathcal B \subset S\),
        \item[2.] a maximal atlas of \(\rp2\)-valued charts on \(S \setminus \mathcal B\), and
        \item[3.] a maximal collection branched charts centred at points of \(\mathcal B\), 
    \end{itemize}
    such that the transition maps on \(S\) are restrictions of elements of \(\pslr3\), \textit{i.e.}, for each pair of (possibly branched) charts, say \((U_1,\varphi_1),\,(U_2,\varphi_2)\in\mathcal U\), there is \(g\in\pslr3\) such that the identity \(\varphi_2=g\circ\varphi_1\) holds.
\end{defn}

\noindent We observe that the order of a branch point is well-defined and it does not depend on the choice of any branched chart around the branch point. In fact, for every branched point \(p\in(S,\sigma)\), any pair of branched charts around \(p\) are equivalent because they differ by a projective transformation. Since the latter is a homeomorphism of \(\rp2\), it readily follows that the order of the branched point does not depend on the choice of the local branched chart.

\begin{rmk}\label{rmk:branchingchain}
    Every branched real projective structure on a closed surface \(S\) yields a \(0\)-chain in \(C_o(S,\Z)\), the group of \(0\)-chains on \(S\), that records the orders of the branch points. More specifically, let \(\sigma\) be such a structure on \(S\), and let \(\mathcal B=\{\,p_1,\dots,p_n\,\}\) be its set of branch points with respective orders \(m_1,\dots,m_n\). Then, the \textit{branching chain} is defined as 
    \begin{equation} 
        D=\sum_{p_i\in \mathcal B}\, m_i\,p_i.
    \end{equation}
    Recall that the degree map is the group homomorphism \(\deg\colon C_o(S,\Z) \longrightarrow \mathbb{Z}\) that assigns to each such a formal sum its total sum of coefficients. The degree \(\deg(\,D\,)\) of \(D\) is defined as \textit{total branching order}.
\end{rmk}

\smallskip

We next recall the notion of moduli spaces of branched projective structures. For this purpose, we first recall the notion of isomorphism for these structures.

\smallskip

\begin{defn}\label{def:projectivemaps}
Let \(\Omega\subset\rp2\) be an open set. A map \(\phi\colon\Omega\longrightarrow \rp2\) is said to be \textit{locally projective} if for each connected components \(C\subset\Omega\), there is a projectivity \(g\colon\rp2\longrightarrow \rp2\) such that \(\phi|_{C}=g|_{C}\). Let \(\sigma_1\) and \(\sigma_2\) be two branched \(\rp2\)-structures on \(S\) and let \(f \colon (S,\sigma_1) \longrightarrow (S,\sigma_2)\) be a smooth map. We shall say that \(f\) is \textit{projective} if it locally projective in local charts. We say it is a \textit{projective isomorphism} if it is also a homeomorphism. Finally, we shall denote by \(\textnormal{Proj}(\,\sigma\,)\) the group of projective isomorphisms of \(\sigma\) on \(S\).
\end{defn}

\begin{defn}
A \textit{marked branched \(\rp2\)-structure} on \(S\) is a triple \((S, \sigma, f)\), where \(\sigma\) is a surface endowed with a branched \(\rp2\)-structure and \(f \colon S \longrightarrow (S,\sigma)\) is an orientation-preserving diffeomorphism. Two marked structures, say \((S, \sigma_1, f_1)\) and \((S, \sigma_2, f_2)\), are declared to be equivalent if \(f_2 \circ f_1^{-1} \colon (S,\sigma_1) \longrightarrow (S, \sigma_2)\) is isotopic to a projective isomorphism \(h \colon (S, \sigma_1)\longrightarrow (S, \sigma_2)\). In the following, we shall denote by \(\mathcal{RP}^2(\,S\,)\) the set of marked branched \(\rp2\)-structures on \(S\) up to this equivalence relation.
\end{defn}

The space defined above acquires the structure of a moduli space once it is equipped with an appropriate topology. We shall recall the details of this topology in \S\ref{ssec:devhol_pair}.

\subsection{Developing-holonomy pair}\label{ssec:devhol_pair} From the classical theory of geometric structures, every unbranched structure is completely captured by a pair known as the developing-holonomy pair; see \cite{goldman2022geometric} for an excellent treatment in the most general framework. In particular, every unbranched real projective structure determines such a pair. Let \(\widetilde S\) be the universal cover of \(S\) and let \(\pi\colon\widetilde{S}\longrightarrow S\) be the covering projection. For any \(\rp2\)-structure on \(S\), there exists a unique \(\rp2\)-structure on \(\widetilde S\) that makes the covering projection a locally projective map as in Definition \ref{def:projectivemaps}. If \(\sigma\) is any real (unbranched) projective structure on \(S\), then its \textit{developing map} is defined as the unique locally projective map \(\textnormal{dev}_\sigma\colon\widetilde S\longrightarrow \rp2\) up to post-composition by an element of \(\pslr3\). Concretely, a developing map is constructed by applying the principle of analytic continuation to any chart of the structure's atlas, which is guaranteed by the fact that real projective structures are (real) analytic geometric structures. Since two charts defined on the same open set on \(S\) must differ by an element in \(\pslr3\), it readily follows that a developing map is uniquely determined up to post-composition by an element in \(\pslr3\). The developing map satisfies an equivariant property with respect to a representation \(\rho\colon\pi_1(\,S\,)\longrightarrow \pslr3\) called \textit{holonomy representation} of \(\sigma\). Clearly, if \(\textnormal{dev}_\sigma\) and \(g\cdot\textnormal{dev}_\sigma\) are two developing maps for \(\sigma\) and the former is equivariant with respect to a representation \(\rho\), then the latter is equivariant with respect to \(\textnormal{Ad}_g(\,\rho\,)\). We stress that the argument just presented applies to any surface, not necessarily closed nor of finite type. See \cite[\S2.5]{alessandrini2019} for more details. The following characterisation has been established in \cite[\S2.2.2]{faracorungi_geom}.

\begin{prop}\label{prop:BPSasmaps}
A branched projective structure $\sigma$ on $S$ is equivalent to the existence of a real analytic map $\textnormal{dev}_\sigma\colon \widetilde S\longrightarrow \rp2$ which is a local diffeomorphism away from a discrete subset $\widetilde{\mathcal B}$, is locally a branched covering at every point of $\widetilde{\mathcal B}$, and is equivariant with respect to a representation $\rho\colon\pi_1(S)\longrightarrow \pslr3$.
\end{prop}

\smallskip

As a direct consequence of Proposition \ref{prop:BPSasmaps}, the deformation space of branched projective structures can be identified with the space of equivalence classes of developing-holonomy pairs modulo orientation-preserving diffeomorphisms and the action of \(\pslr3\) by conjugation. Consequently, the deformation space of branched projective structures can be framed from an alternative viewpoint as the orbit space of developing-holonomy pairs under the simultaneous action of \(\pslr3\) by post-composition and the mapping class group by pre-composition. More precisely, we quotient the space of such pairs by post-composition with elements of \(\pslr3\) and by pre-composition of the developing maps with orientation-preserving diffeomorphisms of \(S\) that are isotopic to the identity. 

\begin{equation}\label{eq:identification}
    \mathcal{RP}^2(\,S\,) = \left\{\,\,\begin{gathered} 
    \text{branched \(\rp2\)-structures}\\ 
    \text{on \(S\) as maximal}\\
    \text{atlases of charts}
\end{gathered}\,\, \right\} \longleftrightarrow \Big\{\,\big[\,g\cdot\textnormal{dev},\,\textnormal{Ad}_g(\,\rho\,)\,\big]\,\,\big\vert\,\,g\in\pslr3\,\Big\}_{\sim \textnormal{Diff}_o(\,S\,)}\,\,.
\end{equation}
\smallskip

\noindent By equipping the set of pairs \((\textnormal{dev}, \rho)\) with the compact-open topology, the deformation space naturally inherits the corresponding quotient topology.

\smallskip

\subsection{The representation space and holonomy map}\label{sec: representation space} We next aim to recall some generalities about the space of representation in \(\pslr3\) and the associated character variety. Let \(S\) be a closed oriented surface of genus \(g \ge 2\). The representation space \( \mathrm{Hom}\big(\pi_1(\,S\,),\,\pslr3\big) \) is defined as the space of representations of the fundamental group of \(S\) into \(\pslr3\), which admits a natural structure of a real analytic variety. The \textit{character variety} is defined as the quotient
\begin{equation}
    \mathfrak X_3(\,S\,) = \mathrm{Hom}^+\Big(\,\pi_1(\,S\,), \pslr3\Big)\Big/\pslr3.
\end{equation}
where \(\pslr3\) acts by conjugation on representations, and \(\mathrm{Hom}^+(\pi_1(\,S\,), \pslr3)\) denotes the subset consisting of completely reducible representations. In this way, the induced topology on the space is Hausdorff. In his seminal work \cite{Hitchin1992}, Hitchin investigated the connected components of character varieties associated with split real Lie groups. In the specific case of \(\pslr3\), a key outcome of his work is the following

\begin{thm}[Hitchin]\label{thm: connected components character variety}
    The character variety \(\mathfrak{X}_3(\,S\,)\) has three connected components: the one containing the class of the trivial representation, the one consisting of representations whose associated flat \(\rp2\)-bundles have non-zero second Stiefel-Whitney class, and the one consisting of representations connected to those arising from uniformization. Moreover, the third one is contained in the smooth locus of \(\mathfrak{X}_3(\,S\,)\) and it is diffeomorphic to \(\R^{16g-16}\).
\end{thm}

\noindent It is worth noting that the only topological invariant for a representation \(\rho \in \mathfrak X_3(\,S\,)\) is a characteristic class with values in \(\mathbb{Z}_2\), namely the second Stiefel-Whitney number \(\mathrm{w}_2(\,\rho\,)\). This invariant determines that two of the connected components have \(\mathrm{w}_2 = 0\), while only one has \(\mathrm{w}_2 \neq 0\). In the following, we denote by \(\mathcal{C}_1\) the unique component associated with representations \(\rho\) such that \(\mathrm{w}_2(\,\rho\,)\neq 0\). The other two components both satisfy the condition \(\mathrm{w}_2 = 0\). We shall denote by \(\mathcal{C}_0\) the components containing the trivial representation, and by \(\mathrm{Hit}_3(\,S\,)\), the component diffeomorphic to \(\mathbb{R}^{16g-16}\) which is nowadays called the \textit{Hitchin component}. We stress that, since they share the same vanishing invariant, the Stiefel-Whitney number does not allow one to distinguish the component \(\mathrm{Hit}_3(\,S\,)\) from \(\mathcal{C}_0\), see \cite[\S 10]{Hitchin1992} or \cite[Example 1.10]{collier2019sigma}. 

\smallskip

\subsubsection{Hitchin component}\label{sssec:hitchincomponent} A key feature of the Hitchin component \(\mathrm{Hit}_3(\,S\,)\) is that it contains a copy of the Teichm\"uller space \(\mathcal{T}(\,S\,)\). Recall that the Teichm\"uller space \(\mathcal{T}(\,S\,)\) can be identified with one of the two connected components of the space of faithful and discrete representations \(\rho_o\colon\pi_1(\,S\,)\longrightarrow \mathrm{PSL}(2,\mathbb{R})\) called \textit{Fuchsian representations}. Every representation in \(\pslr2\), in particular Fuchsian once, can be promoted to a representation into \(\pslr3\) by post-composing it with \( \jmath\colon\mathrm{PSL}(2,\mathbb{R})\longrightarrow \pslr3\) the unique irreducible representation (up to conjugation). The assignment \(\rho_o \longmapsto \jmath \circ \rho_o\) defines a natural embedding of \(\mathcal{T}(\,S\,)\) as a closed sub-manifold inside \(\mathrm{Hit}_3(\,S\,)\). According to Labourie in \cite{Labourie2006Anosov}, we shall henceforth adopt the subsequent terminology: a representation in the Hitchin component is called \textit{\(3\)-Fuchsian representations}, if it is defined as the post-composition of Fuchsian representations into \(\pslr2\) with the irreducible representation.

\smallskip

\subsubsection{Barbot component}\label{sec: Barbot component}
Unlike the Hitchin component, one can consider the reducible embedding 
\begin{equation*}
    \iota \colon \mathrm{SL}(2,\mathbb{R}) \to \mathrm{PSL}(3,\mathbb{R}), \quad A \longmapsto \begin{pmatrix} A & 0 \\ 0 & 1 \end{pmatrix}.
\end{equation*}
Starting from a lift of a Fuchsian representation \(\tilde{\rho}_o \colon \pi_1(S) \to \mathrm{SL}(2,\mathbb{R})\), one can therefore consider the induced representation \(\rho \defeq \iota \circ \tilde{\rho}_o \colon \pi_1(S) \to \mathrm{PSL}(3,\mathbb{R})\). Since the Euler number of \(\tilde{\rho}_o\) is equal to \(g-1\), the second Stiefel--Whitney number of the induced representation satisfies \(w_2(\rho) \equiv g-1 \pmod{2}\). If we denote by \(\mathcal{B}\) the connected component of \(\mathfrak{X}_3(S)\) containing \(\rho\), then \(\mathcal{B} = \mathcal{C}_0\) when \(g\) is odd, and \(\mathcal{B} = \mathcal{C}_1\) when \(g\) is even. We shall refer to \(\mathcal{B}\) as the \emph{Barbot component}, since it was originally studied by Barbot in \cite{barbot2010three}. In conclusion, this component always contains \(2^{2g}\) copies of Teichm\"uller space, arising from the reducible embedding \(\iota\) of \(\mathrm{SL}(2,\mathbb{R})\) into \(\mathrm{PSL}(3,\mathbb{R})\) and the choice of a lift of a Fuchsian representation.

\smallskip

\subsubsection{Holonomy map}\label{sssec:holomap} As previously discussed in \S\ref{ssec:devhol_pair}, every branched \(\rp2\)-structure determines an equivalence class of developing-holonomy pairs. The holonomy map is thus defined as the mapping that associates each branched structure with its corresponding holonomy representation, that is,
\begin{equation}\label{eq:holonomy-map}
    \hol \colon \mathcal{RP}^2(\,S\,) \longrightarrow \mathfrak X_3(\,S\,) = \mathrm{Hom}^+\Big(\,\pi_1(\,S\,), \pslr3\,\Big)\Big/\pslr3.
\end{equation}

\noindent It can be shown that this map is continuous with respect to the topologies on the domain and target. In what follows, we shall adopt this terminology. We shall say that a representation \(\rho\colon\pi_1(\,S\,)\longrightarrow\pslr3\) is said to be \textit{geometrisable} whenever it lies in the image of the holonomy map, that is, if it arises as the holonomy of some, possibly branched, real projective structure. One of the main result of the companion paper \cite{faracorungi_geom} states the holonomy is surjective. More precisely, it states the following.

\begin{thm}
    Let \(S\) be a closed orientable surface. Then every representation \(\rho\colon\pi_1(\,S\,)\longrightarrow\pslr3\) arises as the holonomy of some branched projective structure on \(S\).
\end{thm}

Having recalled the necessary background, in the following sections we shall use Higgs bundles to study the geometrisation of representations and other properties from a holomorphic perspective.

\medskip

\section{Cyclic Higgs bundles}\label{sec:holomorphic point of view}

The present section, is the starting point to explore the connections between real branched projective structures and Higgs bundles on Riemann surfaces, thus introducing a holomorphic perspective. In fact, Higgs bundles are well-established as a fundamental tool in the study of the character varieties associated with real (or complex) semisimple Lie groups. In our case, we are interested in the specific cases of \(\mathfrak{X}_3(S)\) and the corresponding one for \(\slc3\), denoted by \(\mathfrak{X}_3^\mathbb{C}(S)\).

\subsection{Higgs bundles and non-abelian Hodge correspondence}\label{sec:higgsbundles}

We begin by briefly recalling the non-abelian Hodge correspondence in the specific case of \(\mathrm{SL}(3,\mathbb{C})\) and \(\slr3\). More general results, valid for complex (or real) Lie groups, can be found in the literature, \textit{e.g.}, see \cite{simpson1992} and \cite{garcia2009hitchin}. In what follows, let \(X\) be a Riemann surface structure on \(S\), and denote by \(K \defeq (T^{1,0}X)^*\) its canonical bundle. 

\smallskip

An \(\mathrm{SL}(3,\mathbb{C})\)-Higgs bundle on \(X\) is a pair \((\mathcal{E},\Phi)\), where \(\mathcal{E}\) is a holomorphic vector bundle of rank 3 such that \(\det(\mathcal{E}) \cong \mathcal{O}_X\), and \(\Phi\) is a holomorphic section of \(\mathrm{End}(\mathcal{E})\,\otimes\,K\) with trace-free endomorphism part. In particular, the holomorphic bundle \(\mathcal{E}\) of any such pair \((\mathcal{E},\Phi)\) has degree zero by definition. We shall say that \((\mathcal{E},\Phi)\) is \emph{stable} if for every proper non-trivial holomorphic sub-bundle \(\mathcal{F} \subset \mathcal{E}\) which is also \(\Phi\)-invariant, we have \(\mathrm{deg}(\,\mathcal F\,)<0\). Moreover, we say that the pair \((\mathcal{E},\Phi)\) is \emph{poly-stable} if it is a direct sum of stable Higgs bundles of degree zero. Let \(\mathcal{M}_3^{\mathbb{C}}(\,X\,)^{ps}\) be the moduli space of poly-stable \(\mathrm{SL}(3,\mathbb{C})\)-Higgs bundles, where two pairs, say \((\mathcal{E}_1,\Phi_1)\) and \((\mathcal{E}_2,\Phi_2)\), are declared to be equivalent if there exists an isomorphism of holomorphic vector bundles, say \(g \colon \mathcal{E}_1 \to \mathcal{E}_2\) with \(\det(\,g\,) = 1\) such that the identity \(\Phi_2 = g \cdot \Phi_1 \cdot g^{-1}\) holds. In what follows, we shall refer to any such isomorphism \(g\) as a \emph{gauge transformation}. For the reader's convenience, we recall that \(\mathcal{M}_3^{\mathbb{C}}(\,X\,)^{ps}\) carries a natural Hausdorff topology, see \cite[\S 3.2]{collier2019sigma}. The notions of stability and poly-stability naturally extend to Higgs bundles with structure group \(\mathrm{SL}(3,\mathbb{R})\). These are defined as (poly)-stable \(\mathrm{SL}(3,\mathbb{C})\)-Higgs bundles \((\mathcal{E},\Phi)\) equipped with the additional structure determined by a section \(Q \in H^0(X, \mathrm{Sym}^2(\mathcal{E}))\) that preserves the Higgs field. More specifically, \(Q \colon \mathcal{E}\to\mathcal{E}^*\) is a holomorphic isomorphism such that the identity \(\Phi^T Q = Q \Phi\) holds. In this setting as well, it is possible to define a moduli space of (poly)-stable Higgs bundles, denoted by \(\mathcal{M}_3(\,X\,)^{s}\) (respectively \(\mathcal{M}_3(\,X\,)^{ps}\)).

\smallskip

Going back to the \(\slc3\) setting, by fundamental results of Hitchin \cite{hitchin1987} and Simpson \cite{simpson1988, simpson1992}, the condition of polystability is equivalent to the existence of a Hermitian metric \(H\) on \(\mathcal{E}\) satisfying the Hitchin equations
\begin{equation}\label{eq: Hitchin_higgs bundles}
F_{\nabla_H} + \Big[\,\Phi,\, \Phi^{*H}\,\Big] = 0,
\end{equation}
where $F_{\nabla_H}$ is the curvature of the Chern connection associated with $H$, and $\Phi^{*H}$ is the adjoint of $\Phi$ with respect to the Hermitian metric. Moreover, if the Higgs bundle is stable, such a Hermitian metric is unique. Given the Hermitian metric \(H\), it is possible to define a flat connection \(\nabla\) on the smooth complex vector bundle \(E\) underlying \(\mathcal{E}\) by
\begin{equation}\label{eq: flat connection non abelian Hodge}
\nabla\defeq \nabla_H + \Phi + \Phi^{*H}.
\end{equation}
Indeed, the condition $\nabla^2 = 0$ is equivalent to the Hitchin equation \eqref{eq: Hitchin_higgs bundles}. The flat bundle $(E,\nabla)$ determines a representation, say
\(
\rho\colon\pi_1(\,S\,)\to \slc3,
\)
called \emph{holonomy representation}, which is defined via the parallel transport along closed loops. Conversely, given a representation \(\rho\colon\pi_1(\,S\,)\longrightarrow \slc3\), the associated rank-$3$ complex flat bundle $(E_\rho,\nabla)$ can be realised as follows. Consider the trivial flat bundle
\(
\widetilde{E} \defeq \widetilde{S} \times \mathbb{C}^3
\)
over the universal cover of the surface, equipped with the trivial connection. The induced bundle over $S$ is then defined as
\begin{equation}\label{eq: flat bundle}
E_\rho \defeq \frac{\widetilde{S} \times \mathbb{C}^3}{ \sim_\rho},
\end{equation}
where $(p,v) \sim_\rho (q,w)$ if and only if there exists $\gamma \in \pi_1(\,S\,)$ such that $q = \gamma \cdot p$ and $w = \rho(\,\gamma\,)\cdot v$. The trivial connection on \(\widetilde E\) boils down to a flat connection, say $\nabla$, on $E_\rho$. Seminal results of Donaldson \cite{donaldson1987} and Corlette \cite{corlette1988} show that if the flat connection \(\nabla\) on \(E_\rho\) is irreducible (resp. reductive), then there exists a Hermitian metric $H$ on $E_\rho$ such that the triple $(E_\rho,\nabla_H^{0,1},\Psi_H^{1,0})$ defines a stable (resp. poly-stable) $\mathrm{SL}(3,\C)$-Higgs bundle, where $\nabla_H$ is the Chern connection of the Hermitian metric \(H\) and \(\Psi_H^{1,0}\) is a \((1,0)\)-form with values in the \(H\)-Hermitian endomorphisms of $E_\rho$ that recovers the connection $\nabla$ together with $\nabla_H$. The proof of these results relies on the existence of a \(\rho\)-equivariant harmonic map
\(
f\colon \widetilde{X} \to \slc3/\mathrm{SU}(3).
\)
For this reason, Hermitian metrics $H$ as above are often referred to as \emph{harmonic metrics}.

\begin{thm}[\cite{hitchin1987, simpson1988, corlette1988, donaldson1987, garcia2009hitchin}]\label{thm: non abelian hodge}
The correspondence described above defines a real-analytic diffeomorphism between the space $\mathcal{M}_3^{\mathbb{C}}(\,X\,)^{ps}$ (resp. \(\mathcal{M}_3(\,X\,)^{ps}\)) and $\mathfrak{X}_3^{\mathbb{C}}(\,S\,)$ (resp. \(\mathfrak{X}_3(\,S\,)\)) on their smooth loci. Moreover, the moduli space $\mathcal{M}_3^{\mathbb{C}}(\,X\,)^s$ (resp. \(\mathcal{M}_3(\,X\,)^s\)) of stable Higgs bundles corresponds to the space of irreducible representations in $\mathfrak{X}_3^{\mathbb{C}}(\,S\,)$ (resp. \(\mathfrak{X}_3(\,S\,)\)).
\end{thm}

\subsection{Slice of cyclic Higgs bundles}\label{ssec: cyclic higgs bundles} We now turn our focus to cyclic Higgs bundles, a special family of Higgs bundles that will play a pivotal role in this work. In particular, the results established in this section will serve as essentials building blocks for the proofs of Theorems \ref{thm: proper_embedding} and \ref{thm: proper embedding d negative}. For this purpose, let 
\(X\) be a fixed Riemann surface, and consider the \(\slc3\)-Higgs bundle determined by the data
\begin{equation}\label{eq: higgs data}
\mathcal{E} = L \oplus \mathcal{O}_X \oplus L^{-1}, \qquad
\Phi =
\begin{pmatrix}
0 & 0 & \nu \\
\gamma_1 & 0 & 0 \\
0 & \gamma_2 & 0
\end{pmatrix},
\end{equation}
where $L$ is a line bundle of degree $d \in \mathbb{Z}$, \(\gamma_1,\gamma_2\) are holomorphic sections of \(L^{-1}K\) and \(\nu\) is a holomorphic section of $L^2K$. The following Proposition provides sufficient conditions for an \(\slc3\)-Higgs bundle to be stable. 

\begin{prop}\label{prop: stable SL(3;C) Higgs bundle}
    The $\mathrm{SL}(3,\C)$-Higgs bundle $(\mathcal{E},\Phi)$ described above is stable if one of the following holds:
    \begin{itemize}
        \item[(i)] $d\in\big(\,0,\,2g-2\big]$ and $\gamma_1,\gamma_2\neq 0$;
        \item[(ii)] $d\in\big[1-g,0\,\big)$ and $\gamma_2,\nu\neq 0$ .
    \end{itemize}Moreover, in the case \(d\in[1-g,0)\) and \(\gamma_1=\gamma_2=0\), the \(\slc3\)-Higgs bundle \((\mathcal E,\Phi)\) is strictly poly-stable.
\end{prop}

\begin{proof}
    According to a result of Simpson \cite[Lemma 6.4]{simpson_cyclic}, in our framework it suffices to verify stability along the sub-bundles induced by the direct sum decomposition of \(\mathcal{E}\). First, assume that $d\in\big(\,0,2g-2\big]$. If $\gamma_1,\gamma_2,\nu\neq 0$, then no proper non-trivial sub-bundle arising from the splitting of $\mathcal{E}$ is preserved by $\Phi$. In the case \(\nu=0\), the only $\Phi$-invariant sub-bundles of this form are $L^{-1}$ and $\mathcal{O}_X\oplus L^{-1}$; since their degrees are strictly negative by hypothesis, the stability condition is satisfied. Let us now assume $d\in[1-g,0)$. In this case the sub-bundles $L^{-1}$ and $\mathcal{O}_X\oplus L^{-1}$ have non-negative degree, while \(\mathcal O_X\) has zero degree, stability is therefore guaranteed provided that the conditions $\nu\neq 0$ and \(\gamma_2\neq 0\) hold. The last statement follows from the fact that if \(\gamma_1=\gamma_2=0\), then we can decompose \((\mathcal{E},\Phi)\) as the direct sum of two degree zero Higgs bundles.
\end{proof}

The complex multiplicative group \(\mathbb{C}^{*}\) acts naturally on the moduli space \(\mathcal{M}_3^{\mathbb{C}}(\,X\,)^{ps}\) of polystable \(\slc3\)-Higgs bundles. This action scales the Higgs field linearly while preserving the isomorphism class of the underlying holomorphic vector bundle. Explicitly, for any \(t \in \mathbb{C}^*\), the action is given by:
\begin{equation}
    t \cdot [\mathcal{E},\Phi] \defeq  [\mathcal{E},t\,\Phi].
\end{equation}
An \(n\)-cyclic Higgs bundle is defined as a Higgs bundle whose isomorphism class is fixed under the restriction of the natural \(\mathbb{C}^*\)-action to the subgroup \(\langle\,\zeta_n\,\rangle \cong \mathbb{Z}_n\) of \(n\)-th roots of unity. In our framework, a direct computation shows that
\begin{equation}
    \zeta_3\Phi = g \cdot \Phi \cdot g^{-1}
\end{equation}
for the order 3 gauge matrix \(g = \operatorname{diag}(\,\zeta_3^{-1}, 1, \zeta_3\,) \in \slc3\), thus showing that Higgs bundles as in Proposition \ref{prop: stable SL(3;C) Higgs bundle} are 3-cyclic. Compare with the proof of Lemma \ref{lem: gauge equivalent Higgs bundles}.

\begin{rmk}\label{rem: cyclic_higgs_bundles}
    We may notice that when \(\nu=0\), the Higgs bundle is in fact fixed by the whole \(\mathbb C^*\)-action. Due to the specific symmetry of $\mathcal{E}$ and $\Phi$, up to a suitable gauge choice, the Hermitian metric $H$ solving the Hitchin equations \eqref{eq: Hitchin_higgs bundles} is diagonal, see \cite[Theorem 4.1.9]{collier2016finite}. 
\end{rmk}

In the case \(\deg(\,L\,)=2g-2\), the isomorphism $L \cong K$ holds provided that at least one of $\gamma_1$ and $\gamma_2$ is not identically zero. When instead $\deg(\,L\,)=1-g$, the isomorphism
\(
L \cong K^{-\frac{1}{2}}
\)
holds provided that $\nu \neq 0$, for some square root of the canonical bundle. In particular, this condition is always satisfied under assumptions (i) and (ii) of Proposition \ref{prop: stable SL(3;C) Higgs bundle}, which ensures stability of the Higgs bundle.
Let us now consider the case in which $\gamma\defeq\gamma_1 = \gamma_2$, so that the Higgs field $\Phi$ and the bilinear form
\begin{equation}
    Q =
    \begin{pmatrix}
        0 & 0 & 1 \\
        0 & 1 & 0 \\
        1 & 0 & 0
    \end{pmatrix}
\end{equation}

are compatible, namely $\Phi^TQ=Q\Phi$, and hence the triple $(\mathcal{E},\Phi,Q)$ defines an $\mathrm{SL}(3,\mathbb{R})$-Higgs bundle. From now on, we shall always assume that $\gamma \neq 0$, so that the bundle $(\mathcal{E},\Phi,Q)$ is stable whenever $d \in (0,2g-2]$, by Proposition \ref{prop: stable SL(3;C) Higgs bundle}. In the case \(d = 2g-2\), and hence \(L \cong K\), then \(\nu\) is a holomorphic section of \(K^3\), that means, \(\nu\) is a holomorphic cubic differential on \(X\). Under this identification, the holonomy of the associated flat connection belongs to the Hitchin component \(\mathrm{Hit}_3(\,S\,)\), see \S\ref{sec: representation space}.

Whenever \(0\neq d \in [1-g, 2g-2)\), the holonomy of \((\mathcal{E}, \Phi, Q)\) lies in one of the other two connected components of \(\mathfrak{X}_3(\,S\,)\), namely \(\mathcal{C}_0\) or \(\mathcal{C}_1\). More specifically, by Proposition~\ref{prop: stable SL(3;C) Higgs bundle}, the Higgs bundle \((\mathcal{E}, \Phi)\) is stable; hence, the Hitchin equations~\eqref{eq: Hitchin_higgs bundles} admit a unique solution, defining a Hermitian metric~\(H\) on \(\mathcal{E}\). The compatibility condition between the Higgs field $\Phi$ and $Q$ implies a compatibility condition between the Hermitian metric $H$ and $Q$; \textit{i.e.}, \(Q H^t Q = H^{-1}\), see \cite{collierLi2019}. Moreover, a conjugate-linear involution, say \(\tau \colon \mathcal{E} \longrightarrow \mathcal{E}\), on the underlying complex vector bundle can be defined by setting \(\tau(\,s\,) \defeq H^{-1}Q\bar{s}\) for any section \(s \in \mathcal{E}\). The fixed-point locus of~\(\tau\) defines a real sub-bundle~\(\mathcal{E}_{\mathbb{R}}\) of rank~\(3\), carrying the flat connection induced by \(\nabla = \nabla_H + \Phi + \Phi^{*_H}\). Indeed, a straightforward computation shows that
\(
\nabla\big(\,\tau(\,s\,)\,\big) = \tau\big(\,\nabla s\,\big)
\), for any smooth section $s$ of $E$. 
\begin{lem}\label{lem: rho_connected_component_Higgs_bundles}
Let $(\mathcal E,\Phi,Q)$ be the $\mathrm{SL}(3,\R)$-Higgs bundle determined by the data
\begin{equation}\label{eq: cyclic Higgs bundles}
    \mathcal{E} = L \oplus \mathcal{O}_X \oplus L^{-1}, \quad
    \Phi =
    \begin{pmatrix}
        0      & 0      & \nu \\
        \gamma & 0      & 0   \\
        0      & \gamma & 0
    \end{pmatrix}, 
    \quad Q =
    \begin{pmatrix}
        0 & 0 & 1 \\
        0 & 1 & 0 \\
        1 & 0 & 0
    \end{pmatrix}.
\end{equation} Assume \(\gamma \neq 0\), and that \(\nu \neq 0\) only if \(d \in [1-g, 0)\). Let \(\rho\) be the corresponding holonomy representation. Then the following condition hold: 
\begin{itemize}
    \item[(i)] $\rho\in\mathcal C_0$ if and only if $1-g<d\neq 0<2g-2$ and $d$ is even, or $d=1-g$ and $g$ is odd; 
    \item[(ii)] $\rho\in\mathcal C_1$ if and only if $1-g<d\neq 0<2g-2$ and $d$ is odd, or $d=1-g$ and $g$ is even.
\end{itemize}

\end{lem}
\begin{proof}
    First, by Proposition \ref{prop: stable SL(3;C) Higgs bundle}, the Higgs bundle \((\mathcal{E}, \Phi, Q)\) is stable and therefore determines an associated holonomy representation \(\rho \in \mathfrak{X}_3(\,S\,)\). Moreover, since $d$ is always less than $2g-2$, it follows that $\rho$ does not lie in the Hitchin component, see \cite{Hitchin1992}. It remains to determine the connected component in which the representation lies. By construction of the real structure induced by the involution~\(\tau\), there is an isomorphism of real vector bundles \(E_{\mathbb{R}} \cong L_{\mathbb{R}} \oplus \underline{\mathbb{R}}\), where \(L_{\mathbb{R}}\) denotes the underlying rank-\(2\) real bundle of~\(L\), and \(\underline{\mathbb{R}}\) is the trivial real line bundle, see \cite[\S 10]{Hitchin1992}. The topological invariant of~\(\rho\) is characterised by the second Stiefel--Whitney class of~\(E_{\mathbb{R}}\), yielding
    \begin{equation}
        \mathrm{w}_2(\,E_{\mathbb{R}}\,) = \mathrm{w}_2(\,L_{\mathbb{R}} \oplus \underline{\mathbb{R}}\,) = \mathrm{w}_2(\,L_{\mathbb{R}}\,) \equiv c_1(\,L\,) \pmod{2},
    \end{equation}
    from which the claim follows.
\end{proof}

\begin{rmk}
Higgs bundles with $d=1-g$ and as in Lemma \ref{lem: rho_connected_component_Higgs_bundles} have been recently studied in \cite{bronstein2025anosov}. One of the main results was to show that the associated holonomy representation is Anosov, see \cite[Theorem 6.5]{bronstein2025anosov}, in particular it is faithful and discrete. Such Higgs bundles form a subset of the \emph{Slodowy slice} associated with the reducible representation $\mathfrak{sl}(2,\mathbb{C}) \to \mathfrak{sl}(3,\mathbb{C})$, see \cite[\S 2]{bradlow2024general}. %
\end{rmk}

The Higgs bundles considered in Lemma~\ref{lem: rho_connected_component_Higgs_bundles} are uniquely determined by the line bundle~\(L\) and the sections~\(\gamma\) and \(\nu\). To abridge the notation, in what follows we shall denote them by~\((L,\gamma,\nu)\).

\begin{lem}\label{lem: gauge equivalent Higgs bundles}
    Two stable $\slr3$-Higgs bundles $(L,\gamma,\nu)$ and $(L',\gamma',\nu')$ with $d\in\big[\,1-g,\,2g-2\,\big]\setminus\{\,0\,\}$ are gauge equivalent if and only if $L\cong L'$ and there exists $\lambda\in\C^*$ such that $\gamma'=\lambda\gamma$ and $\nu'=\lambda^{-2}\nu$. 
\end{lem}

\begin{proof}
    According to our discussion in \S\ref{sec:higgsbundles}, we notice that if $\theta$ is a gauge transformation realising the identity $\theta\cdot (L,\gamma,\nu)=(L',\gamma',\nu')$, then $L$ is necessarily isomorphic to $L'$. By regarding \(\theta\) as a \(3 \times 3\) matrix, its block structure is constrained to be upper triangular for \(d \in \big(\,0,2g-2\big]\), and lower triangular for \(d \in \big[1-g,0\,\big)\). In what follows we shall assume to be in the former case, \textit{i.e.} $d\in\big(\,0,2g-2\big]$; since the latter one follows \textit{mutatis mutandis}. Let us write
    
    \begin{equation}
        \theta=
        \begin{pmatrix}
        a_{11} & a_{12}      & a_{13} \\ 
        0 & a_{22} & a_{23} \\ 
        0 & 0 & a_{33}
        \end{pmatrix}.
    \end{equation}
    By definition, $\theta$ conjugates the Higgs field associated with the pair $(\gamma,\nu)$ to the one associated with $(\gamma',\nu')$. This readily implies that $a_{12}=a_{13}=a_{23}=0$. Moreover, since $\theta$ preserves the bilinear form $Q$ we conclude that $\theta=\mathrm{diag}(\lambda,1,\lambda^{-1})$ for some $\lambda\in\C^*$, hence the desired result.
\end{proof}

We study the moduli subspace of stable \(\mathrm{SL}(3,\mathbb{R})\)-Higgs bundles as in Lemma~\ref{lem: rho_connected_component_Higgs_bundles}, showing that it carries the structure of a vector bundle over a suitable symmetric product of the surface~\(X\), and computing its dimension. The argument proceeds in close analogy with those presented in \cite[\S 10]{hitchin1987} and \cite[\S 5.2]{rungi2025complex}. For any fixed integer \(d \in [1-g,2g-2] \setminus \{\,0\,\}\), we define
\begin{equation}
    \begin{aligned}\label{eq: subspace of cyclic higgs bundles}&\mathcal{CS}_d^+(\,X\,)\defeq\left\{\,\big[L,\gamma,\nu\big]\,\ \middle|\ \,\deg(\,L\,)=d\in(0,2g-2],\ \gamma\in H^0(X,L^{-1}K)\setminus\{\,0\,\},\ \nu\in H^0(X,L^2K)\,\,\right\}, \\& \mathcal{CS}_d^-(\,X\,)\defeq\left\{\,\big[L,\gamma,\nu\big]\ \,\middle|\ \,\deg(\,L\,)=d\in[1-g,0),\ \gamma\in H^0(X,L^{-1}K)\setminus\{\,0\,\},\ \nu\in H^0(X,L^2K)\setminus\{\,0\,\}\,\,\right\}.
    \end{aligned}
\end{equation}

We observe that every element of these two spaces corresponds to an isomorphism class of stable \(\mathrm{SL}(3,\mathbb{R})\)-Higgs bundles, by Proposition \ref{prop: stable SL(3;C) Higgs bundle}. 

\smallskip

\subsection{Geometry and topology of slices}\label{ssec:geo&top of slices} 

This section is devoted to the study of the slices introduced in \S\ref{ssec: cyclic higgs bundles}. Before proceeding, it is convenient to recall some standard terminology concerning divisors on Riemann surfaces that will be used in what follows. 

\begin{defn}\label{def: Comparison of divisors}
    Given two divisors on a Riemann surface \(X\), say
    \begin{equation}
    D_1=\sum_{p\in X} m_p p \qquad \text{and} \qquad D_2=\sum_{p\in X} n_p p,
    \end{equation}
    the relation \(D_1 \leq D_2\) holds if \(m_p \leq n_p\) for every \(p \in X\), or equivalently if the difference \(D_2 - D_1\) is an effective divisor. Analogously, the strict inequality \(D_1 < D_2\) is defined if \(D_1 \leq D_2\) and the coefficient inequality is strict for at least one point~\(p \in X\). Finally, the \emph{support} of a divisor~\(D\), denoted by \(\operatorname{Supp}(\,D\,)\), is the finite subset of~\(X\) consisting of all points appearing in~\(D\) with non-zero coefficients.
\end{defn}

Given a non-vanishing holomorphic section~\(s\) of a line bundle~\(L\) over~\(X\), the associated divisor of zeroes, denoted by \(\operatorname{div}(\,s\,)\), is defined as the formal sum
\begin{equation}
\text{div}(\,s\,)\defeq\sum_{i=1}^k \operatorname{ord}_{p_i}(\,s\,)\cdot p_i,
\end{equation}
where $\{p_1,\dots,p_k\}\subset X$ is the set of zeros of $s$ and $\operatorname{ord}_{p_i}(\,s\,)\in \mathbb{Z}_{\geq 1}$ denotes the local vanishing order of the section at the point $p_i$. In other words, any non-zero element $s\in H^0(X,L)$ yields an inclusion $\mathcal{O}_X\hookrightarrow L$ and a (non canonical) isomorphism $L\cong\mathcal O_X(\,D\,)$ with $D\defeq\text{div}(\,s\,)$, see \cite[Proposition 2.3.18]{huybrechts2005complex} for more details. We now obtain a parametrisation for \(\mathcal{CS}_d^+(\,X\,)\) and \(\mathcal{CS}_d^-(\,X\,)\) similar to the one originally introduced by Hitchin for \(\pslr2\), see \cite[Theorem~10.8]{hitchin1987}. The proof is similar to those of \cite[Theorem~5.4]{collier2024holomorphic} and \cite[Proposition~2.7]{diaf2025domination}, see also \cite[Theorem~5.7]{rungi2025complex}.

\begin{prop}\label{prop: parametrization of subspace higgs bundles}
    Let $X$ be a Riemann surface of genus $\ge 2$, then for every $d\in\big(\,0,2g-2\big]$, the slice $\mathcal{CS}_d^+(\,X\,)$ is biholomorphic to the total space of a holomorphic vector bundle of complex rank $2d+g-1$ over the $(2g-2-d)$-th symmetric product of $X$.
\end{prop}

\begin{prop}\label{prop: parametrization of subspace higgs bundles2}
Let $X$ be a Riemann surface of genus $\ge 2$. Then for every $d\in\big[1-g,0\,\big)$, the slice $\mathcal{CS}_d^-(\,X\,)$ is biholomorphic to a bundle over an $H^1(X,\Z_2)$-cover of the $(2d+2g-2)$-th symmetric product of $X$ whose fibre is $(\C^{g-1-d}\setminus\{\,0\,\})/\{\,\pm\mathrm{Id}\,\}$. 
\end{prop}

\begin{proof}[Proof of Proposition \ref{prop: parametrization of subspace higgs bundles}]

Let \(d\) be a positive integer, \(d\le2g-2\), and consider the space
\begin{equation}\label{eq: holomoprhic bundle cubic twisted}
    \mathcal{QD}_d(\,X\,)\defeq\Big\{\,(q_3,D)\mid D\in\operatorname{Sym}^{2g-2-d}(X),\, q_3\in H^0\big(X,K^3(-2D)\big)\Big\}.
\end{equation}

We may observe that $\mathcal{QD}_d(\,X\,)$ is a holomorphic vector bundle over $\mathrm{Sym}^{2g-2-d}(\,X\,)$ with fibre isomorphic to the vector space $H^0(X,K^3(-2D))$ and a straightforward application of the Riemann--Roch theorem shows that the fibre has complex dimension $2d+g-1$. Moreover, whenever we have a pair \((q_3,D)\in\mathcal{QD}_d(X)\) with \(q_3\neq 0\), then \(\operatorname{div}(q_3)\ge2D\). Let 
\begin{equation}\label{eq: holomorphic_map_CS_d}
    \psi_d \colon \mathcal{CS}_d^+(\,X\,) \longrightarrow \mathcal{QD}_d(\,X\,),
\end{equation}
be the map that associates to an isomorphism class \([L, \gamma, \nu]\) the pair \(\big(\,\nu\gamma^2, \operatorname{div}(\,\gamma\,)\,\big)\). 

\begin{claim}
    The mapping \(\psi_d\) is well-defined.
\end{claim}

\begin{proofclaim}
    For every triple \((L, \gamma, \nu)\), the product \(q_3 = \nu \gamma^2\) is a holomorphic section of \(L^2 K \otimes (L^{-1} K)^2 \cong K^3\). Furthermore, since \(\gamma\) is a non-zero holomorphic section and \(\deg(\,L\,)=d\), then \(D = \operatorname{div}(\,\gamma\,)\) is an effective divisor of degree \(\deg(L^{-1}K) = 2g - 2 - d\). Next, we may observe that \(\nu\gamma^2\in H^0\big(X,K^3(-2D)\big)\) since \(\nu\) is holomorphic, even in the case \(\nu=0\). Finally, if two triples \((L_1, \gamma_1, \nu_1)\) and \((L_2, \gamma_2, \nu_2)\) are gauge equivalent via an isomorphism \(L_1 \cong L_2\), then their sections scale as \(\gamma_2 = \lambda \gamma_1\) and \(\nu_2 = \lambda^{-2} \nu_1\) for some \(\lambda \in \mathbb{C}^*\) by Lemma \ref{lem: gauge equivalent Higgs bundles}. As a consequence, both \(D = \operatorname{div}(\,\gamma_1\,) = \operatorname{div}(\,\gamma_2\,)\) and \(q_3 = \nu_1 \gamma_1^2 = \nu_2 \gamma_2^2\) are left invariant, by showing that \(\psi_d\) is independent of the choice of representative \([L, \gamma, \nu]\) within the isomorphism class.
\end{proofclaim}

Since the map \(\psi_d\) is holomorphic, the desired result readily follows once we prove that \(\psi_d\) is a bijection, and hence a biholomorphism. We begin by showing the injectivity of the mapping \(\psi_d\), \textit{i.e.}, we show that the isomorphism class of every triple $[L,\gamma,\nu]$ is uniquely determined by the pair $\big(\,\nu\gamma^2,\text{div}(\,\gamma\,)\,\big)$. For this purpose, let \([L_1,\gamma_1,\nu_1]\) and \([L_2,\gamma_2,\nu_2]\) be two triples and assume that
\begin{equation}
    \psi_d\big([L_1,\gamma_1,\nu_1]\big)=\psi_d\big([L_2,\gamma_2,\nu_2]\big)=(q_3,D).
\end{equation}

It readily follows that $D\defeq\text{div}(\,\gamma_1\,)=\text{div}(\,\gamma_2\,)$. Moreover, since the holomorphic sections \(\gamma_1\) and \(\gamma_2\) share the same zero divisor \(D\), each section induces a canonical isomorphism between its respective line bundle and \(\mathcal{O}_X(\,D\,)\), namely \(L_1^{-1}K \cong \mathcal{O}_X(\,D\,) \cong L_2^{-1}K\). Consequently, both \(L_1\) and \(L_2\) are isomorphic to the line bundle \(K \otimes \mathcal{O}_X(-D)\) and hence \(L_1\cong L_2\). In particular, since two holomorphic sections with the same divisors differ by some \(\lambda\in \C^*\), it readily follows that \(\gamma_2=\lambda\gamma_1\). Finally, by substituting this latter into the cubic differential relation yields \(\nu_2 (\lambda \gamma_1)^2 = \nu_1 \gamma_1^2\), which implies \(\nu_2 = \lambda^{-2} \nu_1\). By Lemma \ref{lem: gauge equivalent Higgs bundles}, the action of the constant \(\lambda \in \mathbb{C}^*\) corresponds precisely to the gauge equivalence for the triples \([L, \gamma, \nu]\), thus showing that the two triples represent the same isomorphism class \([L_1, \gamma_1, \nu_1] = [L_2, \gamma_2, \nu_2]\). In order to prove the surjectivity, let $(q_3,D)\in\mathcal{QD}_d(X)$ be any element. This datum determines a line bundle \(L\defeq K\otimes\mathcal{O}_X(-D)\) of degree \(d=2g-2-\deg(D)\in(0,2g-2]\), together with a holomorphic section, say \(\gamma\in H^0(X,\mathcal{O}_X(D))\) such that \(D=\text{div}(\gamma)\). Therefore, the section \(\nu\) is recovered by setting \(\nu\defeq q_3\gamma^{-2}\). Since \(q_3\in H^0\big(X,K^3(-2D)\big)\) by assumption, with \(D=\operatorname{div}(\gamma)\), the section \(\nu\) is holomorphic as desired.
\end{proof}

\begin{proof}[Proof of Proposition \ref{prop: parametrization of subspace higgs bundles2}] 
In what follows, we denote by $\mathrm{Sym}^n(\,X\,)$ the $n$-th symmetric product of the Riemann surface $X$, and by $\mathrm{Pic}^n(\,X\,)$ the set of line bundles of degree $n$. There is a biholomorphism between the set of pairs $(\mathcal{N},[\alpha])$, with $\mathcal{N}\in \mathrm{Pic}^n(\,X\,)$ and $[\alpha]\in \mathbb PH^0(X,\mathcal{N})$, and $\mathrm{Sym}^n(\,X\,)$, given by associating to $[\alpha]$ its divisor of zeros. Let \(d\) be a negative integer, $1-g\le d<0$. Recall that the space $\text{Sym}^{2d+2g-2}(\,X\,)$ admits an $H^1(X,\Z_2)$-cover denoted by ${\text{Sym}^{2d+2g-2}(\,X\,)}_0$ which parametrises pairs, say $(\mathcal N,[\delta])$, where $\mathcal N$ is a line bundle of degree $d+g-1$ and $[\delta]\in\mathbb{P}\big(H^0(X,\mathcal{N}^2)\setminus\{\,0\,\}\big)$. Once a square root of $K$ is fixed, we have a well-defined injective map 

\begin{equation}\label{eq: phi_d_map}
    \begin{aligned} 
    \phi_d \colon \mathcal{CS}_d^-(X) &\longrightarrow \mathrm{Sym}^{2d+2g-2}(X)_0 \times H^0(X, K^3) \\ [L, \gamma, \nu] &\longmapsto \left(\,\big(\,L K^{1/2}, [\nu]\,\big), \nu\gamma^2\,\right). 
    \end{aligned}
\end{equation}

Let $W$ be the image of $\phi_d$. Projection onto the first factor defines a surjective map $p\colon W \rightarrow \text{Sym}^{2d+2g-2}(\,X\,)_0$. The fibre $p^{-1}\left(\big(\,\mathcal{N},[\delta]\,\big)\right)$ comprises of all holomorphic cubic differentials that can be written as $\alpha^{2}\delta$ for some section $\alpha \in H^{0}\big(X, \mathcal{N}^{-1}K^{3/2}\big)$. Therefore, the fibre can be identified with the image of the quadratic map $H^{0}(X, \mathcal{N}^{-1}K^{3/2}) \setminus \{\,0\,\} \rightarrow H^{0}(X,\mathcal{N}^{-2}K^{3})$. In particular, it is biholomorphic to $\big(\,\C^{g-1-d} \setminus \{\,0\,\}\,\big) / \{\,\pm \mathrm{Id}\,\}$. 
\end{proof}

\noindent We are now ready to introduce the spaces of primary interest for the purposes of this section. 

\begin{defn}
    We shall define \(\mathcal{R}_d^+(\,X\,)\) and \(\mathcal{R}_d^-(\,X\,)\) as the subspaces of $\mathfrak{X}_3(S)$ corresponding to \(\mathcal{CS}_d^+(X)\) and \(\mathcal{CS}_d^-(X)\), respectively, via the non-abelian Hodge correspondence.
\end{defn}

\begin{rmk}\label{rem: slice R_d(X) contained in connected component}
    For $d=2g-2$, the space $\mathcal{CS}_d^+(\,X\,)$ is naturally identified with the space of holomorphic cubic differentials on $X$, and therefore the corresponding subspace of representations $\mathcal{R}_d^+(\,X\,)$ lies in the Hitchin component. In other words, it coincides with the fibre of the vector bundle of holomorphic cubic differentials over $\mathcal{T}(\,S\,)$, which is diffeomorphic to $\mathrm{Hit}_3(\,S\,)$, see \cite{labourie2007flat} and \cite[Theorem 2]{loftin2001}. More generally, as a direct application of Lemma \ref{lem: rho_connected_component_Higgs_bundles}, it readily follows that the spaces $\mathcal{R}_d^+(\,X\,)$ and $\mathcal{R}_d^-(X)$ are contained in $\mathcal{C}_0$ if and only if
    \(
    d\in (1-g,2g-2)\setminus\{\,0\,\}\cap 2\mathbb{Z},
    \)
    or $d=1-g$ with $g$ odd. In the same fashion, the spaces $\mathcal{R}_d^+(\,X\,)$ and $\mathcal{R}_d^-(X)$ are contained in $\mathcal{C}_1$ if and only if $d\in(1-g,2g-2)\setminus\{\,0\,\}$ is odd, or $d=1-g$ with $g$ even.   
\end{rmk}

The next result establishes that \(\mathcal{R}_d^+(\,X\,)\) is a properly embedded sub-variety in the character variety \(\mathfrak{X}_3(\,S\,)\) for every \(0<d \le 2g-2\). To fix notation, we recall for the readers' convenience that for such integer \(d\), the space \(\mathcal{QD}_d(\,X\,)\) is biholomorphic to \(\mathcal{CS}_d^+(\,X\,)\) via the map \(\Theta_d(q_3,D) \defeq [L,\gamma,\nu]\), where \(L \defeq K(-D)\), \(\gamma\) is the canonical section of \(\mathcal{O}_X(\,D\,)\), and \(\nu \defeq q_3\gamma^{-2}\). This map is nothing else than the inverse of the one introduced in \eqref{eq: holomorphic_map_CS_d}.

\begin{thm}\label{thm: proper_embedding}
    For every \(0<d\le 2g-2\), the map $\Theta_d\colon\mathcal{QD}_d(X)\to\mathcal{M}_3(X)^s$ is a proper embedding.
\end{thm}

\begin{proof}
According to Proposition \ref{prop: parametrization of subspace higgs bundles}, the mapping \(\Theta_d\) is a biholomorphism onto its image, which coincides with \(\mathcal{CS}_d^+(\,X\,)\). Therefore, it is sufficient to prove that \(\Theta_d\) is proper map. In his work \cite{Hitchin1992}, Hitchin shows the existence of a natural fibration map
\begin{equation}
    \pi\colon\mathcal{M}_3(X)^s\longrightarrow H^0(X,K^2)\oplus H^0(X,K^3)
\end{equation}
that applies to an isomorphism class \([\mathcal{E},\Phi]\) the pair \(\left(\mathrm{tr}(\Phi^2),\mathrm{tr}(\Phi^3)\right)\). Let \(\mathrm{pr}_3\) be the projection on the second factor, \textit{i.e.}, the projection onto the space of holomorphic cubic differentials. Let \(C \subset \mathcal{M}_3(\,X\,)^s\) be any compact subset. Since the projection map is continuous, it readily follows that \(B_C \defeq \mathrm{pr}_3(\pi(C))\) is compact in the space \(H^0(X, K^3)\). Moreover, we may observe that, if \((q_3,D)\in \Theta_d^{-1}(\,C\,)\), then \(q_3\in B_C\). As a consequence, the fibre \(\Theta_d^{-1}(\,C\,)\) lies in the intersection \(\mathcal{QD}_d(\,X\,)\,\cap\, \big(B_C\times \mathrm{Sym}^N(\,X\,)\big)\), where \(N\defeq2g-2-d\). Since \(X\) is compact, so is the space \(\operatorname{Sym}^N(\,X\,)\), and hence the space \(B_C \times \operatorname{Sym}^N(\,X\,)\) is compact as well. The following technical statement holds.

\begin{claim}\label{claim: closure in product space}
    The space \(\mathcal{QD}_d(\,X\,)\) is closed in \(H^0(X, K^3) \times \operatorname{Sym}^N(\,X\,)\).
\end{claim}

Assume Claim \ref{claim: closure in product space} holds. Then the intersection
\(
        \mathcal{QD}_d(\,X\,) \cap \bigl(B_C \times \operatorname{Sym}^N(\,X\,)\bigr)
\)
is compact. Moreover, as \(C\) is compact in the Hausdorff space \(\mathcal{M}_3(X)^s\), the preimage \(\Theta_d^{-1}(C)\) is closed in \(\mathcal{QD}_d(X)\) by continuity of the mapping \(\Theta_d\). As a consequence, \(\Theta_d^{-1}(C)\) is closed in the compact set \(\mathcal{QD}_d(X)\cap\bigl(B_C\times \operatorname{Sym}^N(X)\bigr).\) It readily follows that \(\Theta_d^{-1}(C)\) is compact. To complete the proof, it remains to show Claim~\ref{claim: closure in product space}.

\begin{proofclaim}[Proof of Claim~\ref{claim: closure in product space}]
This can be verified in local holomorphic coordinates and is a consequence of routine complex analysis. Let \((q_{3,n}, D_n) \longrightarrow (q_{3,\infty}, D_\infty)\) be a convergent sequence such that \(q_{3,n}\in H^0\big(X,K^3(-2D_n)\big)\) for every \(n\). We aim to show that \(q_{3,\infty} \in H^0\bigl(X, K^3(-2D_\infty)\bigr)\). If \(q_{3,\infty}=0\) the conclusion is immediate. Thus, we may assume that $q_{3,\infty}\neq 0$. In particular, $q_{3,n}\neq 0$ for all sufficiently
large $n$. Fix \(p \in \operatorname{supp}(D_\infty)\), and let \(r\) be the multiplicity of \(p\) in \(D_\infty\). Choose a coordinate disc \(U\) centred at \(p\) whose closure contains no points of \(\operatorname{supp}(D_\infty)\) other than \(p\). Since \(D_n \longrightarrow D_\infty\) in \(\operatorname{Sym}^N(X)\), after shrinking \(U\) if necessary, for all sufficiently large \(n\) the restriction \({D_n|_U}\) has degree \(r\), with all points in its support converging to \(p\). In the corresponding local trivialisation of \(K^3\), write \(q_{3,n} = f_n(z) \, dz^3\) and \(q_{3,\infty} = f_\infty(z) \, dz^3\). The convergence \(q_{3,n} \longrightarrow q_{3,\infty}\) implies that \(f_n\) converges to \(f_\infty\) uniformly on compact subsets of \(U\). The standing assumption \(q_{3,n}\in H^0\big(X,K^3(-2D_n)\big)\) implies that \({2D_n|_U}\) contributes \(2r\) zeroes to \(f_n\), counted with multiplicity, whose support converges to \(p\). For the sake of contradiction suppose that \(s \defeq \operatorname{ord}_p(f_\infty) < 2r\). Choose a smaller coordinate disc \(V\) centred at \(p\), with \(\overline{V} \subset U\), such that \(p\) is the only zero of \(f_\infty\) in \(\overline{V}\) and \(f_\infty\) does not vanish on \(\partial V\). For every sufficiently large \(n\), all points of \({D_n|_U}\) lie in \(V\), and therefore \(f_n\) has at least \(2r\) zeroes in \(V\), counted with multiplicity. On the other hand, since \(f_n \to f_\infty\) uniformly on \(\partial V\) and \(f_\infty\) does not vanish on \(\partial V\), for all sufficiently large \(n\) we have
\begin{equation}
    \lvert\,f_n(\,z\,) - f_\infty(\,z\,)\, \rvert < \lvert\,f_\infty(\,z\,)\, \rvert
\end{equation}
for every \(z \in \partial V\). By Rouché's theorem \cite[Chapter~V, Theorem~3.8]{conway1973functions}, \(f_n\) and \(f_\infty\) have the same number of zeroes in \(V\), counted with multiplicity. Since \(p\) is the only zero of \(f_\infty\) in \(V\), the number of zeroes of \(f_\infty\) in \(V\), counted with multiplicity, is precisely \(\operatorname{ord}_p(f_\infty) = s\). By our assumption \(s < 2r\), this contradicts the fact that \(f_n\) has at least \(2r\) zeroes in \(V\). Hence, \(\operatorname{ord}_p(f_\infty) \ge 2r\). Since this holds for every \(p \in \operatorname{supp}(D_\infty)\), we conclude that \(q_{3,\infty} \in H^0\bigl(X, K^3(-2D_\infty)\bigr)\).
\end{proofclaim}
Since Claim~\ref{claim: closure in product space} holds, the desired conclusion follows.
\end{proof}

The next result establishes that \(\mathcal{R}_d^-(\,X\,)\) is a properly embedded sub-variety in the irreducible locus of the character variety \(\mathfrak{X}_3(S)\) for every \(1-g\leq d < 0\), extending a result for the minimal degree \(d = 1-g\) recently obtained by Bronstein--Davalo in \cite[Proposition~2.5]{bronstein2025anosov}. Unlike the previous scenario, the present one is more subtle. Let \(\mathcal{F}_d\) denote the holomorphic bundle over the \(H^1(X,\mathbb{Z}_2)\)-cover of the \((2d+2g-2)\)-th symmetric product of \(X\), which has been already considered in the proof of Proposition \ref{prop: parametrization of subspace higgs bundles}. In what follows, we shall denote such a cover as \(B_d\). After fixing a square root \(K^{1/2}\) of the canonical bundle, a point of \(\mathcal{F}_d\) is specified by a triple, say \((\mathcal{N},[\nu],q_3)\) where: $\mathcal{N}$ is a line bundle of degree $d+g-1$, $[\nu]\in\mathbb{P}\big(H^0(X,\mathcal{N}^2)\setminus\{\,0\,\}\big)$ and \(q_3=\gamma^2\nu\in H^0(X,K^3)\), for some \(\gamma\in H^0(X,\mathcal{N}^{-1}K^{3/2})\).

\begin{thm}\label{thm: proper embedding d negative}
    For every $1-g\le d <0$, the map 
    \begin{align*}
        \Psi&_d: \ \mathcal{F}_d\longrightarrow\mathcal{M}_3(X)^s \\ & (\mathcal{N},[\nu],q_3)\mapsto[\mathcal{N}K^{-\frac{1}{2}},\gamma,\nu]
    \end{align*} 
    is a proper embedding onto $\mathcal{CS}_d^-(X)$.
\end{thm}

\begin{proof}
By Proposition \ref{prop: parametrization of subspace higgs bundles2}, the map \(\Psi_d\) is a biholomorphism onto \(\mathcal{CS}_d^-(X)\). Therefore, it sufficient to prove that \(\Psi_d\) is proper map. Let \(C\subset\mathcal M_3(X)^s\) be compact, and let \(x_n=(N_n,[\nu_n],q_{3,n})\in\Psi_d^{-1}(C)\) be a sequence. Since \(B_d\to\operatorname{Sym}^{2d+2g-2}(X)\) is a finite covering and \(\operatorname{Sym}^{2d+2g-2}(X)\) is compact, the space \(B_d\) is compact as well. Thus, after passing to a subsequence, we may assume that \((N_n,[\nu_n])\to(N_\infty,[\nu_\infty])\) in \(B_d\). Choose local trivialisations around \((N_\infty,[\nu_\infty])\), and choose representatives \(\nu_n\in H^0(X,N_n^2)\) such that \(\nu_n\to\nu_\infty\neq0.\) For every \(n\), choose \(\gamma_n\in H^0(X,N_n^{-1}K^{3/2})\setminus\{\,0\,\}\) satisfying \(q_{3,n}=\nu_n\gamma_n^2.\) Since the Hitchin fibration is continuous and \(\Psi_d(x_n)\in C\), the compactness of \(C\) implies that the sequence of cubic differentials \(q_{3,n}\in H^0(X,K^3)\) is bounded. We claim that the same holds for \(\gamma_n\). Indeed, suppose otherwise. After passing to a subsequence, we would have \(\|\gamma_n\|\to+\infty.\) Define \(\widehat\gamma_n\defeq\frac{\gamma_n}{\|\gamma_n\|}\) and after passing to a further subsequence,
\begin{equation}
    \widehat\gamma_n\longrightarrow\widehat\gamma_\infty,\qquad\|\widehat\gamma_\infty\|=1.
\end{equation}

Therefore,
\begin{equation}
    \frac{q_{3,n}}{\,\,\|\gamma_n\|^2}=\nu_n\widehat\gamma_n^2\longrightarrow\nu_\infty\widehat\gamma_\infty^2.
\end{equation}

Since both \(\nu_\infty\) and \(\widehat\gamma_\infty\) are non-zero, their product
\(\nu_\infty\widehat\gamma_\infty^2\) is a non-zero holomorphic cubic differential. Consequently,
\begin{equation}
    \|q_{3,n}\|=\|\gamma_n\|^2\left\|\frac{q_{3,n}}{\,\,\|\gamma_n\|^2}\right\|\longrightarrow+\infty,
\end{equation}

contradicting the boundedness of \(q_{3,n}\). Therefore \(\gamma_n\) must be bounded. After passing to another subsequence, we may assume that \(\gamma_n\to\gamma_\infty\) for some \(\gamma_\infty\in H^0(X,N_\infty^{-1}K^{3/2}).\) We claim that \(\gamma_\infty\neq0.\)
Indeed, if \(\gamma_\infty=0\), the corresponding Higgs bundles would converge in the poly-stable moduli space to the Higgs bundle
\begin{equation}
    L_\infty\oplus\mathcal O_X\oplus L_\infty^{-1},\quad 
    \Phi_\infty=\begin{pmatrix}
    0&0&\nu_\infty\\
    0&0&0\\
    0&0&0
    \end{pmatrix},
    \quad
    L_\infty\defeq N_\infty K^{-1/2}.
\end{equation}

By Proposition \ref{prop: stable SL(3;C) Higgs bundle}, this limit is strictly poly-stable. On the other hand, since \(\Psi_d(x_n)\in C\) and \(C\) is compact, after passing to a subsequence \(\Psi_d(x_n)\) converges to a point \(y\in C\subset\mathcal M_3(X)^s.\) The poly-stable moduli space is Hausdorff, so the limit is unique. This is impossible, since one limit is stable and the other is strictly poly-stable. Hence \(\gamma_\infty\neq0\) and it follows that
\begin{equation}
x_n\longrightarrow
x_\infty\defeq
\bigl(N_\infty,[\nu_\infty],
\nu_\infty\gamma_\infty^2\bigr)
\end{equation}
in \(\mathcal F_d\). Moreover \(\Psi_d(x_\infty)\in C\), by continuity. Therefore every sequence in \(\Psi_d^{-1}(C)\) admits a convergent subsequence whose limit still belongs to \(\Psi_d^{-1}(C)\). Hence \(\Psi_d^{-1}(C)\) is compact, and \(\Psi_d\) is proper as desired.
\end{proof}
\noindent The proof shows explicitly that the submanifolds $\mathcal{R}_d^-(X)$ are properly embedded in the irreducible locus of the character variety, but not in $\mathfrak{X}_3(S)$ itself, in contrast with $\mathcal{R}_d^+(X)$. Moreover, we observe that $\mathcal{R}_d^+(\,X\,)$ is a properly embedded topological disks in $\mathfrak{X}_3(\,S\,)$ if and only if $d=2g-2$, in which case \(\mathcal{R}_d^+(X)\) has real dimension $10g-10$. For all \(d \in [1-g, 2g-2) \setminus \{\,0\,\}\), the topology of $\mathcal{R}_d^+(\,X\,)$ and $\mathcal{R}_d^-(X)$ is non-trivial, as established in Proposition~\ref{prop: parametrization of subspace higgs bundles} and \ref{prop: parametrization of subspace higgs bundles2}. In the next section, we shall see how the representations in the subspaces $\mathcal{R}_d^+(\,X\,)$ and $\mathcal{R}_d^-(X)$ admit a geometric description in terms of branched \(\rp2\)-structures on \(S\), endowed with a special conical metric at the branch points.

\section{Construction of branched projective structures of hyperbolic affine sphere type}\label{sec:cyclic Higgs bundles and BGS}

Since the seminal work of Baraglia in his doctoral dissertation \cite{baraglia2009g2}, Higgs bundles have been used to realise geometric \((G, \mathbb{X})\)-structures on manifolds, as already alluded to in \S~\ref{ssec:higgs bundles}. Following the classical approach, given a Higgs bundle and its associated flat bundle, a geometric structure is determined by a section transverse to the parallel foliation; see \cite{alessandrini2019}. In the present section we aim to extend this construction to the branched setting. More precisely, in \S~\ref{sec: from cyclic to geometric structres} we shall illustrate how Higgs bundles can be employed to realise branched \(\mathbb{R}\mathrm{P}^2\)-structures on surfaces, together with a Riemannian metric having singularities at the branch points and arising from a conical hyperbolic affine sphere (\S~\ref{sec: affine sphere singular}). This provides a geometric description of the representations forming the subspaces \(\mathcal{R}_d^+(X)\) and \(\mathcal{R}_d^-(X)\) defined in \S\ref{ssec: cyclic higgs bundles}. In \S\ref{sec: prescribed data problem}, we finally address the existence problem for branched structures with prescribed branching data.

\subsection{Affine spheres with conical singularities}\label{sec: affine sphere singular}
In this section, we briefly recall the definition of a hyperbolic affine sphere in \(\mathbb{R}^3\), one of the central objects in affine differential geometry \cite{nomizu1994affine} and in the Labourie--Loftin parametrisation of the Hitchin component \(\operatorname{Hit}_3(S)\); see \cite{loftin2001,labourie2007flat}. Starting from this classical notion, we introduce a version with conical singularities and show that, also in this setting, one can naturally associate a holomorphic cubic differential with the surface, subject to restrictions on the orders of its zeros imposed by the branching data. 

\smallskip

Let \(S\) be a connected and oriented surface. Let \(f\colon S\longrightarrow\mathbb{R}^3\) be a locally strictly convex affine immersion and let \(\nabla\) denote the standard flat connection on \(\mathbb{R}^3\). Given a transverse vector field \(\xi\) along \(f\), one can write
\begin{equation}
    \nabla_Xf_*Y=f_*(\widehat\nabla_XY)+G(X,Y)\xi \qquad \nabla_X\xi=-f_*(BX)+\tau(X)\xi,
\end{equation}

where \(\widehat{\nabla}\) is a torsion-free connection on \(S\), \(G\) is a symmetric bilinear form, \(B\) is the affine shape operator, and \(\tau\) is a one-form. The transverse field \(\xi\) is called \emph{equiaffine} if \(\tau=0\). The \emph{Blaschke affine normal} is the unique equiaffine transverse vector field such that the induced volume form \(\theta_\xi(X,Y)\defeq\det\bigl(f_*X,f_*Y,\xi\bigr)\) coincides with the area form of \(G\), namely \(\theta_\xi = \mathrm{d}A_G\). The corresponding metric \(G_B \defeq G\) is called the Blaschke metric. The immersion \(f\) is called an affine sphere if \(B=\lambda\,\mathrm{Id}\) for some constant \(\lambda\), and it is called hyperbolic if \(\lambda<0\). Normalising \(\lambda=-1\) and placing the centre at the origin yields \(\xi=f\). The structure equations therefore become
\begin{equation}\label{eq: structure equations hyperbolic affine spheres}
    \nabla_Xf_*Y=f_*(\widehat\nabla_XY)+G_B(X,Y)f,
    \qquad
    \nabla_Xf=f_*(X).
\end{equation}

\begin{defn}\label{def:singular affine sphere}
    Let $S$ be a connected oriented surface, let $\mathcal P\subset S$ be a discrete subset, and let $m\colon\mathcal P\to\mathbb Z_{\geq 1}$. A continuous map $f\colon S\to\mathbb R^3\setminus\{0\}$ is called a \emph{hyperbolic affine sphere with conical singularities} along $\mathcal P$ of orders $m$ if the following conditions hold:
    \begin{enumerate}
        \item[(i)] The restriction $f|_{S\setminus \mathcal P}:S\setminus \mathcal P\to\mathbb R^3$ is an immersed hyperbolic affine sphere in the classical sense, normalised so that its affine shape operator is $B=-\mathrm{Id}$ and its centre is the origin.
        \item[(ii)] For every $p\in \mathcal P$, there exists a smooth orientation-preserving local coordinate $z$ centred at $p$ such that the Blaschke metric takes the form \(G_B=|\,z\,|^{2m(p)}e^{2u(z)}|\,dz\,|^2,\) where $u$ is smooth in a neighbourhood of the origin. The extended tensor on \(S\) will be called \emph{conical Blaschke metric}.
    \end{enumerate}
\end{defn}

Notice that the conical condition on the metric \(G_B\) implies that the induced complex structure on \(S_0 \defeq S \setminus \mathcal{P}\) extends to the whole surface. Moreover, for the sake of brevity, we shall simply refer to such a map as a \emph{conical hyperbolic affine sphere} whenever the branch locus and the orders are clear from the context. 

\smallskip

From now on, we restrict our attention to the case where \(S\) is a closed oriented surface of genus \(g \geq 2\). As in the classical case, the surface \(S_0\) naturally inherits the \emph{Pick tensor} \(C \defeq \widehat{\nabla} G_B\). Denoting by \(J\) the complex structure induced by the conformal class of the Blaschke metric \(G_B\), we obtain a \(J\)-holomorphic cubic differential on \(S_0\), called the \emph{Pick differential} and denoted by  \(q_{\mathrm{Pick}}\). It is well known that the Pick tensor satisfies the identity \(C = 2\operatorname{Re}(q_{\mathrm{Pick}})\); see, for instance, \cite[\S4.3]{benoist2013cubic}.

\begin{lem}\label{lem: singular hyperbolic affine spheres}
    Let \(f:S\to\mathbb{R}^3\) be a hyperbolic affine sphere with conical singularities at \(\mathcal P=\{p_1,\ldots,p_k\}\) of orders \(\underline m=\{m_1,\ldots,m_k\}\). Then the Pick differential \(q_{Pick}\) extends to a holomorphic cubic differential on the Riemann surface \(X=(S,J)\) such that \(q_{Pick}\in H^0\big(X,K^3(-2D)\big)\), being \(D=\sum_{i=1}^km_i\cdot p_i\).
\end{lem}

\begin{proof}
    Let \(p_i \in \mathcal{P}\) be any point and let \(z\) be a complex coordinate centred at \(p_i\) such that \(G_B = |\,z\,|^{2m_i} e^{2u_i} |\,\mathrm{d}z\,|^2\), where \(u_i\) is a smooth function defined in a neighbourhood \(U_i\) of \(p_i\). On \(U_i\setminus\{p_i\}\), we may write \(q_{\mathrm{Pick}} = Q(\,z\,) \, \mathrm{d}z^3\), where \(Q\) is a holomorphic function. As a consequence of Wang's equation, see \cite[\S6.2]{li2019introduction}, the identity
    \begin{equation}
        K_{G_B} = -1 + 2 \|\,q_{\mathrm{Pick}}\,\|_{G_B}^2
    \end{equation}
    follows. On the other hand, since \(K_{G_B} = -4|\,z\,|^{-2m_i} e^{-2u_i} \partial_z \partial_{\bar{z}} u_i\) on \(U_i\) and \(\|\,q_{\mathrm{Pick}}\,\|_{G_B}^2 = |\,Q(\,z\,)\,|^2 |\,z\,|^{-6m_i} e^{-6u_i}\), the following identity readily follows
    \begin{equation}
        -4|\,z\,|^{-2m_i} e^{-2u_i} \partial_z \partial_{\bar{z}} u_i = -1 + 2|\,Q(\,z\,)\,|^2 |\,z\,|^{-6m_i} e^{-6u_i}.
    \end{equation}
    Equivalently,
    \begin{equation}
        2|\,Q(\,z\,)\,|^2 = |\,z\,|^{6m_i} e^{6u_i} - 4|\,z\,|^{4m_i} e^{4u_i} \partial_z \partial_{\bar{z}} u_i.
    \end{equation}
    Since \(u_i\) is smooth, the right-hand side is \(O(|\,z\,|^{4m_i})\), and hence \(|\,Q(\,z\,)\,| = O(|\,z\,|^{2m_i})\). Let \(F(z) \defeq Q(\,z\,)\,z^{-2m_i}\). We may notice that \(F(z)\) is holomorphic and bounded on the punctured disc \(U_i\setminus\{\,p_i\,\}\). Recall that \(z(\,p_i\,)=0\). The removable singularity theorem applies and then the function \(F\) extends holomorphically across \(0\). As a consequence, \(Q(\,z\,) = z^{2m_i} F(\,z\,)\) with \(F\) holomorphic near \(0\), and therefore \(\operatorname{ord}_{p_i}(q_{\mathrm{Pick}}) \geq 2m_i\). Since \(p_i\) has been chosen arbitrarily, the desired claim follows.
\end{proof}

These new geometric objects allow us to define a refined version of branched geometric structures as follows.

\begin{defn}\label{def: branched conical projective}
    Let $\sigma$ be a branched $\mathbb{RP}^2$-structure on $S$, with developing--holonomy pair $(\operatorname{dev},\rho)$ and branching divisor \( D=\sum_{i=1}^k m_i p_i.\) We say that $\sigma$ is of \emph{hyperbolic affine sphere type} if its developing map admits a $\rho$-equivariant lift \( f\colon \widetilde S\longrightarrow\mathbb R^3\setminus\{0\}\) which is a hyperbolic affine sphere with conical singularities of orders $m_i$ along the lifts of the points $p_i$, in the sense of Definition \ref{def:singular affine sphere}.
\end{defn}

By Definition~\ref{def: branched conical projective}, any hyperbolic affine-sphere lift of a branched \(\rp2\)-structure of the above type induces a singular Riemannian metric on \(S\) with integral conical singularities along the branch locus, namely its conical Blaschke metric \(G_B\). We provide a basic, yet motivating, example that led us to Definition~\ref{def: branched conical projective}. As we shall explain, hyperbolic structures with integral conical singularities can be regarded as a special case of the structures introduced above.

\begin{ex}
        We briefly recall that the hyperboloid model $\mathbb{H}^2=\{x\in\mathbb{R}^3\mid x_0^2+x_1^2-x_2^2=-1,\ x_2>0\}$ realises the hyperbolic plane as a submanifold of $\mathbb{R}^3$, in such a way that the natural inclusion is equivariant with respect to the irreducible representation $\jmath\colon\mathrm{PSL}(2,\mathbb{R})\to\mathrm{SL}(3,\mathbb{R})$. Suppose now that $S$ is a closed surface of genus $g\ge2$ endowed with a branched hyperbolic structure with branched points at \(\mathcal{P} = \{p_1, \ldots, p_k\}\subset S\). In terms of its developing--holonomy pair, this amounts to a \(\rho_o\)-equivariant map \(\operatorname{dev} \colon \widetilde{S} \longrightarrow \mathbb{H}^2\) which is a local diffeomorphism away from a discrete set of points and is locally a branched covering of order \(m_i \in \mathbb{Z}_{\geq 1}\) at each pre-image of \(p_i \in \mathcal{P} = \{p_1, \ldots, p_k\}\); see \cite{faraco2021, faracomaret2025, nguyen2021hyperbolic} for further details. By composing the developing map with the natural inclusion $\mathbb{H}^2\to\mathbb{R}^3$, we obtain a map $f\colon\widetilde{S}\to\mathbb{R}^3$ which is equivariant with respect to the representation $\rho\defeq \jmath\circ\rho_o\colon\pi_1(S)\to\mathrm{SL}(3,\mathbb{R})$ and defines a conical hyperbolic affine sphere in the sense of Definition~\ref{def:singular affine sphere}. In this case, the conical Blaschke metric $G_B$ has constant negative curvature away from the branch points, while the Pick tensor vanishes identically. Finally, by projectivising the hyperboloid model and the map $f$, one obtains a branched projective structure in the sense of Definition~\ref{def: branched conical projective}.
\end{ex}

\subsection{Generalising Baraglia's construction}\label{sec: from cyclic to geometric structres} This section contains the main construction of the paper, in which we use Higgs bundles to build families of branched projective structures of hyperbolic affine-sphere type on the surface. Let \(X\) be a fixed Riemann surface structure on \(S\), and consider the family of stable cyclic \(\mathrm{SL}(3,\mathbb{R})\)-Higgs bundles
\begin{equation}
    \mathcal{E} = L \oplus \mathcal{O}_X \oplus L^{-1}, \quad \Phi = \begin{pmatrix}
0 & 0 & \nu \\
\gamma & 0 & 0 \\
0 & \gamma & 0
\end{pmatrix}, \quad \operatorname{deg}(L) \in [1-g, 2g-2] \setminus \{0\}
\end{equation}
which parametrise the subspaces \(\mathcal{C}_d^+(X)\) and \(\mathcal{C}_d^-(X)\) of \(\mathcal{M}_3(X)^s\) defined in \eqref{eq: subspace of cyclic higgs bundles}. By Theorem~\ref{thm: proper_embedding} and Theorem~\ref{thm: proper embedding d negative}, we obtain properly embedded submanifolds \(\mathcal{R}_d^+(X)\) (resp.\ \(\mathcal{R}_d^-(X)\)) of the character variety \(\mathfrak{X}_3(S)\) (resp.\ the irreducible locus), whose points correspond to the holonomies of the associated flat connections. The main result of this section is the following:

\begin{thm}\label{thm: every rho is holonomy of conical branched} 
    Let \(X\) be a Riemann surface structure on $S$. Then every representation \(\rho\) in $\mathcal{R}_d^+(X)$ and $\mathcal{R}_d^-(X)$ arises as the holonomy of a branched projective structure of hyperbolic affine sphere type.    
\end{thm} 

In the following we split the proof of Theorem~\ref{thm: every rho is holonomy of conical branched} into Lemma \ref{lem: branche structures from Higgs bundles} and Proposition \ref{thm: developing map lifts to a singular affine spheres}. Firstly, we recall for the readers' convenience the setting and introduce the objects required for the proof. By definition, see \eqref{eq: subspace of cyclic higgs bundles}, the moduli spaces \(\mathcal{C}_d^+(X)\) and \(\mathcal{C}_d^-(X)\) consist of equivalent classes of triples \((L,\gamma,\nu)\) for which \(\gamma \neq 0\) in all cases, whereas \(\nu \neq 0\) is required only when \(d = \operatorname{deg}(L) \in [1-g,0)\). After a suitable gauge transformation, the Hermitian metric \(H\) on \(\mathcal{E}\) solving the Hitchin equations \eqref{eq: Hitchin_higgs bundles} is given by \(H = \operatorname{diag}(h,1,h^{-1})\), where \(h\) is a Hermitian metric on the line bundle \(L\), see Remark~\ref{rem: cyclic_higgs_bundles}. In this case, consider the involution \(\tau \colon \mathcal{E} \to \mathcal{E}\) defined as 
\begin{equation}
    \tau(v_1,v_2,v_3) \defeq \big(h^{-1}\overline{v}_3, \overline{v}_2, h\overline{v}_1\big),
\end{equation}
where we use the \(\mathbb{C}\)-anti-linear identification induced by the Hermitian metric
\begin{equation}
    h \colon L \longrightarrow L^{-1}, \qquad v \longmapsto h(v,\cdot).
\end{equation}
The fixed-point locus of \(\tau\) defines a real rank-\(3\) sub-bundle of \(\mathcal{E}\), given by
\begin{equation}
E_{\mathbb{R}} \defeq \operatorname{Fix}(\tau) = \left\{ (v,\lambda,\overline{v}) \in \mathcal{E} \ \middle|\ v \in L,\ \lambda \in \underline{\mathbb{R}} \right\}.
\end{equation}
It is endowed with the flat connection \(\nabla = \nabla^H + \Phi + \Phi^{*_H}\), where \(\nabla^H\) denotes the Chern connection with respect to \(H\), and
\begin{equation}
    \Phi^{*_H} = H^{-1}\overline{\Phi}^{\,t}H = 
    \begin{pmatrix}
        0 & h^{-1}\overline{\gamma} & 0 \\
        0 & 0 & h^{-1}\overline{\gamma} \\
        h^2\overline{\nu} & 0 & 0
    \end{pmatrix}.
\end{equation}

Let \((U,z)\) be a local holomorphic chart on \(X\). The Hermitian metric \(h\) on \(L\) is locally represented with respect to a local holomorphic section \(\ell \colon U \subset X \to L\) by \(h(\ell,\ell) = h_0 = e^{2w}\), where \(w\) is a smooth real-valued function defined on \(U\). By using this notation, the Chern connection takes the form \(\nabla^H = \nabla^{1,0} + \overline{\partial}_{\mathcal{E}} + \mathcal{A}\), where  \( \mathcal{A} = \operatorname{diag}(2\partial w, 0, -2\partial w) \) is the local matrix of \(1\)-forms. By denoting by \(s\) the section spanning the real rank-\(1\) sub-bundle \(\underline{\mathbb{R}} \subset E_{\mathbb{R}}\), namely \(s(p) \defeq (0,1,0)\) for any \(p \in U \subset X\), we readily obtain an induced section, say \(\sigma \colon U \subset X \to \mathbb{P}(E_{\mathbb{R}})\) defined by \(\sigma(p) \defeq [s(p)]\). Since \((E_{\mathbb{R}},\nabla)\) is a flat vector bundle with holonomy \(\rho\), its projectivisation is naturally identified with the flat projective bundle \(\mathbb{P}(E_{\mathbb{R}}) \cong (\widetilde{X} \times \rp2) / \sim_{\rho}\). By pulling back \(\sigma\) to the universal cover and using this identification, we obtain a section of \(\widetilde{X} \times \rp2\), the trivial projective bundle. This section is the graph of a smooth map \(\operatorname{dev}_{\sigma} \colon \widetilde{X} \longrightarrow \rp2\) which is \(\rho\)-equivariant by construction, see also \cite[Proposition~4.1]{alessandrini2019}. Furthermore, the differential of \(\operatorname{dev}_{\sigma}\) is represented by \(\nabla s\) modulo the line spanned by \(s\).

\smallskip

Recall that \(\gamma\neq0\) is a holomorphic section of \(L^{-1}K\) and hence it has exactly \(2g-2-d\ge0\) zeroes counted with multiplicity. We shall denote by \(Z(\gamma)\defeq \{p_i\}\subset X\) the finite subset of its zeroes and by \(\widetilde Z(\gamma)\subset\widetilde X\) the discrete subset of the universal cover obtained by lifting \(Z(\gamma)\).

\begin{lem}\label{lem: branche structures from Higgs bundles}
    The section \(\sigma\) is transverse to the parallel foliation of \(\mathbb{P}(E_{\mathbb{R}})\) precisely on \(X \setminus Z(\gamma)\). Moreover, the \(\rho\)-equivariant map \(\operatorname{dev}_\sigma:\widetilde X\longrightarrow\mathbb{RP}^2\) satisfies the following properties:
    \begin{itemize}
        \item[(i)]\(\mathrm{dev}_\sigma\) is a local diffeomorphism on \(\widetilde X\setminus \widetilde Z(\gamma)\); 
        \item[(ii)] at each point \(\widetilde{p_i}\in\widetilde Z(\gamma)\) it is a branched covering of order \(m_i\in\mathbb Z_{\ge 1}\), where \(m_i\) is the multiplicity of \(p_i\) as a zero of \(\gamma\) and \(\widetilde{p_i}\) is any of its lift.
    \end{itemize}
\end{lem}

\begin{proof}
Let \(p \in \widetilde{X} \setminus \widetilde{Z}(\gamma)\) and let \((U,z)\) be a local holomorphic coordinate centred at \(p\). According to \cite[Proposition~4.5]{alessandrini2019}, the section \(\sigma\) is transverse to the parallel foliation at \(p\) if and only if \(\{s, \nabla_{\partial_z}s, \nabla_{\partial_{\bar{z}}}s\}\) are linearly independent on \(U\). With respect to the local coordinate \(z\), we may write \(h, \gamma\) and \(\nu\) as
\begin{equation}
    h=h_0\,{\overline{\ell}} \ ^{-1}\otimes\ell^{-1},\qquad\gamma=\gamma_0\,\ell^{-1}\otimes dz,\qquad
\nu=\nu_0\,\ell^2\otimes dz,
\end{equation}

where \(\ell\) a local holomorphic section of \(L\) and \(h_0=e^{2w}\). By combining the computations carried out in the previous paragraph, it follows that
\begin{equation}
    \mathcal A+\Phi+\Phi^{*^H}=
    \begin{pmatrix}
        \,\,2\partial w                        & h^{-1}_0\overline\gamma_0\mathrm d\bar z & \nu_0\mathrm d z \,\,\\
        \,\,\gamma_0\mathrm dz                 & 0                                        & h^{-1}_0\overline\gamma_0\mathrm d\bar z \,\,\\ 
        \,\,h^2_0\overline\nu_0\mathrm d\bar z & \gamma_0\mathrm d z                      & -2\partial w \,\,
    \end{pmatrix}.
\end{equation}

Moreover, since \(s(p)\defeq(0,1,0)\), an explicit computation shows that
\begin{equation}
    \nabla_{\partial_z}s=\big(\mathcal A(\partial_z)+\Phi(\partial_z)\big)\cdot s=
    \begin{pmatrix}
        0 \\
        0 \\ 
        \gamma_0
    \end{pmatrix},\quad
    \nabla_{\partial_{\overline z}}s=\Phi^{*^H}(\partial_{\overline z})\cdot s=
    \begin{pmatrix}
        h^{-1}_0\overline\gamma_0 \\
        0 \\
        0
    \end{pmatrix}.
\end{equation}

Therefore, \(\{s,\nabla_{\partial_z}s,\nabla_{\partial_{\overline z}}s\}\) at linearly independent on \(U\). Since \(p \in \widetilde{X} \setminus \widetilde{Z}(\gamma)\) has been chosen arbitrarily, we readily conclude that \(\{s,\nabla_{\partial_z}s,\nabla_{\partial_{\overline z}}s\}\) are linearly independent at every point \(p\in X\) such that \(\gamma(p)\neq 0\). On the other hand, \(\sigma\) is tangent to the parallel foliation along the finite subset \(Z(\gamma)\). As a consequence, the smooth map \(\mathrm{dev}_\sigma\) is a local diffeomorphism on \(\widetilde X\setminus\widetilde Z(\gamma)\), thus proving \((i)\) as desired. 
It remains to study the local behaviour of \(\mathrm{dev}_\sigma\) on \(\widetilde{Z}(\gamma)\). For this purpose, let \(p_i \in Z(\gamma)\) be a zero of multiplicity \(m_i\). Let \((U,z)\) be a local holomorphic chart centred at \(p_i\), and let \(\widetilde{p}_i\) be a lift of \(p_i\) to \(\widetilde{S}\), so that \(\gamma|_U = z^{m_i}\gamma_0(\,z\,)\,dz \otimes \ell^{-1}\) for some non-vanishing holomorphic function \(\gamma_0\). With a slight abuse of notation, we retain \(z\) to denote a local coordinate on \(\widetilde{U}\) centred at \(\widetilde{p}_i\). We then select a local affine chart on \(\rp2\) sending \(\operatorname{dev}_\sigma(\widetilde{p}_i)\) to the origin. Under these chart identifications, \(\operatorname{dev}_\sigma\) is locally represented by a smooth map
\begin{equation}
F \colon (\widetilde{U},0) \longrightarrow (\mathbb{R}^2,0)
\end{equation}

satisfying \(F(0)=0\). After complexifying the tangent space of the target, the expressions for \(\nabla_{\partial_z}s\) and \(\nabla_{\partial_{\overline z}}s\) give 
\begin{equation}\label{eq: differential dev map proof} 
    dF(\partial_z) = z^{m_i}u(z), \quad dF(\partial_{\overline z}) = \overline z^{\,m_i}\overline{u(z)}, 
\end{equation} 

where \(u(z)\) is the image, under the differential of the projectivisation map, of \(\gamma_0(z)\ell^{-1}\). Since \(\gamma_0(0)\neq 0\) and the vector \(\ell^{-1}\) is not contained in the line generated by \(s\), we have \(u(0)\neq 0\). The second equality in \eqref{eq: differential dev map proof} follows also from the fact that \(F\) is a real map. Furthermore, \(u(0)\) and \(\overline{u(0)}\) are linearly independent in the complexified tangent space of the target. Integrating along the radial segment from $0$ to $z$, we obtain 
\begin{equation} 
    F(z)=\frac{1}{m_i+1}\left(z^{m_i+1}u(0)+\bar z^{m_i+1}\overline{u(0)}\right)+O(|z|^{m_i+2}).
\end{equation} 

If we now consider the real-linear map 
\begin{equation} 
    A\colon\mathbb C\longrightarrow\mathbb R^2,\quad A(w)=wu(0)+\bar w\,\overline{u(0)},
\end{equation} 
using the linear independence of $u(0)$ and $\overline{u(0)}$ we get that $A$ is an isomorphism. Consequently, we have that  \begin{equation} 
    F(z)=\frac{1}{m_i+1}A(z^{m_i+1})+O(|z|^{m_i+2}).
\end{equation}
In particular, $0$ is an isolated preimage of $F(0)$ and the absolute value of the local topological degree of $F$ at $0$ is $m_i+1$. Since $F$ is a local diffeomorphism away from $0$, its restriction to sufficiently
small punctured neighbourhoods defines a finite covering of degree $m_i+1$.
Hence $\widetilde p_i$ is a branch point of order $m_i$, as desired.
\end{proof}

As a direct consequence of Lemma~\ref{lem: branche structures from Higgs bundles} and Proposition~\ref{prop:BPSasmaps}, we readily conclude that every representation in either of the spaces \(\mathcal{R}_d^+(X)\) and \(\mathcal{R}_d^-(X)\) is indeed the holonomy of a branched projective structure whose branching data is precisely determined by the divisor associated with the holomorphic section \(\gamma\). This generalises a result of Baraglia \cite[\S 3.4.2]{baraglia2009g2} to a wider class of Higgs bundles outside the Hitchin component. In this broader context, the next objective is to construct a Riemannian metric with conical singularities along \(Z(\gamma)\) arising from a conical hyperbolic affine sphere, thereby extending the classical construction of the Blaschke metric on convex domains.

\smallskip

Before proceeding, we need to introduce a torsion-free connection on the punctured surface \(X_0\) as follows. Let \(V \subset E_{\mathbb{R}}\) denote the real rank-\(2\) sub-bundle defined as the fixed locus of \(\tau\) in \(L \oplus L^{-1}\). Notice that \(V\) is naturally identified with the quotient bundle \(E_{\mathbb{R}}/\underline{\mathbb{R}}\). The construction above then yields an isomorphism
\begin{equation}\label{eq: singular affine connection}
    \begin{aligned}
        T_pX_0 &\longrightarrow V_p|_{ X_0}, \\
        Y &\longmapsto \nabla_Y s \quad(\text{mod} \ \underline{\mathbb R}),
    \end{aligned}
\end{equation}

for every \(p\in X_0\defeq  X\setminus Z(\gamma)\). The aforementioned torsion-free connection, denoted with \(\widehat\nabla\), on \(TX_0\), is completely determined by the identity
\begin{equation}\label{eq: affine connection Higgs bundles}
    \nabla_Y\nabla_U s=\nabla_{\widehat\nabla_Y U}s+G_B(Y,U)s,\quad Y,U\in\Gamma(T X_0)
\end{equation}

where \(G_B\) is a bilinear form on \(T X_0\). Notice that the flatness of \(\nabla\) translates into the torsion-freeness of \(\widehat\nabla\) and the symmetry of the bilinear form \(G_B\). As we shall see, the objects \(\widehat\nabla\) and \(G_B\) are in fact, respectively, the affine connection and the Blaschke metric of a conical hyperbolic affine sphere.
\begin{prop}\label{thm: developing map lifts to a singular affine spheres}
The developing map $\operatorname{dev}_\sigma:\widetilde X\to\RP^2$, defined by the section \(\sigma\), lifts to a $\rho$-equivariant conical hyperbolic affine sphere $ f:\widetilde X\longrightarrow\mathbb R^3 $ whose singular locus is precisely the lift of $Z(\gamma)$ to the universal cover of the surface. The conical Blaschke metric and the Pick differential of this lift are given by \begin{equation}\label{eq: Blaschke metric e Pick differential developing map}
G_B=2h^{-1}|\gamma|^2,\quad q_{Pick}=-2\nu\gamma^2. 
\end{equation}In particular, both are uniquely determined by the isomorphism class of the stable Higgs bundle.
\end{prop}
\begin{proof}
Arguing as in the construction of $\operatorname{dev}_\sigma$ before Lemma~\ref{lem: branche structures from Higgs bundles}, but working with $E_{\mathbb R}$ instead of its projectivization, the section $s$ determines a smooth $\rho$-equivariant map \(f_{s}:\widetilde X\longrightarrow\mathbb R^3\) whose projectivization coincides with $\operatorname{dev}_\sigma$. For simplicity, we shall henceforth denote $f_{s}$ by $f$. Moreover, by \(\rho\)-equivariance, we may work directly on the surface \(X\). In accordance with Definition \ref{def:singular affine sphere}, we now have to verify that \(f\) satisfies the defining properties of a singular hyperbolic affine sphere. We begin by showing that \(f|_{X_0}\) defines a hyperbolic affine sphere in the classical sense, where \(X_0=X\setminus Z(\gamma)\). In the process, we will see that the associated Blaschke metric has conical singularities precisely along the zeroes of \(\gamma\). We first observe that, by construction, \(f_*(Y)=\nabla_Y s\) for every vector field \(Y\) on \(X_0\).
Now let $(\ell,e,\ell^{-1})$ be a local holomorphic frame adapted to the decomposition \(\mathcal E=L\oplus\mathcal O_X\oplus L^{-1}.\) Writing $\gamma=\gamma_0 \ \ell^{-1}\otimes dz$ and $\nu=\nu_0 \ \ell\otimes \ell\otimes dz$ in the corresponding local trivializations, the flat connection gives
\[
f_z=\nabla_{\partial_z}s=\gamma_0\ell^{-1},
\qquad
f_{\bar z}=\nabla_{\partial_{\bar z}}s
=h_0^{-1}\overline{\gamma_0}\ell,
\] being \(h_0\) the smooth positive local function representing \(h\). A direct computation using the flat connection as in Lemma \ref{lem: branche structures from Higgs bundles}, yields
\( \nabla_{\partial_z}f_{\bar z}=\nabla_{\partial_{\bar z}}f_z=\lambda f,\)
where \(\lambda=h_0^{-1}|\gamma_0|^2\). On the other hand, we have 
\begin{align*}
&\nabla_{\partial_z}f_z=(\partial_z\log\lambda)f_z+\frac{h_0\nu_0\gamma_0}{\overline{\gamma_0}}f_{\bar z},\\ & \nabla_{\partial_{\bar z}}f_{\bar z}=\frac{h_0\overline{\nu_0}\,\overline{\gamma_0}}{\gamma_0}f_z+(\partial_{\bar z}\log\lambda)f_{\bar z}.
\end{align*}
Therefore, the the bilinear form \(G_B\) in Equation \eqref{eq: affine connection Higgs bundles} is completely determined by
\[
G_B(\partial_z,\partial_z)=G_B(\partial_{\bar z},\partial_{\bar z})=0,\qquad G_B(\partial_z,\partial_{\bar z})=\lambda.
\]
Since \(|dz|^2=dx^2+dy^2,\) this is equivalent to \(G_B=2h_0^{-1}|\gamma_0|^2|dz|^2.\) Moreover, by construction we have \(\nabla_Yf=f_*(Y),\) which implies that \(f\) satisfies Equations \eqref{eq: structure equations hyperbolic affine spheres}. It only remains to check the Blaschke normalization \(\theta_f=d A_{G_B}\). Let \(\varepsilon_{\mathbb C}=\ell^*\wedge e^*\wedge (\ell^{-1})^*\) be the induced holomorphic volume form on \(\mathcal E\). The $\tau$-fixed $H$-orthonormal frame
\[
u_1=\frac{1}{\sqrt2}\left(h^{-1/2}_0\ell+h_0^{1/2}\ell^{-1}\right),\qquad u_2=e,\qquad u_3=
\frac{i}{\sqrt2}\left(h_0^{-1/2}\ell-h_0^{1/2}\ell^{-1}\right)
\]
satisfies \(\varepsilon_{\mathbb C}(u_1,u_2,u_3)=-i.\) Choosing the orientation for which $(u_1,u_2,u_3)$ is positive, the parallel real volume form on $E_{\mathbb R}$ is therefore \(\Omega_{\mathbb R}=\left.i\,\varepsilon_{\mathbb C}\right|_{E_{\mathbb R}}.\) Since
\[
f_x=h_0^{-1}\overline{\gamma_0}\ell+\gamma_0\ell^{-1},\qquad f_y=-ih_0^{-1}\overline{\gamma_0}\ell+i\gamma_0\ell^{-1},
\]
we obtain
\[
\Omega_{\mathbb R}(f_x,f_y,f)=i\det\begin{pmatrix}h_0^{-1}\overline{\gamma_0}&
-i h_0^{-1}\overline{\gamma_0}&0\\0&0&1\\\gamma_0&i\gamma_0&0\end{pmatrix}=2h_0^{-1}|\gamma_0|^2.
\]
On the other hand, \(dA_{G_B}=2h_0^{-1}|\gamma_0|^2\,dx\wedge dy\) and thus \(\theta_f=\Omega_{\R}(\cdot,\cdot,f)=dA_{G_B}.\) Therefore $f|_{X_0}$ is a hyperbolic affine sphere centred at the origin, with affine shape operator $B=-\operatorname{Id}$ and Blaschke metric \(G_B=2h_0^{-1}|\gamma_0|^2|dz|^2.\) In particular, \(f\) is well defined on the entire surface \(X\) and satisfies all the conditions in Definition \ref{def:singular affine sphere}. Finally, we identify the Pick differential of the affine sphere. Recall that, with our convention, the Pick tensor is \(C\defeq \widehat\nabla G_B,\) and the Pick differential $q_{Pick}$ is determined by \(C=2\operatorname{Re}(q_{Pick}),\) or, equivalently, \(q_{Pick}=C^{3,0}=C(\partial_z,\partial_z,\partial_z).\)
From the computation above the induced affine connection satisfies
\[
\widehat\nabla_{\partial_z}\partial_z=(\partial_z\log\lambda)\partial_z+\frac{h_0\nu_0\gamma_0}{\overline{\gamma_0}}\partial_{\bar z},\qquad\widehat\nabla_{\partial_z}\partial_{\bar z}=0.
\]
Therefore we obtain
\begin{align*}
C(\partial_z,\partial_z,\partial_z)&=(\widehat\nabla_{\partial_z}G_B)(\partial_z,\partial_z)\\&=-2G_B\left(\widehat\nabla_{\partial_z}\partial_z,\partial_z\right)\\&=-2\frac{h_0\nu_0\gamma_0}{\overline{\gamma_0}}G_B(\partial_{\bar z},\partial_z)\\
&=-2\frac{h_0\nu_0\gamma_0}{\overline{\gamma_0}}h_0^{-1}|\gamma_0|^2\\&=-2\nu_0\gamma_0^2,
\end{align*} and it follows that the Pick differential is \(q_{Pick}=-2\nu\gamma^2\).
\end{proof}

\smallskip 

\begin{rmk}
    Surprisingly, upon fixing a square root of the canonical bundle, we may observe that the same construction works also when \(\deg(L)=1-g\) and hence \(L\cong K^{-1/2}\). In this case, the Higgs bundle is given by
    \begin{equation}
        \mathcal{E}=K^{-1/2}\oplus\mathcal O_X\oplus K^{1/2},\qquad\Phi=\begin{pmatrix}0 & 0 & 1 \\ \gamma & 0 & 0 \\0 & \gamma & 0\end{pmatrix},
    \end{equation}
where \(\gamma\neq 0\) is a holomorphic section of \(K^{3/2}\). This specific case is not new in literature as it has been already studied by Bronstein--Davalo \cite[\S 6.2]{bronstein2025anosov} where the authors construct a fibred geometric structure, without branching, modelled on the complete flag variety in \(\R^3\). We observe that the expression of the Higgs bundle used in \cite{bronstein2025anosov} differs from ours, but the two are isomorphic by the gauge transformation with entries equal to \(-1\) on the anti-diagonal and \(0\) elsewhere, thus representing the same point in the moduli space, with the only difference being that in the present work the assumption \(\gamma\neq 0\) is necessary for the stability to be guaranteed. 
\end{rmk}

\subsection{Prescribing branching orders and the total branching degree}\label{sec: prescribed data problem}

As discussed in the previous section, upon fixing a Riemann surface structure, it is possible to realise a branched projective structure of hyperbolic affine sphere type on \(X\) by starting from a certain family of Higgs bundles, see Definition \ref{def: branched conical projective}. On the other hand, in the study of branched hyperbolic structures, remarkable results of Troyanov and McOwen, see \cite{troyanov1991prescribing, mcowen1993prescribed} respectively, show that, upon fixing a conformal structure on \(S\) and a set of marked points, say \(\mathcal P\) on \(S\), there exists a unique hyperbolic metric with prescribed conical singularities at \(\mathcal P\), provided that the cone angles satisfy a Gauss--Bonnet type condition. The aim of this section is to show that an analogous phenomenon also occurs for the projective setting, provided that the degree of the branching data satisfies a bound from above and each multiplicity is an integer. 

\begin{thm}\label{thm: prescribed problem data}
    Let \(D=\sum_{i=1}^km_ip_i\) be a \(0\)-chain on \(S\) such that \(\deg(D)\in(0,2g-2)\). Then, for any Riemann surface structure \(X\) on \(S\), there exists a family of branched projective structures of hyperbolic affine sphere type on \(S\) with branching divisor \(D\), parametrised by \(H^0\bigl(X,K^3(-2D)\bigr).\) Moreover, this family form a moduli space of complex dimension \(5g-5-2\deg(D)\).
\end{thm}

A similar statement holds for zero-chains of degree \(\deg(D)\in(2g-2,3g-3]\). More specifically, the following holds.

\begin{thm}\label{thm: prescribed problem data2}
    Let \(D=\sum_{i=1}^km_ip_i\) be a \(0\)-chain on \(S\) such that \(\deg(D)\in(2g-2,3g-3]\). Let \(X\) be any Riemann surface structure on \(S\) such that \(H^0\big(X,K^3(-2D)\big)\neq0\) holds. Then there exists a family of branched projective structures of hyperbolic affine sphere type on \(S\) with branching divisor \(D\), parametrised by \(H^0\bigl(X,K^3(-2D)\bigr)\setminus\{0\}\).
\end{thm}

\begin{proof}[Proof of Theorem \ref{thm: prescribed problem data}]
    Let \(X\) be a Riemann surface structure on \(S\), and define the line bundle \(L\defeq K(-D),\) of degree \(d\defeq 2g-2-\deg(D)\in(0,2g-2).\) Since \(L^{-1}K\cong\mathcal O_X(D)\), we denote by \(\gamma\neq 0\) the canonical holomorphic section of \(\mathcal O_X(D)\), whose zero divisor is precisely \(D\). As a consequence of the Riemann--Roch theorem, a direct computation readily implies that \(\dim_{\mathbb C} H^0\big(X,K^3(-2D)\big)=5g-5-2\deg(D)>0\). For every \(\nu\in H^0\bigl(X,L^2K\bigr)\cong H^0\bigl(X,K^3(-2D)\bigr),\) we consider the stable \(\mathrm{SL}(3,\mathbb R)\)-Higgs bundle defined by the data
    \begin{equation}\label{eq: fibrato di Higgs stabile nella dimostrazione}
        \mathcal{E}=L\oplus\mathcal O_X\oplus L^{-1},\qquad\Phi=\begin{pmatrix}0&0&\nu\\
        \gamma&0&0\\0&\gamma&0\end{pmatrix}.
    \end{equation}
    By Lemma \ref{lem: branche structures from Higgs bundles} and Proposition \ref{thm: developing map lifts to a singular affine spheres}, each such Higgs bundle gives rise to a branched projective structure of hyperbolic affine sphere type \((\operatorname{dev}_\nu,\rho_\nu)\) whose branching divisor is \(\operatorname{div}(\gamma)=D\). We may observe that the vector space \(H^0\bigl(X,K^3(-2D)\bigr)\) is naturally identified with the fibre over \(D\) of the projection
    \begin{equation}
        \mathcal{QD}_d(X)\longrightarrow \operatorname{Sym}^{2g-2-d}(X),\quad(q_3,D')\longmapsto D',
    \end{equation}
    already considered in the proof of Proposition \ref{prop: parametrization of subspace higgs bundles}. Through the biholomorphism \(\Theta_d:\mathcal{QD}_d(X)\to\mathcal{CS}_d(X)\), this fibre thus corresponds to the family of branched projective structures of hyperbolic affine sphere type in the slice \(\mathcal{CS}_d(X)\) with fixed branching divisor \(D\). 
\end{proof}   
    
\begin{proof}[Proof of Theorem \ref{thm: prescribed problem data2}]
    Let \(X\) be a Riemann surface structure on \(S\) such that \(H^0\big(X,K^3(-2D)\big)\neq\{0\}.\) Once again, we define the line bundle \(L\defeq K(-D),\) which now has negative degree \(d\in[1-g,0).\) By proceeding as in the proof of Theorem \ref{thm: prescribed problem data}, for every non-zero holomorphic section \(\nu\in H^0\big(X,K^3(-2D)\big)\), we can realise the stable \(\mathrm{SL}(3,\mathbb R)\)-Higgs bundle as in equation \eqref{eq: fibrato di Higgs stabile nella dimostrazione}. By applying Lemma \ref{lem: branche structures from Higgs bundles} and Proposition \ref{thm: developing map lifts to a singular affine spheres} once again, we obtain the desired conclusion.
\end{proof}

For clarity, we point out that the assumption \(H^0\big(X,K^3(-2D)\big)\neq\{0\}\) when \(\deg(D)\in(2g-2,3g-3]\) has been imposed to ensure the stability of the Higgs bundle, as proved in Proposition \ref{prop: stable SL(3;C) Higgs bundle}. In particular, in the extremal case \(\deg(D)=3g-3,\) we obtain an isomorphism \(K^3(-2D)\cong\mathcal O_X,\) and therefore \(\mathcal O_X(D)\) must be a square root of \(K^3\). In this case, \(\nu\) is a complex constant.   

\begin{rmk}
    We may notice that Theorems \ref{thm: prescribed problem data} and \ref{thm: prescribed problem data2} are sharp. In fact, weaker versions require only the total branching degree to be prescribed, \textit{i.e.}, the sum of the local branching orders rather than the full branching divisor. In fact, a similar argument shows that for every \textit{non-zero} integer $d \in [1-g,2g-2)$ and every fixed Riemann surface $X$, there exists a family of branched projective structures of hyperbolic affine sphere type with total branching degree equal to $2g-2-d$. Each of these families is parametrised by either the space $\mathcal{CS}_d^+(X)$ or $\mathcal{CS}_d^-(X)$, where the sign is taken according to the value of $d$ as explained in \S\ref{ssec: cyclic higgs bundles}, see Equation \eqref{eq: subspace of cyclic higgs bundles}. Each of these families has complex dimension \(5g-5-\deg(D)\), where $D=\operatorname{div}(\gamma)$ is the branching divisor determined by the holomorphic section $\gamma$, and $D$ is allowed to vary in the corresponding base of the holomorphic vector bundle.
\end{rmk}

\section{Analysis of branched projective structures of hyperbolic affine sphere type}\label{sec:analysis_and_geometry}

Building upon the geometric realisation established earlier, this chapter is devoted to the analytical analysis of branched projective structures of hyperbolic affine sphere type, with a primary focus on the Hitchin equations and the properties of the conical Blaschke metrics.

\smallskip

\subsection{Analysis of the Hitchin equations}\label{sec: analysis Hitchin equations}
Having established how a certain class of stable \(\mathrm{SL}(3,\mathbb R)\)-Higgs bundles gives rise to branched geometric structures on \(X\), we now turn to the analytic study of the Hitchin equations. We recall that every such a structure determines a holonomy representation, namely an element of \(\mathfrak{X}_3(S)\). When \(\deg(L)=2g-2\), then the corresponding representation is a Hitchin representation, see the discussion after Remark \ref{rem: cyclic_higgs_bundles}, and \(L\cong K\). In this case, the Hermitian metric \(2h^{-1}\) on \(K^{-1}\) has strictly negative Gaussian curvature and coincides with the Blaschke metric of the corresponding hyperbolic affine sphere, see \cite[\S 3.4.2]{baraglia2009g2}. In this section, we show that the same conclusion holds in the more general framework where \(\deg(L)=d\in(0,2g-2)\), provided that the divisors of \(\gamma\) and \(\nu\) satisfy a suitable condition. We also extend our analysis to negative values of \(d\), \textit{i.e.}, \(d\in[1-g,0)\), where the scenario differs significantly; see Corollary \ref{cor: Curvature Blaschke metric singular case}.

\smallskip

\subsubsection{General analysis} In what follows, we denote by \([L,\gamma,\nu]\) the isomorphism class of a Higgs bundle of the form \eqref{eq: cyclic Higgs bundles}, lying in either \(\mathcal{C}_d^+(X)\) or \(\mathcal{C}_d^-(X)\), according to the degree of the line bundle \(L\). Since the Higgs bundle is stable, it follows that \(H=\operatorname{diag}(h,1,h^{-1})\) solves the Hitchin equation \eqref{eq: Hitchin_higgs bundles}. In the present setting, this reduces to a single equation for \(h\), which can be written locally as
\begin{equation}\label{eq: local equation for h}
    \partial_z\partial_{\overline z}\log h_0+|\gamma_0|^2h_0^{-1}-|\nu_0|^2h_0^2=0,
\end{equation}

where \(h_0\), \(\gamma_0\), and \(\nu_0\) denote the local expressions of \(h\), \(\gamma\), and \(\nu\), respectively, in a holomorphic chart \((U,z)\). Let \(g_0\) denote the unique hyperbolic metric in the conformal class of \(X\). By multiplying both sides of Equation \eqref{eq: local equation for h} by \(g_0^{-1}\), we obtain globally defined equation
\begin{equation}\label{eq: global Hitchin equation}
\frac{1}{4}\Delta_0\log h_0+\|\gamma\|^2-\|\nu\|^2=0,
\end{equation}
where \(\Delta_0\) denotes the Laplace operator with respect to \(g_0\). Although \(h_0\) denotes the local expression of the Hermitian metric \(h\), the quantity \(\Delta_0\log h_0\) is globally defined, since two local expressions of \(h\) differ by the squared modulus of a nowhere-vanishing holomorphic function. Before turning to the main result, we state a technical lemma that will be useful later, see \cite[Proposition 3.2]{diaf2025domination} for a similar setting.

\begin{lem}\label{lem: divisor condition norms}
    The following statements hold: 
    \begin{itemize}
        \item[(i)] in the case \(d\in(0,2g-2]\), if \(\text{div}(\gamma)<\text{div}(\nu)\) then \(||\nu||^2<||\gamma||^2\) on the punctured surface \(X_0=X\setminus Z(\gamma)\);
        \item[(ii)] in the case \(d\in[1-g,0)\), if \(\text{div}(\nu)<\text{div}(\gamma)\) then \(||\gamma||^2<||\nu||^2\) on the punctured surface \(X_0=X\setminus Z(\nu)\).
    \end{itemize}
\end{lem}

\begin{proof}
    We prove item (i), since item (ii) follows by exchanging the roles of \(\gamma\) and \(\nu\). Set \( D\defeq \operatorname{div}(\nu)-\operatorname{div}(\gamma)\) so that, by assumption, \(D\) is a non-zero effective divisor. Consider the function
    \begin{equation}\label{eq: aux function}
        u\defeq \log\frac{\|\nu\|^2}{\|\gamma\|^2}.
    \end{equation}
    Although the auxiliary function \eqref{eq: aux function} \(u\) is defined away from the zeros of \(\gamma\) and \(\nu\), it extends smoothly across every common zero at which \(\gamma\) and \(\nu\) have the same order. In fact, if \(p\in X\) and, in a local coordinate centred at \(p\), \(\gamma_0=z^n\tilde{\gamma}_0, \  \nu_0=z^m\tilde{\nu}_0,\) with \(\tilde{\gamma}_0(0),\tilde{\nu}_0(0)\neq0\) and \(m\ge n\), then locally
    \begin{equation}
    u=2(m-n)\log|z|+\log\left|\frac{\tilde{\nu}_0}{\tilde{\gamma}_0}\right|^2+3\log h_0.
    \end{equation}
    Hence \(u\) extends smoothly when \(n=m\), whereas \(u\to-\infty\) as one approaches \(p\) whenever \(m>n\). It follows that \(u\) is smooth on \(X\setminus\operatorname{Supp}(D)\) and tends to \(-\infty\) near \(\operatorname{Supp}(D)\). On \(X\setminus\operatorname{Supp}(D)\), the function \(u\) satisfies
    \begin{equation}\label{eq: differential equation u}
        \Delta_0u=\frac{4}{g_0}\partial_z\partial_{\bar z}\bigl(\log|\nu_0|^2+3\log h_0-\log|\gamma_0|^2\bigr)=3\Delta_0\log h_0=12\bigl(\|\nu\|^2-\|\gamma\|^2\bigr)=12\|\gamma\|^2(e^u-1),
    \end{equation}
    where we have used Equation \eqref{eq: global Hitchin equation}. We first show that \(u\le0\). Suppose, by contradiction, that \(u>0\) somewhere, and let \(\Omega\) be a connected component of \( \{x\in X\setminus\operatorname{Supp}(D):u(x)>0\}\). Since \(u\to-\infty\) near \(\operatorname{Supp}(D)\), the closure \(\overline{\Omega}\) is a compact subset of \(X\setminus\operatorname{Supp}(D)\) and \(u=0 \ \text{on} \ \partial\Omega.\) Instead, on \(\Omega\), Equation \eqref{eq: differential equation u} gives
    \begin{equation}
        \Delta_0u=12\|\gamma\|^2(e^u-1)\ge0.
    \end{equation}
    Thus \(u\) is sub-harmonic on \(\Omega\), and the weak maximum principle yields
    \begin{equation}
        \sup_{\Omega}u\le\sup_{\partial\Omega}u=0,
    \end{equation}
    that contradicts \(u>0\) on \(\Omega\). Therefore \(u\le0\) on \(X\setminus\operatorname{Supp}(D)\). We next prove that the inequality is strict. Since \(e^u-1\ge u\), Equation \eqref{eq: differential equation u} implies \(\Delta_0u\ge12\|\gamma\|^2u.\) Since \(u\le0\), if \(C\defeq 12\max_X\|\gamma\|^2>0,\) then \(12\|\gamma\|^2u\ge Cu,\)     and therefore
    \begin{equation}
        (\Delta_0-C)u\ge0 \qquad\text{on }X\setminus\operatorname{Supp}(D).
    \end{equation}
    Since \(X\setminus\operatorname{Supp}(D)\) is connected, the strong maximum principle, see \cite[Theorem 1]{minda1987strong}, implies that either 
    \begin{itemize}
        \item[1.] \(u<0\) everywhere on \(X\setminus\operatorname{Supp}(D)\), or
        \item[2.] \(u\equiv0\). 
    \end{itemize}
    The latter is impossible, since \(D\neq0\) and \(u\to-\infty\) near \(\operatorname{Supp}(D)\). Consequently, \(u<0 \ \text{on} \ X\setminus\operatorname{Supp}(D).\) It remains to consider points of \(X_0\cap\operatorname{Supp}(D)\). If \(p\) is such a point, then \(\gamma(p)\neq0\), whereas \(\nu(p)=0\). Thus the inequality \(\|\nu\|^2(p)=0<\|\gamma\|^2(p)\) is trivial. Combining the two cases gives \(\|\nu\|^2<\|\gamma\|^2 \ \text{on} \ X_0,\) as claimed.
\end{proof}

In fact, the preceding argument establishes the stronger inequality \(u<0\) on \(X\setminus\operatorname{Supp}(D)\), including the common zeros of \(\gamma\) and \(\nu\) having the same vanishing order.

\begin{rmk}
    Similar, though not necessarily identical, conditions on the divisors of \(\gamma\) and \(\nu\) have been considered in other contexts. In \cite{dai2019cyclic}, and more generally in \cite{sagman2024hitchin}, it is used in the study of maximum principles for systems of PDEs arising from the Hitchin equations for cyclic Higgs bundles, and in the study of the curvature of the metric on the surface induced by the harmonic map. In \cite{biswas2021uniformization, m2024conformal}, it is employed in the study of a special class of \(\mathrm{SL}(2,\mathbb R)\)-Higgs bundles in relation to branched hyperbolic structures and conformal limits. In \cite{bronstein2025anosov, zhang2025non} such a condition is used to prove the Anosov property for the corresponding holonomy representation. In \cite{collier2024holomorphic, nie2024cyclic, rungi2025complex}, analogous assumptions are imposed in the study of infinitesimally rigid minimal surfaces in homogeneous spaces. Finally, in \cite{evans2025geometric}, similar conditions are used to construct geometric structures associated with \(\mathrm{G}_2'\)-Hitchin representations.
\end{rmk}

In the following we shall denote by \(2h^{-1}|\gamma|^2\) the Riemannian metric with conical singularities constructed in Proposition \ref{thm: developing map lifts to a singular affine spheres} from the Higgs bundle \([L,\gamma,\nu]\), see also Equation \eqref{eq: affine connection Higgs bundles}.

\begin{thm}\label{thm: gaussian curvatura singular metrics}
The following conditions hold:
\begin{itemize}
    \item[(i)] If \(d\in (0,2g-2]\) and \(\text{div}(\gamma)<\text{div}(\nu)\), then the singular metric \(2h^{-1}|\gamma|^2\) has strictly negative Gaussian curvature on \(X\setminus Z(\gamma)\). Moreover, the curvature \(K_{2h^{-1}|\gamma|^2}\) can be extended to a strictly negative continuous function on \(X\).
    \item[(ii)] If \(d\in[1-g,0)\) and \(\text{div}(\nu)<\text{div}(\gamma)\), then the singular metric \(2h^{-1}|\gamma|^2\) has strictly positive Gaussian curvature on \(X\setminus Z(\gamma)\). Moreover, the curvature \(K_{2h^{-1}|\gamma|^2}\) goes to \(+\infty\) approaching either the (possibly empty) finite subset \(Z(\gamma)\setminus Z(\nu)\) or a point \(p\in Z(\nu)\) such that \(\operatorname{ord}_p(\nu)<\operatorname{ord}_p(\gamma)\).
\end{itemize}
\end{thm}
\begin{proof}
We prove point (i), corresponding to the case \(d\in(0,2g-2].\) According to our convention for the Laplacian, \(\Delta_0=\frac{4}{g_0}\partial_z\partial_{\bar z},\) the Gaussian curvature of a Riemannian metric \(G\) defined on \(TX\) is given by 
\begin{equation}
    K_G=-\frac{2}{G_0}\partial_z\partial_{\bar z}\log G_0,
\end{equation}

where \(G_0\) denotes the local positive conformal factor representing \(G\). On \(X_0\defeq X\setminus Z(\gamma)\), there is a well-defined Riemannian metric \(2h^{-1}|\gamma|^2\) on \(TX_0\), whose local expression is given by \(2h_0^{-1}|\gamma_0|^2|\mathrm dz|^2\). A straightforward computation then yields   
\begin{align*}
    K_{2h^{-1}|\gamma|^2}&=-\frac{1}{h_0^{-1}|\gamma_0|^2}\partial_z\partial_{\bar z}\log(h^{-1}_0|\gamma_0|^2)\\&=\frac{1}{h_0^{-1}|\gamma_0|^2}\partial_z\partial_{\bar z}\log h_0-\frac{1}{h_0^{-1}|\gamma_0|^2}\partial_z\partial_{\bar z}\log|\gamma_0|^2 \\ & =\frac{g
    _0}{4h_0^{-1}|\gamma_0|^2}\Delta_0\log h_0 \tag{\(\gamma_0\) is holomorphic} \\ &=\frac{1}{\|\gamma\|^2}(\|\nu\|^2-\|\gamma\|^2) \tag{Equation \eqref{eq: global Hitchin equation}} \\ &=\frac{\|\nu\|^2}{\|\gamma\|^2}-1.
\end{align*} 

The conclusion follows from Lemma \ref{lem: divisor condition norms}. Let \(p\in Z(\gamma)\). Upon fixing some local coordinate \(z\), both sections \(\gamma\) and \(\nu\) can be locally written as \(\gamma=z^n\tilde\gamma_0(z) l^{-1}\otimes \mathrm d z\) and \(\nu=z^m \tilde\nu_0(z)l^2\otimes\mathrm dz\), with \(\tilde\gamma_0,\tilde\nu_0\) local holomorphic function such that \(\tilde\gamma_0(0),\tilde\nu_0(0)\neq 0\). By hypothesis on the divisor we know that \(n\le m\) and, if the inequality is strict, then clearly \(K_{2h^{-1}|\gamma|^2}(p)=-1\). If \(n=m\), then the curvature admits the finite limit 
\begin{equation}
    h_0(p)^3\frac{|\tilde{\nu}_0(p)|^2}{|\tilde{\gamma}_0(p)|^2}-1.
\end{equation}

To see that this limit is strictly negative, recall that the proof of Lemma \ref{lem: divisor condition norms} establishes 
\begin{equation}
    u\defeq \log\frac{\|\nu\|^2}{\|\gamma\|^2}<0 
\end{equation}

on $X\setminus\operatorname{Supp}(D)$, where $D=\operatorname{div}(\nu) -\operatorname{div}(\gamma)$. Since $m=n$, we have $p\notin\operatorname{Supp}(D)$, and hence 
\begin{equation}
    u(p)=\log\left(h_0(p)^3\frac{|\tilde{\nu}_0(p)|^2}{|\tilde{\gamma}_0(p)|^2}\right)<0.
\end{equation}

Therefore, the extended curvature is strictly negative at $p$. Consequently, the Gaussian curvature extends continuously to a function on $X$ taking values in $[-1,0)$. The statement (ii) easily follows. As for the first part, the same procedure applies \textit{mutatis mutandis} by using the expression of the curvature and Lemma \ref{lem: divisor condition norms}. Instead, if there exists \(p\in Z(\gamma)\setminus Z(\nu)\) then clearly \(K_{2h^{-1}|\gamma|^2}(p_n)\to+\infty\) as \(p_n\to p\), and the same happens for the other case.
\end{proof}

For \(d=2g-2\), and hence \(L\cong K\), the result was already known. Indeed, the Riemannian metric \(2h^{-1}\) is well defined on all of \(TX\) and coincides with the Blaschke metric of the associated hyperbolic affine sphere, see \cite[\S 6.2]{li2019introduction}. In this case, the condition on the divisors appearing in Lemma \ref{lem: divisor condition norms} is automatically satisfied, since \(\gamma=1\in\mathbb C.\) The same happens when \(d=1-g\) since \(\nu=1\in\C\). 

\begin{rmk}\label{rem: curvature of induced metric from Higgs bundles}
    It is interesting to observe that, in the range \(d\in[1-g,0)\), the metric arising from the construction of the branched \(\mathbb{RP}^2\)-structure via the developing map described in \S \ref{sec: from cyclic to geometric structres} has positive Gaussian curvature under the assumptions of Lemma \ref{lem: divisor condition norms}, and it always diverges in a finite set of points. Indeed, a point \(p\in Z(\nu)\) such that \(\operatorname{ord}_p(\nu)<\operatorname{ord}_p(\gamma)\) or a point \(p\in Z(\gamma)\setminus Z(\nu)\), always exists according to the hypothesis on the divisors. In this case, it is possible to consider the singular metric \(h^2|\nu|^2,\) which has strictly negative Gaussian curvature on \(X\setminus Z(\nu)\) and whose singular locus is precisely \(Z(\nu)\). However, this metric does not arise from the developing-map construction of the branched \(\mathbb{RP}^2\)-structure, whose naturally induced metric is \(2h^{-1}|\gamma|^2\). This distinction will play an important role in the next sections, where we discuss the relation with branched hyperbolic structures.
\end{rmk}

\begin{cor}\label{cor: boundendess cubic differential}
    Let \([L, \gamma, \nu]\) be a Higgs bundle of the form \eqref{eq: cyclic Higgs bundles}, and let \(q_3 = \gamma^2 \nu\) be the holomorphic cubic differential associated with \([L, \gamma, \nu]\). Then the following hold:
    \begin{itemize}
        \item[(i)] If \(d\in (0,2g-2]\) and \(\text{div}(\gamma)<\text{div}(\nu)\), then $\|q_3\|_{h^{-1}|\gamma|^2}^2$ extends smoothly to a function on $X$ which is strictly smaller than $1$ everywhere;
        \item[(ii)] If \(d\in[1-g,0)\) and \(\text{div}(\nu)<\text{div}(\gamma)\), then \(\|q_3\|^2_{h^{-1}|\gamma|^2}>1\) on \(X_0\). Moreover, this quantity tends to \(+\infty\) as one approaches any point of the possibly empty finite set \(Z(\gamma)\setminus Z(\nu)\), or any point \(p\in Z(\nu)\) satisfying \(\operatorname{ord}_p(\nu)<\operatorname{ord}_p(\gamma). \)
    \end{itemize}
\end{cor}

\begin{proof}
    The claim follows directly from the expression of the curvature \(K_{2h^{-1}|\gamma|^2}=\frac{\|\nu\|^2}{\|\gamma\|^2}-1\), Theorem \ref{thm: gaussian curvatura singular metrics} and the computation below: 
    \begin{equation}
        \|q_3\|^2_{h^{-1}|\gamma|^2}=|\nu_0|^2|\gamma_0|^4(h^{-1}_0|\gamma|^2)^{-3}=|\nu_0|^2|\gamma_0|^{-2}h^{3}_0=\frac{\|\nu\|^2}{\|\gamma\|^2}=\|s\|_{h^3}^2,
        \end{equation}
being \(s\defeq \displaystyle\frac{\nu}{\gamma}\in H^0(X,L^3)\).
\end{proof}

Once again, when \(d=2g-2\) the inequality was already known and proven to be always strict, see for example \cite[\S 6.2]{li2019introduction} for more details. By contrast, the case \(d\in[1-g,0)\) exhibits a substantial difference from other settings, in which boundedness of the norm of the \(k\)-differential is typically expected, see \cite{dai2024bounded} for several results in this direction.

\subsubsection{Analysis of conical hyperbolic affine spheres}
According to Theorem \ref{thm: every rho is holonomy of conical branched}, we recall that to every Higgs bundle \([L,\gamma,\nu]\) with \(d\in[1-g,2g-2)\setminus\{0\}\), we may associate its developing-holonomy pair \((\operatorname{dev},\rho)\) defining the corresponding branched projective structure of hyperbolic affine sphere type. In particular, there is a \(\rho\)-equivariant conical hyperbolic affine sphere, say \(f_{\operatorname{dev}}\colon\widetilde X\to\R^3\); see Definition \ref{def:singular affine sphere}. We shall denote by \(K_{G_B}\) and \(q_{\mathrm{Pick}}\) the Gaussian curvature of the conical Blaschke metric \(G_B\) and the Pick differential of the associated conical hyperbolic affine sphere, respectively. In particular, they are both uniquely determined by \([L,\gamma,\nu]\) according to Proposition \ref{thm: developing map lifts to a singular affine spheres}. By combining Equation \eqref{eq: Blaschke metric e Pick differential developing map} with the analysis carried out in Theorem \ref{thm: gaussian curvatura singular metrics} and Corollary \ref{cor: boundendess cubic differential}, we directly obtain the following interesting result.

\begin{cor}\label{cor: Curvature Blaschke metric singular case}
The Blaschke metric and the Pick differential of
$f_{\operatorname{dev}}$ satisfy the following properties:
\begin{itemize}
    \item[(i)] If $d\in(0,2g-2)$ and $\operatorname{div}(\gamma)<\operatorname{div}(\nu)$, then $K_{G_B}$ extends to a continuous strictly negative function on $X$. Moreover, $\|q_{Pick}\|_{G_B}^2$ extends to a smooth function on $X$ satisfying \(\|q_{Pick}\|_{G_B}^2 < \frac12\);
    \item[(ii)] If $d\in[1-g,0)$ and \(\operatorname{div}(\nu)<\operatorname{div}(\gamma),\) then $K_{G_B}$ is strictly positive on \(X_0,\)
and tends to $+\infty$ at every point of the possibly empty finite set
\(Z(\gamma)\setminus Z(\nu)\) or a point \(p\in Z(\nu)\) such that \(\operatorname{ord}_p(\nu)<\operatorname{ord}_p(\gamma)\). Moreover, \(\|q_{Pick}\|_{G_B}^2>\frac12 \) on $X_0$, and \(\|q_{Pick}\|_{G_B}^2\) goes to \(+\infty\) approaching either the (possibly empty) finite subset \(Z(\gamma)\setminus Z(\nu)\) or a point \(p\in Z(\nu)\) such that \(\operatorname{ord}_p(\nu)<\operatorname{ord}_p(\gamma)\).
\end{itemize}
\end{cor}

\begin{rmk}
    Although seemingly technical, the result stated in Corollary \ref{cor: Curvature Blaschke metric singular case} highlights an important distinction. The singular hyperbolic affine spheres arising from the construction in the case $d>0$ exhibit a behaviour analogous to that of the Hitchin case $d=2g-2$. By contrast, when $d<0$, a genuinely new phenomenon appears, both in terms of the sign of the curvature of the Blaschke metric and of the boundedness properties of the Pick differential. This difference between the boundedness of \(q_{Pick}\) when \(d>0\) and its possible divergence when \(d<0\) has the following explanation. In the former case, the Pick differential satisfies the additional condition \(\operatorname{ord}_{p_i}(q_{Pick})\ge 3m_i,\) which ensures boundedness, see \eqref{eq: explanation}. By contrast, in the case \(d<0\), the condition \(3m_i>\operatorname{ord}_{p_i}(q_{Pick})\ge 2m_i,\) holds for at least one point \(p_i\), which allows the norm of the Pick differential to diverge. In fact, this can be explicitly seen with a straightforward computation in local coordinates centred in $p_i$. By writing \( q_{Pick}=Q(z)\,dz^3=z^{l_i}Q_0(z)\,dz^3\) with \(Q_0(0)\neq 0, \) we obtain 
    \begin{equation}\label{eq: explanation}
        \|q_{Pick}\|_{G_B}^2 = |Q(z)|^2|z|^{-6m_i}e^{-6u_i} = |Q_0(z)|^2|z|^{2l_i-6m_i}e^{-6u_i}, 
    \end{equation}
    where $u_i$ is smooth near $p_i$. It readily follows that $\|q_{Pick}\|_G^2$ is bounded in a neighbourhood of $p_i$ if and only if \( l_i\geq 3m_i. \)
\end{rmk}

\subsection{Branched hyperbolic structures}\label{sec:extension of branched projective}

In the present section we extend our discussion to branched hyperbolic structures. These geometric structures on surfaces are determined by singular hyperbolic metrics with conical singularities obtained by pulling back the constant negative curvature metric on \(\mathbb H^2\) through an associated branched harmonic map and have been extensively studied over the past fifteen years by several authors, see \cite{derointholozan2026holonomies, faraco2021, faracomaret2025, mathews2011, mathews2012}. In \cite{biswas2021uniformization}, \(\mathrm{SL}(2,\mathbb R)\)-Higgs bundles has been recently used to construct such structures on a Riemann surface \(X\) of genus \(g\ge 2\) from the perspective of non-abelian Hodge theory. In general, the resulting singular metric does not belong to the same conformal class as \(X\). 

\smallskip

By contrast, our construction of branched projective structures of hyperbolic affine sphere type, as in \S \ref{sec: from cyclic to geometric structres}, bypasses harmonic maps into symmetric spaces, relying instead on the developing map into \(\mathbb{RP}^2\). The aim of this section is to explain how the points in the slices \(\mathcal{R}_d^+(X)\), as \(d\) varies in \((0,2g-2)\), partially generalise the branched hyperbolic structures studied in \cite{biswas2021uniformization}. For this purpose, we begin by recalling that every representations \(\rho\in\mathcal{R}_d^+(X)\) corresponds to isomorphism classes of cyclic stable Higgs bundles, say \([L,\gamma,\nu]\) with \(d\in(0,2g-2)\). Let \(\Phi_0\) denote the Higgs field, as in \eqref{eq: higgs data}, obtained by setting \(\nu=0\), and consider the bilinear form
\begin{equation}
    B\defeq \begin{pmatrix}0 & 0 & 1\\0 & -1 & 0\\1 & 0 & 0\end{pmatrix}.
\end{equation}

A direct computation shows that the identity \(B\Phi_0+\Phi_0^tB=0\) holds. In the same fashion, it can be seen that the same condition is satisfied by \((\Phi_0)^{*H}\). It follows that the flat connection \(\nabla_0\) defined as \begin{equation}
    \nabla_0=\nabla_0^H+\Phi_0+(\Phi_0)^{*H}
\end{equation} 
on \(E_{\mathbb R}\) preserves the bilinear form \(B\). As a consequence, the inclusions \(\operatorname{Hol}(\nabla_0)\subset\mathrm{SO}(B|_{E_{\mathbb R}})\subset\mathrm{SL}(3,\mathbb R)\) holds. It is straightforward computation show that \(B\) has signature \((2,1)\) once restricted to the real sub-bundle \(E_{\mathbb R}\). In fact, the real structure is induced by the involution \(\tau(v_1,v_2,v_3)=\bigl(h^{-1}\overline{v}_3,\overline{v}_2,h\overline{v}_1\bigr),\) so every element of \(E_{\mathbb R}\) can be written in the form \(t\defeq \bigl(h^{-1}\overline{v},\lambda,v\bigr),\) with \(v\in\mathbb C\) and \(\lambda\in\mathbb R\) and a direct computation gives \(B(t,t)=2h^{-1}|v|^2-\lambda^2\). Since \(h^{-1}\) is represented by a positive real-valued function, by writing \(v=x+iy\), it readily follows that \(B(t,t)=2h^{-1}(x^2+y^2)-\lambda^2,\)
and hence \(B\) has the signature \((2,1)\) on \(E_{\mathbb R}\) as stated. As a consequence, we may conclude that a Higgs bundle of the form \([L,\gamma,\nu]\) with \(d\in(0,2g-2)\) admits a reduction to an \(\mathrm{SO}_0(2,1)\)-Higgs bundle provided \(\nu=0\). 

In \S\ref{sec: from cyclic to geometric structres} we introduced the smooth section \(s\colon X\longrightarrow E_{\mathbb R}\) given by \(s(p)=(0,1,0)\), which enabled the construction of a branched projective structure of hyperbolic affine sphere type on \(X\). More precisely, this yielded a developing-holonomy pair \((\operatorname{dev}_s,\rho)\), see Proposition \ref{prop:BPSasmaps}. We next observe that \(s\) pointwise defines a time-like vector in each fibre of \((E_{\mathbb R},B)\), as \(B(s,s)=-1\). Since the Klein--Beltrami model of the hyperbolic plane is realised as \(\mathbb H^2=\{[t]\in\mathbb P(E_{\mathbb R})\mid B(t,t)<0\}\), and the developing map is defined by \(\operatorname{dev}_s(p)=[F(p,s(p))]\) via the trivialisation \(F\colon\widetilde{E_{\mathbb R}}\to\widetilde X\times\mathbb R^{2,1}\), it readily follows that \(\operatorname{dev}_s(\widetilde X)\subset\mathbb H^2\subset\mathbb{RP}^2\). In other words, we have proved the following result.

\begin{thm}
    Let \([L,\gamma,0]\) be a cyclic Higgs bundle with \(d\in(0,2g-2)\), and let \((\operatorname{dev}_s,\rho)\) be the developing-holonomy pair associated with the section \(s(p)=(0,1,0)\). Then \(\rho\) takes values in \(\operatorname{SO}_0(2,1)\) and the developing map \(\operatorname{dev}_s\) maps into the Klein--Beltrami model \(\mathbb H^2\subset\mathbb{RP}^2\). Consequently, the pair \((\operatorname{dev}_s,\rho)\) defines a branched hyperbolic structure on \(X\) conformal to the underlying complex structure.
\end{thm}

The same result applies when \(d=2g-2\), in which case the convex projective structure constructed by Baraglia \cite[\S 3.4.2]{baraglia2009g2} reduces to a complete hyperbolic structure conformal to the fixed Riemann surface.

\smallskip

This description provides a particularly transparent description of the construction developed by Biswas--Bradlow--Dumitrescu--Heller in \cite{biswas2021uniformization}. Let \(D=\operatorname{div}(\gamma),\) so that $L\simeq K(-D)$, see \S \ref{sec: prescribed data problem}. As we have done several times, we shall distinguish two scenarios depending on the value of \(d\). 

\subsubsection{Positive values of \(d\)} In the first place, suppose that $d=\deg(L)\in (0,2g-2)$. We also distinguish two sub-cases according to the parity of \(d\). Assume \(d\) positive and even, and choose a square root, say $M$, of $L$. Then \(M^2\simeq K(-D),\) and the $\mathrm{SL}(2,\mathbb R)$-Higgs bundle considered in \cite{biswas2021uniformization}, with vanishing upper-right entry, is given by
\begin{equation}\label{eq: chosen higgs bundle}
    \left(M\oplus M^{-1},\begin{pmatrix}0&0\\\delta&0\end{pmatrix}\right),
\end{equation}

where $\delta$ is the canonical section of \(M^{-2}K\simeq\mathcal O_X(D)\) with divisor $D$. Upon identifying \(\delta\) with \(\gamma\), the Higgs bundle \((\mathcal{E},\Phi_0)\) is, up to a  normalisation of the Higgs field, the image of this rank-two Higgs bundle under the irreducible representation
\begin{equation}
    \operatorname{Sym}_2:\mathbb{P}\mathrm{SL}(2,\mathbb R)\longrightarrow\mathrm{SL}(3,\mathbb R).
\end{equation}

For the vanishing quadratic differential chosen above, the construction of \cite[\S 3]{biswas2021uniformization} yields the hyperbolic metric in the conformal class of \(X\), with conical singularities along \(D\) whose orders are determined by the multiplicities of the zeros of \(\gamma\). Thus, the branched hyperbolic structure arising from the rank-three Higgs bundle as in Equation \eqref{eq: chosen higgs bundle} is the projective realisation, through \(\mathbb P\mathrm{SL}(2,\mathbb R)\cong \mathrm{SO}_0(2,1)\subset \mathrm{SL}(3,\mathbb R)\), of the branched hyperbolic structure arising from the rank-two construction of \cite{biswas2021uniformization}. 

For odd values of $d$ a similar interpretation holds. According to \cite[\S 4]{biswas2021uniformization}, upon passing to an unramified double cover, it is always possible to choose the required square root and then perform the $\mathrm{SL}(2,\mathbb R)$ construction done for even values of \(d\). The induced projective flat connection descends to the original surface $X$, yielding a flat $\mathbb{P}\mathrm{SL}(2,\mathbb R)$-bundle. Since the irreducible three-dimensional representation factors through $\mathbb P\mathrm{SL}(2,\mathbb R)$, the $\mathrm{SO}_0(2,1)$-valued flat connection on $X$ is defined without choosing either the double cover or an $\mathrm{SL}(2,\mathbb R)$ lift. In this sense, the rank-three description provides a uniform formulation of the construction of \cite{biswas2021uniformization} which does not distinguish between the even and odd degree cases. 

\subsubsection{Negative values of \(d\)} The scenario changes noticeably when we consider negative values of the degree, namely \(d\in[1-g,0)\). More precisely, let
\begin{equation}
    (E,\Phi)=\left(L\oplus\mathcal O_X\oplus L^{-1},\begin{pmatrix}0 & 0 & \nu\\\gamma & 0 & 0\\0 & \gamma & 0\end{pmatrix}\right),\qquad 1-g\le\deg(L)<0,\quad \nu\neq 0.
\end{equation}

For \(\gamma\neq 0\), the Higgs bundle is stable, see Proposition \ref{prop: stable SL(3;C) Higgs bundle}, and the projective structure obtained from the section \(s=(0,1,0)\) is branched precisely along \(\operatorname{div}(\gamma)\). On the other hand, for \(\gamma\equiv0\), the Higgs bundle becomes strictly poly-stable. In order to compare the resulting rank-two factor with the conventions used in \cite{biswas2021uniformization}, set \(M\defeq L^{-1},\) so that \(\deg(M)=-\deg(L)>0\), and exchange the two line bundles in the rank-two summand. Then, the strictly poly-stable Higgs bundle can be written as
\begin{equation}
    (\mathcal E,\Phi_0)=(\mathcal O_X,0)\oplus\left(M\oplus M^{-1},\begin{pmatrix}
 0& 0\\\nu & 0\end{pmatrix}\right),
\end{equation}

where \(\nu\in H^0(X,M^{-2}K)\). Thus, the rank-two factor is a stable \(\mathrm{SL}(2,\mathbb R)\)-Higgs bundle (since \(\nu\neq 0\)) of the form considered in \cite{biswas2021uniformization}. Let \(\rho\) be the associated \(\mathrm{SL}(3,\mathbb R)\)-representation. Then \(\rho\) splits as \(\rho=\rho_0\oplus 1\), being \(\rho_0\) the holonomy of a branched hyperbolic structure conformal to \(X\), whose branching divisor is \(\operatorname{div}(\nu)\). In fact, when \(\gamma=0\) the section \(s=(0,1,0)\) becomes parallel for the rank-three flat connection, so that the corresponding developing map degenerates to a constant map. As a consequence, the projective structures obtained for \(\gamma\neq0\), which are branched along \(\operatorname{div}(\gamma)\), do not extend to the case \(\gamma=0\) to a projective structure of the same type. Instead, the reducible limit carries the branched hyperbolic structure arising from the rank-two factor, whose branching locus is determined by the zeros of \(\nu\). In particular, according to Remark \ref{rem: curvature of induced metric from Higgs bundles}, when \(\gamma=0\), the singular negatively curved Riemannian metric \(h^2|\nu|^2\) induced on the surface is precisely the metric with conical singularities obtained in \cite{biswas2021uniformization}.

\section{Anosov representations into \(\mathrm{SL}(3,\R)\)}\label{sec: Anosov}

We finally aim to show Theorem \ref{thm: Anosov Barbot Introduction}, \textit{i.e.}, we aim to show that, if \(\rho\colon\pi_1(S)\longrightarrow\mathrm{SL}(3,\mathbb R)\) is Anosov and non-Hitchin, then necessarily its second Stiefel-Whitney number satisfies the identity \(w_2(\rho)\equiv g-1\pmod 2\). In particular, the Anosov locus in the \(\mathrm{SL}(3,\mathbb R)\)-character variety of a surface group is contained in the union of the Hitchin component and the Barbot component. As a direct consequence, we shall derive that all representations associated with Higgs bundles of the form \([L,\gamma,\nu]\) with \(d=\deg(L)\not\equiv g-1\pmod 2\) are not Anosov. The question remains open in the complementary case.

\subsection{Limit curve and hyperconvex representations}\label{sec: Barbot limit curve}

The proof of our main result relies crucially on two main building blocks, the first of which is the work of Barbot \cite{barbot2010three} whose relevant notions we now briefly recall. We begin with the following

\begin{rmk}\label{rmk:boundary at infinity}
    For a closed connected oriented surface \(S\) of genus \(g\ge 2\), let \(T^1S\) denote its unit tangent bundle. Since \(\pi_1(S)\) is a hyperbolic group, it has a natural boundary \(\partial \pi_1(S)\), which can be identified \(\pi_1(S)\)-equivariantly with \(\partial_\infty \widetilde S\) via the orbit map from \(\pi_1(S)\) to \(\widetilde S\). Upon fixing a hyperbolic metric on \(S\), the boundary at infinity \(\partial_\infty\widetilde S\) can in turn be identified with \(T^1_{\tilde x}\widetilde S\) as follows. Given \(v\in T^1_{\tilde x}\widetilde S\), let \(c_v\) denote the unique geodesic at \(\tilde x\) tangent to \(v\). Then the geodesic ray \(c_v|_{[0,+\infty)}\) converges to a unique point \(v^+\in\partial_\infty\widetilde S\) and the map \(v\mapsto v^+\) is a homeomorphism between \(T^1_{\tilde x}\widetilde S\) and \(\partial_\infty \widetilde S\). Similarly, the ray \(c_{-v}|_{[0,+\infty)}\) determines a second point \(v^-\in\partial_\infty\widetilde S\). We shall use these identifications throughout without further comment, see also \cite[\S\S 1.1 - 1.2]{guichard08}.
\end{rmk}

The theory of Anosov representations has expanded significantly, and several distinct notions are now available in the context of Lie groups. Pioneered by Labourie for $\mathrm{SL}(n,\mathbb R)$ in \cite{Labourie2006Anosov}, the concept was later generalised for more general Lie groups; \textit{e.g.}, see \cite{guichard2012anosov, kapovich2017dynamics, gueritaud2017anosov, bochi2019anosov} and references therein (this list is far from exhaustive). For \(A\in \mathrm{SL}(3,\mathbb R)\), we denote by \(\sigma_1(A)\geq \sigma_2(A)\geq \sigma_3(A)\) its singular values, ordered decreasingly. In the present work we shall adopt the following

\begin{defn}
    A representation \(\rho\colon\pi_1(S)\to \slr3\) is called \emph{Anosov} if there exist positive constants \(C\) and \(C'\) such that \[\log\frac{\sigma_1\big(\rho(\gamma)\big)}{\sigma_2\big(\rho(\gamma)\big)}\ge C|\,\gamma\,|-C',\quad \text{for all} \ \gamma\in\pi_1(S)\] where \(|\gamma|\) denotes the word length of an element with respect to a word metric \(|\,\cdot\,|\) on \(\pi_1(S)\).
\end{defn}

\begin{rmk}
    The definition of Anosov representation is independent of the word metric \(|\,\cdot\,|\). Moreover, since for \(A\in \mathrm{SL}(3,\mathbb R)\) we have \(\sigma_i(A^{-1})=\sigma_{4-i}(\,A\,)^{-1}\) for \(i\in\{1,2,3\}\), an Anosov representation \(\rho\) equivalently satisfies \[\log\frac{\sigma_2(\rho(\gamma))}{\sigma_3(\rho(\gamma))}\geq C|\,\gamma\,|-C'.\] For this reason, such representations are also referred to as \emph{Borel--Anosov}.
\end{rmk}

Among several features enjoyed by Anosov representations, the most relevant for our purposes is the following: for an Anosov representation \(\rho\) then there exist two continuous, injective, \(\rho\)-equivariant maps, denoted by \(\xi_1\colon\partial\pi_1(S)\to\mathbb{RP}^2\) and \(\xi_1^*\colon\partial\pi_1(S)\to(\mathbb{RP}^2)^*\), such that
\begin{itemize}
    \item \(\xi_1(t)\subset\xi_1^*(t)\) for every \(t\in\partial\pi_1(S)\);
    \item \(\xi_1(t)\oplus\xi_1^*(t')=\mathbb{R}^3\) for every \(t\neq t'\).
\end{itemize}

Under the conditions above, by a slight abuse of notation, we may regard \(\xi_1(t)\) as a line in \(\mathbb{R}^3\), and \(\xi_1^*(t)\) as a plane in \(\mathbb{R}^3\). Another important notion, introduced by Labourie \cite{Labourie2006Anosov} in arbitrary dimension, is that of a \emph{hyperconvex} representation. In the case of \(\mathrm{SL}(3,\mathbb R)\), an Anosov representation \(\rho\) is hyperconvex if and only if it preserves a strictly convex domain in \(\mathbb{RP}^2\), and hence if and only if it belongs to the Hitchin component. In what follows, we shall make use of the following characterisation provided by Barbot in \cite[Corollary 6.6]{barbot2010three}. 

\begin{prop}[Barbot]\label{prop: limit curve}
    An Anosov representation \(\rho\colon\pi_1(S)\longrightarrow \mathrm{SL}(3,\mathbb R)\) is hyperconvex if and only if the limit curve \(\xi_1\bigl(\partial\pi_1(S)\bigr)\subset\mathbb{RP}^2\) is null-homotopic. Equivalently, the dual limit curve \(\xi_1^*\bigl(\partial\pi_1(S)\bigr)\subset(\mathbb{RP}^2)^*\) is null-homotopic.
\end{prop}

\subsection{Cohomology of the unit tangent bundle}\label{sec: cohomoogy unit tangent bundle}
The second building block for the proof of Theorem \ref{thm: Anosov Barbot Introduction} is a computation in the cohomology of the unit tangent bundle. In the first part, we recall and state the results needed from \cite[Appendix B.1]{guichard2010}. In the second part, we prove a fundamental proposition, obtained as a consequence of the work of Guichard--Wienhard, which plays a crucial role in the proof of the main theorem.

\smallskip

Let $\pi\colon T^1S\to S$ denote the natural projection, and let \(\eta\in H^2(S;\mathbb Z)\) be the orientation class. We denote by $\overline\eta\in H^2(S;\mathbb Z_2)$ its reduction modulo $2$. In what follows, we shall denote by $\langle\,\cdot,\cdot\,\rangle$ the natural pairing between cohomology and homology, without specifying the coefficients, as these will always be clear from the context. By construction, we know that $\langle\eta,[S]\rangle=1$ and $\langle\overline{\eta},[S]_2\rangle=1$, where $[S]$ and $[S]_2$ denote the fundamental classes in $H_2(S,\mathbb Z)$ and $H_2(S,\mathbb Z_2)$, respectively. For $p\in S$, let $j_p\colon T^1_pS\hookrightarrow T^1S$ be the inclusion of the fibre. Given $\alpha\in H^1(T^1S;\mathbb Z_2)$, we set
\begin{equation}\label{eq:fibre-evaluation}
    \varepsilon(\alpha)\defeq \big\langle j_p^*\alpha,[T^1_pS]_2\big\rangle\in\mathbb Z_2,
\end{equation}

where $[T^1_pS]_2$ is the mod-$2$ fundamental class of the circle. This definition is clearly independent of $p$, since the fibres over any two points are homologous.

\begin{lem}[{\cite[Proposition B.1]{guichard2010}}]\label{lem:cohomology-unit-tangent}
    The following statements hold.
    \begin{enumerate}
        \item[(i)] The map \(\pi^*\colon H^1(S;\mathbb Z)\to H^1(T^1S;\mathbb Z)\) is an isomorphism. Moreover,
        \begin{equation}
            H^2(T^1S;\mathbb Z)\cong\mathbb Z^{2g}\oplus\mathbb Z/(2g-2)\mathbb Z,
        \end{equation}
        and its torsion subgroup is generated by $\tau\defeq \pi^*\eta$.
        \item[(ii)] With $\mathbb Z_2$ coefficients, there are short exact sequences
        \begin{equation} 
        \label{eq:gysin-mod-two-degree-one}
            0\to H^1(S;\mathbb Z_2) \xrightarrow{\pi^*}H^1(T^1S;\mathbb Z_2) \xrightarrow{\varepsilon}\mathbb Z_2 \to 0
        \end{equation}
    and
    \begin{equation} \label{eq:gysin-mod-two-degree-two}
        0\to H^2(S;\mathbb Z_2)\xrightarrow{\pi^*}H^2(T^1S;\mathbb Z_2)\xrightarrow{}H^1(S;\mathbb Z_2)\to0 .
    \end{equation}
    In particular, $\pi^*\colon H^2(S;\mathbb Z_2)\to H^2(T^1S;\mathbb Z_2)$ is injective, and
    \begin{equation}
        H^1(T^1S;\mathbb Z_2)\cong(\mathbb Z_2)^{2g+1},\quad H^2(T^1S;\mathbb Z_2)\cong(\mathbb Z_2)^{2g+1}.
    \end{equation}
    \end{enumerate}
\end{lem}

With a slight abuse of notation, we denote the pull-back map induced in cohomology again by $\pi^*$, since the coefficients will always be clear from the context. The two short exact sequences, as well as the proof of the previous lemma, are obtained by applying the Gysin sequence to the $S^1$-bundle $T^1S\to S$; see, for instance, \cite[Proposition 14.33]{bott1982differential}. Let $\smile$ denote the cup product on $H^*(T^1S,\mathbb Z_2)$ and recall that, for a generic element $\alpha\in H^1(T^1S,\mathbb Z_2)$, then $\alpha\smile\alpha$ is a well-defined element of $H^2(T^1S,\mathbb Z_2)$, see \cite[Chapter 3 - \S 3.2]{Hatcher2002AlgebraicTopology}. The following result is a consequence of Lemma \ref{lem:cohomology-unit-tangent} and standard arguments in algebraic topology.

\begin{prop}\label{prop: cup product itself}
    For every $\alpha\in H^1(T^1S;\mathbb Z_2)$, one has
        \begin{equation}\label{eq: cup product square}
            \alpha\smile\alpha=(g-1)\varepsilon(\alpha)\,\pi^*\overline\eta\quad\text{in }H^2(T^1S;\mathbb Z_2).
        \end{equation}
    In particular, if the restriction of $\alpha$ to a fibre is non-zero, then
    \begin{equation}
        \alpha\smile\alpha=(g-1)\pi^*\overline\eta.
    \end{equation} 
\end{prop}
\begin{proof}

We denote by $r_S$ and $r_{T^1S}$ the reduction maps modulo $2$ on the cohomology of $S$ and $T^1S$, respectively. We also denote by $\beta$ and $\beta_2$ the Bockstein homomorphisms associated with the short exact sequences 
\begin{equation}
    0\longrightarrow\mathbb Z\xrightarrow{\times2}\mathbb Z\longrightarrow\mathbb Z_2 \longrightarrow0\quad\text{ and } \quad 0\longrightarrow\mathbb Z_2 \xrightarrow{\,1\mapsto2\,}\mathbb Z_4 \longrightarrow\mathbb Z_2\longrightarrow0,
\end{equation}

respectively; see \cite[\S 3.E]{Hatcher2002AlgebraicTopology} for more details. The first chain of equalities that we want to prove is 
\begin{equation}
    \operatorname{Im}(\,r_{T^1S}\,)=\pi^*H^1(S,\mathbb{Z}_2)=\operatorname{ker}(\,\varepsilon\,)\subset H^1(T^1S,\mathbb{Z}_2).
\end{equation}

The second equality follows directly from the short exact sequence \eqref{eq:gysin-mod-two-degree-one}. The first equality, instead, follows from the fact that \(\pi^*\colon H^1(S,\mathbb{Z})\to H^1(T^1S,\mathbb{Z})\) is an isomorphism by Lemma \ref{lem:cohomology-unit-tangent}, point \((i)\), combined with the relation \(r_{T^1S}\bigl(\pi^*a\bigr)=\pi^*\bigl(r_S(a)\bigr)\), for every $a\in H^1(S,\mathbb{Z})$. Recall that the maps $r_{T^1S}$ and $\beta$ are related by the long exact sequence associated with the corresponding short exact sequence of coefficients, see \cite[\S 3.E]{Hatcher2002AlgebraicTopology}. For our purposes, the relevant part is the sequence
\begin{equation}\label{eq:bockstein-unit-tangent}
    H^1(T^1S;\mathbb Z)\xrightarrow{r_{T^1S}}H^1(T^1S;\mathbb Z_2)\xrightarrow{\beta} H^2(T^1S;\mathbb Z)\xrightarrow{\times 2}H^2(T^1S;\mathbb Z).
\end{equation}

In particular, by using the identity $\ker(\varepsilon)=\operatorname{Im}(r_{T^1S})$ and by the exactness of \eqref{eq:bockstein-unit-tangent}, it follows that 
\(\beta(\alpha)=0\) if and only if \(\varepsilon(\alpha)=0\), for every $\alpha\in H^1(T^1S;\mathbb Z_2)$. On the other hand, $\operatorname{Im}\beta$ is the subgroup of $H^2(T^1S;\mathbb Z)$ annihilated by multiplication by $2$. By Lemma \ref{lem:cohomology-unit-tangent} point \((i)\), this subgroup is \(\{0,(g-1)\tau\}\), where \(\tau=\pi^*\eta.\) Consequently,
\begin{equation}\label{eq:integral-bockstein-fibre}
    \beta(\alpha)=(g-1)\varepsilon(\alpha)\tau,
\end{equation} where on the right we regard $\varepsilon(\alpha)$ as either $0$ or $1$. By taking the reduction modulo $2$ of \eqref{eq:integral-bockstein-fibre} and by using the identities \((r_{T^1S}\circ\beta)(\alpha)=\beta_2(\alpha)=\alpha\smile\alpha\), when $\alpha\in H^1(T^1S;\mathbb Z_2)$; see \cite[p. 489 - point (7)]{Hatcher2002AlgebraicTopology}, we directly obtain the desired relation \eqref{eq: cup product square}.
\end{proof}

\subsection{Main result and applications}\label{sec: main result Anosov}

We briefly recall that, given a representation $\rho\colon \pi_1(S)\to \mathrm{SL}(3,\R)$, we denote by $E_\rho$ the associated real flat bundle and by $w_2(E_\rho)$ its second Stiefel--Whitney class; see \S \ref{sec: representation space}. Using the notation introduced in \S \ref{sec: cohomoogy unit tangent bundle}, we define the \textit{second Stiefel--Whitney number} of $\rho$ by
\begin{equation}
    w_2(\,\rho\,)\defeq \langle w_2(E_\rho),[S]_2\rangle.
\end{equation}

Recall also that the Barbot component is identified with $\mathcal C_0$ when the genus of $S$ is odd, and with $\mathcal C_1$ when the genus of $S$ is even; see \S \ref{sec: Barbot component}. The main result of this section is the following.
\begin{thm}\label{thm: Anosov}
    Let \(\rho\colon\pi_1(S)\to \slr3\) be an Anosov and non-Hitchin representation. Then, 
    \begin{equation}
        w_2(\,\rho\,)\equiv g-1 \ (\text{mod} \ 2).
    \end{equation} 
    In particular, the set of Anosov representations is contained in the union of the Hitchin component and the Barbot component.
\end{thm}

Recall that, after fixing a Riemann surface structure \(X\) on \(S\), in \S\ref{sec: from cyclic to geometric structres} we constructed submanifolds \(\mathcal{R}_d^+(X)\) and \(\mathcal{R}_d^-(X)\) of the character variety, whose points arise as holonomies of branched projective structures of hyperbolic affine sphere type. The following result is an immediate consequence of the theorem above and Lemma \ref{lem: rho_connected_component_Higgs_bundles}.

\begin{cor}\label{cor:non-anosov-cyclic-slices}
    Let $X$ be a Riemann surface structure on $S$, and suppose that \begin{equation} 
        d\in(1-g,2g-2)\setminus\{0\},\quad d\not\equiv g-1\pmod2.
    \end{equation}
    Then no representation in $\mathcal R_d^+(X)$ and $\mathcal{R}_d^-(X)$ is Anosov.
\end{cor}

We now focus on the proof of Theorem \ref{thm: Anosov}, relying on the material introduced in \S \ref{sec: Barbot limit curve}.  Although the same construction can also be carried out in the Hitchin case, we restrict to the case of Anosov and non-Hitchin representations \(\rho\colon\pi_1(S)\to \mathrm{SL}(3,\mathbb R)\), since this is the case of major interest for our purposes. Upon fixing a conformal identification \(\widetilde S\cong \mathbb H\), then \(\partial\pi_1(S)\cong\partial_\infty\widetilde S\cong T^1_{\widetilde x}\widetilde{S}\cong\mathbb{RP}^1\) hold by Remark \ref{rmk:boundary at infinity}.  

Let \(\xi_1\) and \(\xi_2\) denote the two boundary maps associated to \(\rho\), where \(\xi_2\) associates to each \(t\in\partial\pi_1(S)\) the \(2\)-plane in \(\mathbb R^3\) dual to the point \(\xi_1^*(t)\in(\mathbb{RP}^2)^*\). We then introduce three real line bundles \(\widetilde{\mathcal L}_+\), \(\widetilde{\mathcal L}_0\), and \(\widetilde{\mathcal L}_-\) over the unit tangent bundle \(T^1\widetilde S\) as follows: the fibres over a point \((p,v)\in T^1\widetilde S\) are defined by
\begin{equation}
\widetilde{\mathcal L}_+(p,v)\defeq \xi_1(v^+),\qquad
\widetilde{\mathcal L}_0(p,v)\defeq \xi_2(v^+)\cap\xi_2(v^-),\qquad
\widetilde{\mathcal L}_-(p,v)\defeq \xi_1(v^-),
\end{equation}

where \(v^\pm\) denote the endpoints in \(\partial_\infty\widetilde S\) of the unique geodesic through \(p\) tangent to \(v\). We may also observe that \(\widetilde{\mathcal L}_0\) is well defined, since the planes \(\xi_2(v^+)\) and \(\xi_2(v^-)\) are distinct by the properties of the boundary maps recalled in \S \ref{sec: Barbot limit curve}. It is straightforward to check that the three line bundles introduced above are in direct sum. We first observe that the two lines \(\widetilde{\mathcal L}_+(p,v)\) and \(\widetilde{\mathcal L}_0(p,v)\) are distinct. Indeed, if they coincided, we would have \(\xi_1(v^+)\subset\xi_2(v^-)\), contradicting the transversality property of the boundary maps. It follows that \(\xi_2(v^+)=\widetilde{\mathcal L}_+(p,v)\oplus\widetilde{\mathcal L}_0(p,v)\). Moreover, by transversality, \(\xi_1(v^-)\not\subset\xi_2(v^+)\), and hence we obtain the direct sum decomposition
\begin{equation}
    \mathbb R^3=\widetilde{\mathcal L}_+(p,v)\oplus\widetilde{\mathcal L}_0(p,v)\oplus\widetilde{\mathcal L}_-(p,v).
\end{equation} 
By the \(\rho\)-equivariance of the boundary maps, the construction descends to the quotient \(T^1S\), yielding a real rank-\(3\) vector bundle
\begin{equation}
    \overline E_\rho\defeq \mathcal L_+\oplus\mathcal L_0\oplus\mathcal L_-\longrightarrow T^1S.
\end{equation}
Recall that, starting from \(\rho\), one can construct the flat bundle \(E_\rho=\big(\widetilde S\times\mathbb R^3\big)/_{\sim_\rho}\) over \(S\). If \(\pi\colon T^1S\to S\) denotes the projection from the unit tangent bundle, then the following result holds.

\begin{lem}\label{lem: first Stiefel-Whitney equal}
    The vector bundle \(\overline E_\rho\) is isomorphic to \(\pi^*E_\rho\). As a consequence, both \(\overline E_\rho\) and \(\mathcal{L}_0\) are orientable and 
    \begin{equation}
        w_1(\mathcal{L}_+)=w_1(\mathcal{L}_-).
    \end{equation}
\end{lem}

\begin{proof}
    Let $q\colon T^1\widetilde S\to T^1S$ be the covering induced by $\widetilde S\to S$, and let $\widetilde\pi\colon T^1\widetilde S\to\widetilde S$ be the natural lifted projection. There is a canonical identification 
    \begin{equation}
        T^1\widetilde S\times\mathbb R^3/\sim_\rho\longrightarrow\pi^*E_\rho,\qquad [v,u]\longmapsto \bigl(q(v),[\widetilde\pi(v),u]\bigr).
    \end{equation}
    Under this identification, the required isomorphism is
        \begin{equation}\label{eq:explicit-splitting-isomorphism}
            \begin{aligned}
                \Phi\colon\overline E_\rho&\longrightarrow\pi^*E_\rho,\\\bigl([v,u_+],[v,u_0],[v,u_-]\bigr)&\longmapsto\bigl(q(v),[\widetilde\pi(v),u_++u_0+u_-]\bigr),
            \end{aligned}
        \end{equation}
    where the three representatives are taken over a common lift $v$. The map is independent of this choice because vector addition commutes with $\rho(\gamma)$ for every $\gamma\in\Gamma$. Moreover, since \(E_\rho\) is orientable, naturality of the first Stiefel--Whitney class yields the equalities \(w_1(\overline E_\rho)=\pi^*w_1(E_\rho)=0\). In order to study the orientability of $\mathcal L_0$, let $\alpha\defeq w_1(\mathcal L_+)$ and consider the map $\iota\colon T^1S\to T^1S$ acting fibre-wise as the antipodal map. Since $S$ is an oriented surface, the Riemannian metric determines a natural fibre-wise rotation action of $S^1$ on $T^1S$. In particular, the mapping $\iota$ is the rotation by $\pi$, and is therefore homotopic to the identity through fibre-wise rotations. Consequently, $\iota^*$ acts as the identity on $H^1(T^1S;\mathbb Z_2)$. Now choose a lift \(\widetilde p\in\widetilde S\) of \(p\). Under the identification \(T_{\widetilde p}^1\widetilde S\to\partial_\infty\widetilde S\) sending \(v\) to \(v^+\), we have \((\iota(v))^+=(-v)^+=v^-\). Indeed, if \(c_v\) denotes the geodesic associated with \(v\), then the geodesic associated with \(-v\) satisfies \(c_{-v}(t)=c_v(-t)\). As a consequence, we obtain \(\iota^*\widetilde{\mathcal L}_+=\widetilde{\mathcal L}_-\), and, by \(\rho\)-equivariance, also \(\iota^*\mathcal L_+=\mathcal L_-\), with a slight abuse of notation. In particular,
    \begin{equation}
        w_1(\mathcal L_-)=\iota^*w_1(\mathcal L_+)=w_1(\mathcal L_+)=\alpha,
    \end{equation}
    since \(\iota\) induces the identity map in cohomology. This implies directly that \(w_1(\mathcal{L}_0)=0\).
\end{proof}

\begin{rmk}
    Although our construction is entirely explicit, it actually arises as a special case of a much more general framework developed by Guichard--Wienhard, see \cite[\S2.1.4]{guichard2010}. Indeed, by considering the two-dimensional subgroup of diagonal matrices \(Z<\mathrm{SL}(3,\mathbb R)\), we obtain a natural \(\rho\)-equivariant map, say \(\sigma\colon T^1\widetilde S\to \mathrm{SL}(3,\mathbb R)/Z\) mapping \(v\) to \((\widetilde{\mathcal L}_+(v),\widetilde{\mathcal L}_0(v),\widetilde{\mathcal L}_-(v))\). This map encodes a \(Z\)-reduction of the associated flat bundle, causing it to split as the direct sum of three real line bundles; cf. the frame variety in \cite[\S 2.3]{barbot2010three}.
\end{rmk}

Having established the first basic properties of the three line bundles $\mathcal L_+$, $\mathcal L_0$, and $\mathcal L_-$, the second step is to relate the fact that the representation $\rho$ is non-Hitchin to the non-orientability of $\mathcal L_+$ when restricted to a fibre of the unit tangent bundle of the surface. As will become clear in the proof, this relies crucially on Barbot's characterisation of the limit curve given in Proposition \ref{prop: limit curve}.
\begin{lem}\label{lem: non orientability on fibres}
If $\rho$ is Anosov and non-Hitchin, then for any $p\in S$ the real line bundle $\mathcal{L}_+|_{T^1_pS}$ is non-orientable.
\end{lem}

\begin{proof}
    The proof is a direct consequence of the definition of $\mathcal L_+$ in terms of the boundary map $\xi_1$. Indeed, let $p\in S$ and let $\widetilde p\in\widetilde S$ be a lift of $p$. We denote by $i_{p}$ the induced identification between $T_p^1S$ and $\partial_\infty\pi_1(S)$ described at the beginning of \S \ref{sec: Barbot limit curve}, and define
    \begin{equation} 
        f_p\defeq \xi_1\circ i_{p}\colon T_p^1S\longrightarrow \mathbb{RP}^2.
    \end{equation}
    Consider the tautological line bundle $\mathcal O_{\mathbb{RP}^2}(-1)$ over the real projective plane, whose fibre over a generic point $[\ell]\in\mathbb{RP}^2$ is the line $\ell\subset\mathbb R^3$ itself. It follows directly from the definition that
    \begin{equation}
        \mathcal L_+|_{T_p^1S}\cong f_p^*\mathcal O_{\mathbb{RP}^2}(-1).
    \end{equation}
    By Proposition \ref{prop: limit curve}, the map $f_p\colon T_p^1S\cong S^1\to\mathbb{RP}^2$ represents the non-trivial element of $\pi_1(\mathbb{RP}^2)\cong\mathbb Z_2$. Since the tautological line bundle $\mathcal O_{\mathbb{RP}^2}(-1)$ is non-orientable, its pull-back by $f_p$ is the non-trivial real line bundle over $S^1$. Therefore, $\mathcal L_+|_{T_p^1S}$ is non-orientable.
\end{proof}

We are now in the right position to conclude the proof of Theorem \ref{thm: Anosov} as a combination the previous discussion with the results of \S \ref{sec: cohomoogy unit tangent bundle}.

\begin{proof}[Proof of Theorem \ref{thm: Anosov}]
    Let $\alpha\defeq w_1(\mathcal L_+)\in H^1(T^1S,\Z_2)$. By Lemma \ref{lem: first Stiefel-Whitney equal}, we have $w_1(\mathcal L_-)=\alpha$ and $w_1(\mathcal L_0)=0$. The Whitney sum formula and naturality of the second Stiefel-Whitney class therefore yield
    \begin{equation}\label{eq:second-stiefel-whitney-square}
        \begin{aligned}
            \pi^*w_2(E_\rho)&=w_2(\mathcal L_+\oplus\mathcal L_0\oplus\mathcal L_-) \\ & =w_1(\mathcal L_+)\smile w_1(\mathcal L_0)+w_1(\mathcal L_+)\smile w_1(\mathcal L_-) \\ & \qquad+w_1(\mathcal L_0)\smile w_1(\mathcal L_-) \\ & =\alpha\smile\alpha.
        \end{aligned}
    \end{equation} By Lemma \ref{lem: non orientability on fibres} we know that $\varepsilon(\alpha)=1$, see Equation \eqref{eq:fibre-evaluation}. Therefore, applying Proposition \ref{prop: cup product itself}, we obtain \(\pi^*w_2(E_\rho)=(g-1)\pi^*\overline\eta\). The injectivity of $\pi^*\colon H^2(S;\mathbb Z_2)\to H^2(T^1S;\mathbb Z_2)$, established in Lemma \ref{lem:cohomology-unit-tangent} point $(ii)$, gives
    \begin{equation}
        w_2(E_\rho)=(g-1)\overline\eta. 
    \end{equation}
    By evaluating on the mod-$2$ fundamental class of $S$ proves the claimed congruence.
\end{proof}

\begin{rmk}\label{rem:topological-invariant-anosov}
    The construction above provides a topological invariant that distinguishes Hitchin representations from non-Hitchin Anosov representations, independently of the genus of the surface. Indeed, let $\rho:\pi_1(S)\to\mathrm{SL}(3,\mathbb R)$ be Anosov, and consider the class \(\alpha\defeq w_1(\mathcal L_+)\in H^1(T^1S;\mathbb Z_2)\). By construction, the restriction of $\mathcal L_+$ to any fibre $T_p^1S$ is isomorphic to the pull-back of the tautological line bundle $\mathcal O_{\mathbb{RP}^2}(-1)$ under the map $\xi_1\circ i_{p}:T_p^1S\to\mathbb{RP}^2$. Consequently, by Proposition \ref{prop: limit curve}, we have
    \begin{equation}
        \varepsilon(\alpha)=
        \begin{cases}
            0, & \text{if }\rho\text{ is Hitchin},\\ 1, & \text{if }\rho\text{ is non-Hitchin}.
        \end{cases}
    \end{equation}
Equivalently, $\mathcal L_+|_{T_p^1S}$ is isomorphic to the trivial line bundle in the Hitchin case, and to the Möbius line bundle otherwise. In particular, the restriction of $\alpha$ to the fibres completely characterises the Hitchin locus among Anosov representations. This is stronger than the information provided by $w_2(E_\rho)$: when $g$ is odd, the latter vanishes for both Hitchin and non-Hitchin Anosov representations, whereas $\varepsilon(\alpha)$ still distinguishes the two cases.
\end{rmk}

\bibliographystyle{amsalpha}
\bibliography{rbpshiggs}

\end{document}